\documentclass[12pt]{amsart}
\usepackage{amsbsy,amssymb,amsmath,amsthm,amscd,amsfonts,latexsym,amstext,delarray,
amsmath,graphicx} 
\usepackage{upgreek}
\usepackage{bm}
\usepackage{qtree}
\usepackage{ytableau}
\usepackage{tikz-cd}
\usepackage[margin=.8in]{geometry}
\usepackage{color}
\usepackage[all]{xy}

\newtheorem{prop}{Proposition}[section]
\newtheorem{thm}{Theorem}[section]
\newtheorem{cor}[prop]{Corollary}
\newtheorem{lem}[prop]{Lemma}

\newtheorem{defn}[prop]{Definition}
\newtheorem{rem}[prop]{Remark}

\newtheorem{ques}[prop]{Question}

\usepackage[normalem]{ulem}

\usepackage[]{mdframed}

\usepackage{bigstrut}

\numberwithin{equation}{section}

\def\F{{\mathbb F}}

\def\R{{\mathbb R}}

\def\P{{\mathbb P}}

\def\cA{{\mathcal A}}
\def\cB{{\mathcal B}}
\def\cC{{\mathcal C}}

\def\cF{{\mathcal F}}
\def\cG{{\mathcal G}}
\def\cH{{\mathcal H}}
\def\cI{{\mathcal I}}

\def\cK{{\mathcal K}}
\def\cL{{\mathcal L}}

\def\cO{{\mathcal O}}
\def\cP{{\mathcal P}}
\def\cQ{{\mathcal Q}}
\def\cR{{\mathcal R}}
\def\cS{{\mathcal S}}
\def\cT{{\mathcal T}}

\def\cV{{\mathcal V}}
\def\cW{{\mathcal W}}

\def\bB{{\mathbb B}}

\def\bE{{\mathbb E}}
\def\bF{{\mathbb F}}

\def\bP{{\mathbb P}}

\def\Tr{{\rm Tr}}

\def\fM{{\mathfrak M}}
\def\fT{{\mathfrak T}}
\def\fF{{\mathfrak F}}
\def\fB{{\mathfrak B}}
\def\fc{{\mathfrak c}}
\def\fb{{\mathfrak b}}
\def\fO{{\mathfrak O}}
\def\fp{{\mathfrak p}}

\title[Merge Dynamics]{Merge on workspaces as Hopf algebra Markov chain}
\author[M.Marcolli, D.Skigin]{Matilde Marcolli and David Skigin}
\date{2025}

\address{Department of Mathematics and Department of Computing and Mathematical Sciences, 
California Institute of Technology, CA 91125, USA}
\email{matilde@caltech.edu}
\address{Mathematical \& Physical Sciences, University College London, London, WC1E 6BT, UK} 
\email{david.skigin.24@ucl.ac.uk}

\begin{document}
\maketitle

\begin{abstract}
We study the dynamical properties of a Hopf algebra Markov chain with state space the binary
rooted forests with labelled leaves. This Markovian dynamical system describes 
the core computational process of structure formation and transformation in syntax via
the Merge operation, according to Chomsky's Minimalism model of generative linguistics. 
The dynamics decomposes into an ergodic dynamical system with uniform stationary
distribution, given by the action of Internal Merge,
while the contributions of External Merge and (a minimal form of) Sideward Merge reduce
to a simpler Markov chain with state space the set of partitions and with combinatorial weights.
The Sideward Merge part of the dynamics prevents convergence to fully
formed connected structures (trees), unless the different forms of Merge are weighted by
a cost function, as predicted by linguistic theory. Results on the asymptotic behavior of the 
Perron-Frobenius eigenvalue and eigenvector in this weighted case, obtained in terms of an 
associated Perron-Frobenius problem in the tropical semiring, show that the usual cost
functions (Minimal Search and Resource Restrictions) proposed in the linguistic literature
do not suffice to obtain convergence to the tree structures, while an additional optimization
property based on the Shannon entropy achieves the expected result for the dynamics. 
We also comment on the introduction of continuous parameters related to semantic
embedding and other computational models, and also on some filtering of the dynamics
by coloring rules that model the linguistic filtering by theta roles and phase structure,
and on parametric variation and the process of parameter setting in Externalization.
\end{abstract}

\tableofcontents

\section{Introduction}

In work of Diaconis, Pang, and Ram in \cite{DiPaRa} and of Pang 
in \cite{Pang}, \cite{Pang2}, a theory of Hopf algebra Markov chains was
developed. These are discrete dynamical systems associated to combinatorial Hopf algebras. 
These Hopf algebras (see \cite{LoRo}) are graded connected commutative Hopf algebras $\cH$
with a basis given by combinatorial objects (in our case, these will be binary forests). 
The dynamical systems called Hopf algebra Markov chains act on the vector space of the 
Hopf algebra as a composition of the form
\begin{equation}\label{HMCgen}
 \cK =  \mu \circ \cL \circ \Delta \, , 
\end{equation} 
where $\mu$ is the Hopf algebra product, $\Delta$ is the Hopf algebra coproduct, and
$\cL$ is a linear endomorphism of $\cH\otimes \cH$.  In other words, one decomposes
combinatorial objects using the coproduct $\Delta$, acts with $\cL$ on the parts of the
decomposition, and recomposes a resulting object using the product $\mu$. 
In the main examples studied so far $\cL$ is a projection. More general forms have
also been studied in  \cite{DiPaRa}, \cite{Pang}, \cite{Pang2} where one uses powers
$\Delta^k$ and $\mu^k$ in the composition \eqref{HMCgen}, namely decompositions into
multiple parts. For our purposes, we will only consider cases of the form \eqref{HMCgen},
where $\cL$ is not just a projection but the composition of a projection and a ``grafting"
operation. While in the cases studied in  \cite{DiPaRa}, \cite{Pang}, \cite{Pang2} the
special form of $\cL$ leads to a combinatorial description of the stationary distribution,
the case we consider does not have that structure, but we will still be able to derive some
general results about the stationary distribution and the dynamics, based on the 
Perron--Frobenius theorem and explicit counting of combinatorial objects, involving
partitions and generalized multinomial coefficients.

\smallskip

The form of the Hopf algebra Markov chain we consider is motivated by 
a mathematical model of the Merge operation in linguistics, as developed
in work of the first author, Chomsky, and Berwick in \cite{MCB}.

\smallskip

The Minimalist Program of syntax in generative linguistics was first introduced by Noam Chomsky 
in the early 1990s, \cite{Chomsky95}. Since then it underwent significant developments, culminating 
in the recent formulation, starting around 2013, where Chomsky introduced the ``free symmetric Merge"
as the fundamental structure-building operation of syntax, and isolated this universal core computational
mechanism, subject to optimality constraints, \cite{ChomskyUCLA}, \cite{ChomskyGK}, \cite{ChomskyElements}. 
In this model, the core Merge action accounts for all the structure-formation and key computational aspects
of syntax (compositional and transformational properties), while an ``externalization interface" accounts 
for parametric syntactic variation across languages and 
embodiment into specific languages, and a ``syntax-semantics interface" accounts for semantic parsing. 
In 2023, joint work of Marcolli, Chomsky, and Berwick showed that the core computational mechanism
of the free symmetric Merge admits a precise mathematical formalization in terms of Hopf algebras and
Hopf algebra Markov chain dynamical systems, \cite{MCB}. The work \cite{MCB}, as
well as subsequent work \cite{MLH}, \cite{MHL}, \cite{ML}, \cite{MarBer}, \cite{SM} and further ongoing developments
show that the two interfaces (externalization and syntax-semantics) also admit good mathematical modeling, 
though providing a full mathematical account of their functioning still requires a significant amount of further 
investigation. 

In this paper we focus on the core computational structure of the free symmetric Merge, as formulated in \cite{MCB}, and we
further investigate its mathematical properties as a dynamical system acting on syntactic workspaces. Namely, 
we study the universal dynamical properties of syntax formation and transformation, as a Hopf algebra Markov chain.

A mathematically precise model of structure formation in language is of particular current interest, because current
artificial systems that simulate human language (such as large language models) are extremely wasteful in terms
of both training and computational resources, and are based purely on statistical data (computing probabilities
of extensions of strings of text), while completely ignoring a wealth of very detailed knowledge that we possess
about the structure of language and its underlying computational mechanism. The astonishing theory-gap that
one can observe in recent developments in the field of artificial intelligence makes it a urgent task to formulate
the existing scientific theory of language-production, broadly known as ``generative grammar",  
in a precise mathematical language. While applications are not the focus of this and other theoretical papers, they are
one of the natural motivations for developing our mathematical models of syntax as a computational system.

\subsection{Summary of the Merge action on workspaces}

We recall here some basic results and notation from \cite{MCB} that we will be using throughout the
paper.

\begin{itemize}
\item $\cS\cO_0$ is a finite set. Elements $\alpha\in \cS\cO_0$ are called {\em lexical items and syntactic features}
(we will just call them ``lexical data" here for brevity).
\item $\cS\cO=\fT_{\cS\cO_0}$ is the countably infinite set of full binary rooted
trees (nonplanar, that is, with no assigned planar embedding) with leaves decorated by
elements of the set $\cS\cO_0$. Elements $T\in \cS\cO$ are called {\em syntactic objects}. We write
$\fT_{\cS\cO_0,n}$ for the set of syntactic objects $T\in \cS\cO$ with $n=\# L(T)$ leaves. 
\item $\cS\cO$ is the free nonassociative commutative magma generated by the set $\cS\cO_0$
$$ \cS\cO = {\rm Magma}_{c,na}(\cS\cO_0, \fM)\cong \fT_{\cS\cO_0} \, , $$
with binary operation $\fM$ that acts on a pair of trees by joining them to a common root (as nonplanar trees),
$$ \fM(T,T') :=\Tree[ $T$ $T'$ ]= \Tree[ $T'$ $T$ ] \, . $$
\item $\fF_{\cS\cO_0}$ is the countably infinite set of binary rooted forests, namely forests
whose connected components are nonplanar full binary rooted trees in $\fT_{\cS\cO_0}$. 
Elements $F\in \fF_{\cS\cO_0}$ are called {\em workspaces}. We
also write $\fF_{\cS\cO_0,n}$ for the set of workspaces $F=T_1\sqcup\cdots\sqcup T_r \in \fF_{\cS\cO_0}$ 
with $n=\# L(F)=\# L(T_1) + \cdots + \# L(T_r)$ leaves. 
\item An {\em accessible term} $T_v\subset T$ of a syntactic objects is
the subtree consisting of a non-root vertex $v$ of $T$ and all it descendants. Thus, the set of
accessible terms of $T$ can be identified with the set of non-root vertices of $T$.
\item The vector space $\cV(\fF_{\cS\cO_0})$ endowed with the product
given by the disjoint union $\sqcup$ is the free commutative algebra (polynomial
algebra) generated by $\fT_{\cS\cO_0}$.
\item There are different forms of coproduct on $\cV(\fF_{\cS\cO_0})$ with
slightly different algebraic properties (see the discussion in \S 1.2 of \cite{MCB}). All of them
are of the form
\begin{equation}\label{coprod}
 \Delta(T)=\sum_{\underline{v}} F_{\underline{v}} \otimes T/F_{\underline{v}} \, , 
\end{equation} 
with $F_{\underline{v}}=T_{v_1}\sqcup\cdots \sqcup T_{v_n}$ a collection of disjoint
accessible terms of $T$ and $T/F_{\underline{v}}$ the resulting quotient. The
difference rests on different ways of obtaining the quotient term $T/F_{\underline{v}}$.
We will focus here on one form of the coproduct, called the ``deletion coproduct" $\Delta^d$ in
\cite{MCB}, where the quotient $T/F_{\underline{v}}$ is taken to be the maximal
full binary tree obtained by edge contractions (that eliminate non-branching vertices) on the
(non-full) binary tree obtained from $T$ by cutting the edge above each non-root vertex $v_i$. 
Another form of the coproduct, called the ``contraction coproduct" in \cite{MCB} uses a form of 
the quotient $T/F_{\underline{v}}$ where each component of the forest $F_{\underline{v}}$ is
contracted to its root vertex. (We will only discuss the contraction coproduct in \S \ref{ContractSec}.)

In addition to the different forms of the quotient $T/F_{\underline{v}}$, different forms of the
coproduct may also differ according to conditions (types of admissible cuts) on the 
extraction of the accessible terms $T_{v_i}$. In the case of the deletion coproduct that
we focus on in this paper, some of the algebraic properties of the coproduct change
according to whether we allow arbitrary admissible cuts that extract the forest $F_{\underline{v}}$
(where admissible means that no two edges of the cut lie on the same path from the root to one 
of the leaves, so that $F_{\underline{v}}$ indeed consists of accessible terms), or else one 
requires the stricter condition that the list $\underline{v}=(v_1,\ldots,v_n)$ of vertices does
not contain any pair of sister-vertices (one never cuts two edges below the same vertex). 
We will not discuss in this paper the difference between these two choices at the level of
algebraic properties, as that will be analyzed elsewhere. However, we will discuss
briefly how these two different choices for admissible cuts affect the properties of the
dynamical system and why, from this perspective, it is preferable to include the
possibility of cutting both edges below the same vertex.

\item The {\em Merge operation on workspaces} is a dynamical system on $\cV(\fF_{\cS\cO_0})$ generated by the
linear operators defined on basis elements $F$ by
\begin{equation}\label{MergeKF}
\cK_2(F) = \sqcup \circ (\fB \otimes {\rm id})\circ \Pi^{(2)} \circ \Delta \ \ \ \text{ and } \ \ \  
\cK_1(F) = \sqcup \circ \Pi^{(1)} \circ \Delta \, , 
\end{equation}
where $\fB$ is the grafting operator
from forests to trees $\fB(T\sqcup T')=\fM(T,T')$, 
and $\Pi^{(k)}$ is the linear projection on $\cV(\fF_{\cS\cO_0})\otimes \cV(\fF_{\cS\cO_0})$ that is the identity
on basis elements $F\otimes F'$ where $F=T_1\sqcup \cdots \sqcup T_k$ has $k$ components and is zero otherwise.
We can write it as a single operator
\begin{equation}\label{MergeKF12}
\tilde\cK(F)= \cK_2(F)+\cK_1(F) = \sqcup \circ (\fB \otimes {\rm id})\circ \tilde\Pi^{(2)} \circ \Delta \, , 
\end{equation}
where $\tilde\Pi^{(k)}=\oplus_{j=1,\ldots,k} \Pi^{(j)}$, and where, for a single tree $T$, we extend $\cB$ as
$\cB(T)=T$, so that the last expression in \eqref{MergeKF12} is indeed equal to $\cK_2(F)+\cK_1(F)$. 

While \eqref{MergeKF12} is the most
natural way of writing the Merge dynamical system from the mathematical perspective, as it has
exactly the form \eqref{HMCgen} of a Hopf algebra Markov chain, it is preferable for adherence
to a more transparent linguistic interpretation, to separate out the operator $\cK$ of \eqref{MergeKF} 
into a a sum
$$ \cK_2 = \sum_{S,S'} \fM_{S,S'} $$
over $S,S'\in \cS\cO$, with
\begin{equation}\label{opMerge}
 \fM_{S,S'} = \sqcup \circ (\fB \otimes {\rm id})\circ \delta_{S,S'} \circ \Delta \, , 
\end{equation} 
with $\delta_{S,S'}$ the Kronecker-delta
operator that is the identity on terms in the range of $\Delta$ that have $S\sqcup S'$ in
the left-hand-side of the coproduct (and zero otherwise). The sum over $S,S'$ always
reduces to a finite sum when applied to a workspace $F$. These operations
$\fM_{S,S'}$ correspond to the linguistic External Merge ($\fM_{S,S'}$ with both $S,S'$ connected
components of the workspace) and Sideward Merge ($\fM_{S,S'}$ where at least one of $S,S'$ is an accessible term). 

The operator $\cK_1$ of \eqref{MergeKF} can be similarly written in the form
$\cK_1=\sum_S \fM_{S,1}$ where $1$ is the unit of the magma $\cS\cO$ (the formal empty tree). 
These operators $\fM_{S,1}$ don't have a direct linguistic interpretation on their own (one can
call them {\em virtual Merge}), but they give rise to a very important operation, that accounts
for ``movement" or ``transformation" in syntax, namely Internal Merge, which is realized here
as a composition $\fM_{S, T/S}\circ \fM_{S,1}$.

In order to make our formulation here as close as possible to the linguistic interpretation, we
will write the Merge Markov chain $\cK$ in terms of the operations of External Merge (EM),
Internal Merge (IM) and Sideward Merge (SM), rather than in terms of EM, SM and virtual Merge (vM)
as in \eqref{MergeKF} and \eqref{MergeKF12}, so that we ensure that each of the arrows
in the Merge graph that describes the Markov chain dynamics has a linguistic interpretation.
All the arguments we present can be easily adapted to the case of \eqref{MergeKF12}. 

\item The Merge dynamical system generated by $\tilde\cK$ as in \eqref{MergeKF12}
defines a Hopf algebra Markov chain in the sense of \cite{DiPaRa}, \cite{Pang}, \cite{Pang2}.

\item The individual Merge operations $\fM_{S,S'}$ and $\fM_{S, T/S}\circ \fM_{S,1}$ can be
weighted according to certain cost functions (Minimal Search, Resource Restrictions) as
discussed in \cite{MCB} and \cite{MLH}. 

\item In addition to the core computational mechanism of syntax described by the 
Merge operations $\fM_{S,S'}$ and $\fM_{S, T/S}\circ \fM_{S,1}$, one can consider
a system of filtering on the syntactic objects constructed by the Merge action on
workspaces, which eliminates structures that are non-viable for semantic interpretation.
There are two main such filters, which are analyzed in \cite{ML} and in \cite{MHL},
one dealing with theta-role assignments and one with the structure of phases. 
As shown in \cite{ML} and in \cite{MHL} these filters can be modeled in terms
of a colored operad with a bud generating system in the sense of \cite{Giraudo}, 
and can equivalently be implemented
by pruning the Merge operations according to coloring rules. From the point of view
of the dynamical system we consider here, this pruning corresponds to eliminating some
of the arrows from the directed graph of the Markov chain dynamics, which correspond
to Merge operations that are not compatible with the coloring restrictions. We will
discuss this in \S \ref{ColorSec}.

\end{itemize}

We take this as the basic structure defining the dynamical system that we analyze in this paper.
Orbits of the action of Merge on workspaces, starting from a forest consisting only of labelled 
leaf vertices and no edges, correspond to linguistic derivations of sentence formation
as repeated applications of Merge operations, building syntactic structures from an unstructured
set of lexical items.

\section{The Merge graphs}\label{GraphSec}

The dynamical system defined by the action of Merge on workspaces can be described in terms
of a directed graph $\cG$ with vertices the forests $F\in \fF_{\cS\cO_0}$ with nonempty set of edges,
$E(F)\neq \emptyset$.
The graph $\cG$ has an edge from a vertex $F$ to a vertex $F'$ 
whenever there is a Merge operation with $F' = \fM_{S,S'}(F)$ or $F'=\fM_{S, T/S}(\fM_{S,1}(F))$. 

The reason why we consider only forests $F$ with non-trivial set of edges is because every
Merge operation (even if applied to a disjoint union of vertices with no edges) produces
a forest with some edges, so the forests consisting only of leaves with no edges are transient
for the dynamics that never returns to them, so they do not contribute to the main properties of the
dynamics like stationary distribution, cycles, etc. 

Since all the Merge operations (with this form of the coproduct) preserve the number of leaves
of the forest, the infinite graph $\cG$ breaks into a disjoint union of finite graphs
$$ \cG= \bigsqcup_{n\geq 2} \cG_n \, . $$
(We do not include here the case $n=1$ as there is no dynamics there.)
Moreover, the Merge operation does not alter the labeling of the leaves so each
graph $\cG_n$ is itself a disjoint union over all such labelings  
$$ \cG_n =\bigsqcup_{A \in {\rm Sym}^n(\cS\cO_0)} \cG_{n,A} \, . $$
The set of vertices of $\cG_{n,A}$ is the finite subset $\fF_{A,n}\subset \fF_{\cS\cO_0}$ of
forests $F\in \fF_{\cS\cO_0}$ with $n=\# L(F)$ for $L(F)$ the set of leaves, with set of labels
$\{ \alpha_\ell \}_{\ell\in L(F)}=A$, and with a non-empty set $E(F)$ of edges.

In our analysis we will restrict to a fixed $n\geq 2$ and a fixed set $A=\{ \alpha_1, \ldots, \alpha_n \}$
of lexical items, hence working only with one $\cG_{n,A}$ graph at a time. 

\subsection{Counting forests} \label{CountSec}

The first step in understanding the structure of the Merge graph $\cG_{n,A}$ is to
have a good counting formula for the number of vertices. Vertices are workspaces,
namely forests whose connected components are non-planar full binary rooted trees
with leaves decorated by the elements $\alpha_i$ of $A$.

\smallskip

In order to count vertices, namely count forests $F\in \fF_{A,n}$, we first
split the set $\fF_{A,n}$ into subsets, according to how the set $L(F)$ of leaves 
of a forest $F=T_1 \sqcup \cdots \sqcup T_r$ splits into the sets $L(T_i)$ of
leaves of the tree components $T_i$ of $F$. 

\smallskip

Let $\cP(n)$ be the set of all partitions of $n$ as a sum $n=k_1+\cdots + k_r$ of integers with $k_i\geq 1$.
We write $\wp=\{ k_1, \ldots, k_r \}$ (with possible multiplicities among the $k_i$) for a partition 
$n=k_1+\cdots+ k_r$.  

We write $\wp_{1^n}=\{ \underline{1}^n \}=\{ 1, \ldots, 1 \}$ for the partition $n=1 + \cdots + 1$ into
a sum of $1$'s and we let
\begin{equation}\label{Pprimen}
\cP'(n):=\cP(n)\smallsetminus \{ \wp_{1^n} \} \, .
\end{equation}
be the set of all partitions $n=k_1+\cdots+k_r$ where at least one of the $k_i$ satisfies $k_i\geq 2$.

\smallskip

\begin{lem}\label{VertPart}
The set of vertices $V(\cG_{n,A})=\fF_{A,n}$ of the directed graph $\cG_{n,A}$ is
partitioned as
\begin{equation}\label{VGpnA}
V(\cG_{n,A}) = \bigsqcup_{\wp\in \cP'(n)} V_{\wp,n,A}\, ,
\end{equation}
where for $\wp=\{ k_1,\ldots, k_r \}$, with $n=k_1+\cdots + k_r$, the set $V_{\wp,n,A}$ consists of 
all the vertices $F$ of $\cG_{n,A}$ that are forests of the form $F=T_1\sqcup \cdots \sqcup T_r$
with $T_i\in \fT_{\cS\cO_0,k_i}$, with nonempty set of edges.
\end{lem}

\proof
Given a workspace $F=T_1 \sqcup \cdots \sqcup T_r$ in $\fF_{A,n}$ with $n$ leaves with labels 
$A=\{ \alpha_1, \ldots, \alpha_n \}$, the set of leaves $L(F)$ also decomposes as a disjoint union
\begin{equation}\label{leafPart}
 L(F)= L(T_1) \sqcup \cdots \sqcup L(T_r) 
\end{equation} 
where $L(T_i)$ is the set of leaves of the connected component of $F$ given by the tree $T_i$.
If $k_i=\# L(T_i)$, we obtain a partition $\wp\in \cP(n)$ of the form
with $n=\# L(F)= k_1+ \cdots + k_r$. Thus, we can split the set $\fF_{A,n}$ into disjoint
subsets $V_{\wp,n,A}$ for $\wp\in \cP(n)$ where $F\in V_{\wp,n,A}$ are forests where
the leaves split in \eqref{leafPart} according to the partition $\wp$. 
Forests with nonempty set of edges, $E(F)\neq \emptyset$, correspond to 
partitions in $\cP'(n)$. 
\endproof

\smallskip

Partitions $\wp=\{ k_1, \ldots, k_r \} \in \cP(n)$ can be represented as Young diagrams,
consisting of $n$ boxes arranged in $r$ rows of lengths $k_i$ (drawn in
non-decreasing order). For example we have 
$$ \wp=\{ 4, 2, 2, 1 \} \in \cP(9) \Leftrightarrow  D(\wp) =\ydiagram{4,2,2,1}  \, .   $$
This partition is the image $\wp=p(F)$ for different workspaces
 $$ F= \Tree[ [ $\bullet$ $\bullet$  ] [ $\bullet$ $\bullet$ ] ]  \sqcup \Tree[ $\bullet$ $\bullet$ ] \sqcup \Tree[ $\bullet$ $\bullet$ ] \sqcup \bullet \ \ \text{ or } \ \  F= \Tree[ $\bullet$ [ $\bullet$ [ $\bullet$ $\bullet$ ] ] ] \sqcup \Tree[ $\bullet$ $\bullet$ ] \sqcup \Tree[ $\bullet$ $\bullet$ ] \sqcup \bullet  $$
 where the $\bullet$ slots at the leaves can be filled with a permutation of the elements of 
 $A=\{ \alpha_1, \ldots, \alpha_9 \}$. 
In particular,  Young diagrams account for partitioning of the leaves over the different components
   of the workspace but not for the different tree topologies of the components.
 To keep track of different possible assignments of labels $\alpha_i$ at the leaves,
   we can also consider Young tableaux on the alphabet $A$, such as
   $$ \begin{ytableau}
   \alpha_{i_1} & \alpha_{i_2} & \alpha_{i_3} & \alpha_{i_4} \\
   \alpha_{i_5} & \alpha_{i_6} \\
   \alpha_{i_7} & \alpha_{i_8} \\
   \alpha_{i_9}
   \end{ytableau} $$
   with $\alpha_{i_1},\ldots,\alpha_{ i_9}$ a permutation of $A=\{ \alpha_1,\ldots, \alpha_9 \}$.

\smallskip

For a partition $\wp=\{ k_1, \ldots, k_m \} \in \cP(n)$, the
{\em multinomial coefficient}
\begin{equation}\label{multinomial}
\mu_{\wp,n}:= \binom{n}{k_1,\ldots, k_r} = \frac{n!}{k_1 ! \, \cdots \, k_m !}  
\end{equation} 
is the size of the coset $S_n/S_\wp$ of the symmetric group $S_n$ by the
Young group of the partition,
\begin{equation}\label{YoungGrp}
 S_\wp := S_{k_1}\times \cdots \times S_{k_r} 
\end{equation} 
with the $k_i$ ordered as in the Young diagram $D(\wp)$.
Note that in \eqref{multinomial} we list the $k_i$ with multiplicities. When 
listing explicitly the multiplicities as $n=a_1 k_1+ \cdots a_r k_r$, we write
\begin{equation}\label{multinomial2}
\mu_{\wp,n}:= \binom{n}{\underbrace{k_1,\ldots, k_1}_{a_1\text{-times}} \ldots, 
\underbrace{k_r, \ldots, k_r}_{a_r\text{-times}}} = \frac{n!}{(k_1 !)^{a_1} \, \cdots \, (k_r !)^{a_r}}  \, . 
\end{equation} 

The {\em generalized multinomial coefficient} of a partition 
\begin{equation}\label{wpmulti}
\wp =\{ \underbrace{k_1, \ldots, k_1}_{a_1}, \ldots, \underbrace{k_r, \ldots, k_r}_{a_r} \}  \ \ \text{ with }\ \  n = \underbrace{k_1 + \cdots + k_1}_{a_1\text{-times}}+ \cdots +\underbrace{k_r+ \cdots + k_r}_{a_r\text{-times}} 
\end{equation}
is given by
\begin{equation}\label{genmultinomial}
\Upsilon_{\wp, n}:= \frac{1}{a_1 !\, \cdots\, a_r !}  \binom{n}{\underbrace{k_1,\ldots, k_1}_{a_1\text{-times}} \ldots, 
\underbrace{k_r, \ldots, k_r}_{a_r\text{-times}}} =\frac{1}{a_1 !\, \cdots\, a_r !} \,\,\,  \frac{n!}{(k_1!)^{a_1} \cdots (k_r!)^{a_r}}  \, .
\end{equation} 
We also consider the quantities
\begin{equation}\label{Lgenmultinomial}
\Lambda_{\wp,n}  := \frac{1}{a_1 !\, \cdots\, a_r !} \binom{n}{\underbrace{k_1,\ldots, k_1}_{a_1\text{-times}} \ldots, 
\underbrace{k_r, \ldots, k_r}_{a_r\text{-times}}} \prod_{i=1}^r ((2k_i -3) !!)^{a_i}
\end{equation}
We also write \eqref{Lgenmultinomial} equivalently as
\begin{equation}\label{Lgenmultinomial2}
\Lambda_{\wp,n}  = \mu_{\wp,n} \cdot \Gamma_{\wp,n}
\end{equation}
where $\mu_{\wp,n}$ is the multinomial coefficient \eqref{multinomial2} and
\begin{equation}\label{multinomialGamma}
\Gamma_{\wp,n}  := \frac{1}{a_1 !\, \cdots\, a_r !}  \prod_{i=1}^r ((2k_i -3) !!)^{a_i} \, . 
\end{equation}

\smallskip

\begin{lem}\label{CountFlem}
The number of vertices in $V(\cG_{n,A})=\fF_{A,n}$ is given by
\begin{equation}\label{ForestsAnCount}
 \Lambda_n := \# \fF_{A,n} =\sum_{\wp=\{ k_1,\ldots, k_r\} \in \cP'(n)} \Lambda_{\wp,n}   \, ,
\end{equation} 
with $\Lambda_{\wp,n}$ as in \eqref{Lgenmultinomial} the sum over the set $\cP'(n)$ of \eqref{Pprimen}.
\end{lem}

\proof
The number of non-planar full binary rooted trees with 
$n$ {\em labelled} leaves is given by the odd double factorial
\begin{equation}\label{NumLabelTrees}
 \# \fT_{A,n} = (2n-3) !!  \, .
\end{equation}
Thus, the number of forests of the form $F=T_1\sqcup \cdots \sqcup T_r$,
with $T_i\in \fT_{A,k_i}$ and $k_1+\cdots+ k_r=n$ is given by
\begin{equation}\label{NumForestsP}
\# \{ F\in \fF_{A,n}\,|\, F=T_1\sqcup \cdots \sqcup T_r, \, \, T_i\in \fT_{A,k_i} \}=
\Upsilon_{\wp,n}\,\, \prod_{i=1}^r ((2k_i -3) !! )^{a_i}
\end{equation} 
where the multinomial coefficients
account for how the labels are split between the components
and the normalization dividing by the factor $a_1 !\, \cdots\, a_r !$ account for the symmetry with
respect to exchanging rows of identical length in the Young diagram of the partition.
Thus, the total number of forests in $\fF_{A,n}$ is a sum 
\begin{equation}\label{LambdaAnp}
 \# \fF_{A,n} = \Lambda_n = \sum_{\wp\in \cP'(n)} \Lambda_{\wp,n} \, , 
\end{equation}
with $\Lambda_{\wp,n}$ as in \eqref{Lgenmultinomial}, with
$\Lambda_{\wp,n}  = \#V_{\wp,n,A}$.
\endproof

\smallskip

Note, for comparison, 
that the counting of the sets of non-planar full binary rootes trees with $n$ {\em undecorated} leaves is more complicated.
It is given by the $n$-th Wedderburn--Etherington numbers $W_n$, for which there is no closed formula, but  
their generating function is the solution to the equation
$$ \cT(x)=x + \frac{1}{2} \cT(x)^2 + \frac{1}{2} \cT(x^2) \, . $$
Correspondingly, the generating function for the binary rooted forests (with
unlabelled leaves) is then
$$ \cF(x) = \exp(\cT(x)) \, . $$

\subsection{Indegrees and outdegrees in the Merge graphs}\label{inoutdegSec}

As edges of the graphs $\cG$ and of its components $\cG_{n,A}$ we will consider
the Merge transformations, subdividing them into External Merge, Internal Merge,
and Sideward Merge arrows. In the case of Internal Merge, we will not consider the
case where the accessible term extracted is one of the two trees $T_1,T_2$ under
the root of a component tree $T=\fM(T_1,T_2)$ because, with the form of the coproduct
we are using, the result of IM on either $T_1$ or $T_2$ will simply reproduce $T$
and we do not include such identity maps as part of the dynamics. In the case of
Sideward Merge, we will restrict to considering only the ``minimal" Sideward Merge
transformations (in the sense discussed in \cite{MLH}), namely those that
only extract and Merge atomic components (single leaves), because these are
the lowest cost SM, in the cost counting of \cite{MCB}, \cite{MLH}. 

\smallskip

We first introduce some notation for convenience. Given a workspace $F$, we write:
\begin{itemize}
\item $\pi_0(F)$ for the set of connected components of $F$ with $\# \pi_0(F)=b_0(F)$,
\item $\pi_{0,E}(F)$ for the subset of connected components of $F$ consisting of a tree with 
nontrivial set of edges, and $c(F):=\# \pi_{0,E}(F)$,
\item $\cC(F)\subset L(F)\times L(F)$ for the subset 
of pairs $(\ell,\ell')$ of leaves of $F$ that belong to a cherry tree in $F$, 
(subtrees $T_v\subset F$ that consist of two sister leaves $\ell,\ell'$ and a vertex $v$ above them and two edges), with 
$d(F)=\# \cC(F)$. 
\item $\cC'(F)=\cC(F)\cap \pi_0(F)$ for the set of cherries of $F$ that are connected components 
(as opposed to accessible terms of components), with $d'(F)=\# \cC'(F)$.
\item $\tilde c(F)=c(F)-d'(F)$ for the number of connected components of $F$ that have at least $3$ leaves. 
\item $d''(F)=d(F)-d'(F)$ for the number of cherries that are accessible terms of components (rather than full components),
\item $\cC(T_i)$ for the set of cherry subtrees of a component $T_i$, with $d(T_i)=\#\cC(T_i)$.
\end{itemize}

\smallskip

We can describe the edges in the graph $\cG_{n,A}$  in the following way.

\begin{prop}\label{EdgesGnA}
The edges of the graph $\cG_{n,A}$ with vertex set $\fF_{A,n}$ are described in the following way.
Consider a vertex $F$ of $\cG_{n,A}$ with $b_0(F)=r$ connected components, $F=T_1\sqcup \cdots \sqcup T_r$,
respectively with $k_i=\# L(T_i)$ leaves, with $\sum_i k_i =n$. Without loss of generality, assume that 
$T_1, \ldots, T_{c(F)}$ are the components with nonempty set of edges, and among them $T_1, \ldots, T_{d'(F)}$
are the components that are cherries.
\begin{enumerate}
\item The number of outgoing External Merge (EM) arrows at the vertex $F$ is equal to
\begin{equation}\label{EMedgesN}
N^{\rm out}_{\rm EM}(F) = \left\{ \begin{array}{ll} \binom{r}{2} & r=b_0(F) \geq 2 \\
0 & r=1\, . \end{array}\right. 
\end{equation}
\item The number of incoming External Merge (EM) arrows at the vertex $F$ is equal to
\begin{equation}\label{EMedgesN2}
N^{\rm in}_{\rm EM}(F) =c(F) \, . 
\end{equation}
\item The number of outgoing Internal Merge (IM) arrows at the vertex 
$F=T_1\sqcup \cdots \sqcup T_r$ with $k_i=\# L(T_i)$ is equal to
\begin{equation}\label{IMedgesN}
N^{\rm out}_{\rm IM}(F) = \sum_{i=1}^{c(F)} (2k_i -4)   \, . 
\end{equation}
\item The number of incoming IM arrows at $F$ is equal to 
\begin{equation}\label{IMedgesN2}
N^{\rm in}_{\rm IM}(F) =\sum_{i=1}^{c(F)} (2k_i-4)  \, . 
\end{equation}
\item The number of outgoing minimal Sideward Merge arrows at $F$ is given by
\begin{equation}\label{minSMedgesN}
N^{\rm out}_{\rm SM,\, min}(F)= \sum_{i=d'(F)+1}^{c(F)} \binom{k_i}{2} + \sum_{i=1}^{c(F)} 
 k_i (r -c(F)) + \sum_{i\neq j, \,\, i,j =1}^{c(F)}  k_i k_j  \, ,
\end{equation}
if extraction of a pair of sister leaves from  a component tree with at least $3$ leaves is allowed  and
\begin{equation}\label{minSMedgesNb}
N^{\rm out}_{\rm SM,\, min}(F) - d''(F)
\end{equation}
if no extraction of both leaves below the same vertex is allowed, 
with $N^{\rm out}_{\rm SM,\, min}(F)$ as in \eqref{minSMedgesN}.
\item The number of incoming minimal Sideward Merge arrows at $F$ is zero if $d'(F)=0$ and otherwise
\begin{equation}\label{minSMedgesN2}
N^{\rm in}_{\rm SM,\, min}(F)= 6 d'(F) (n-c(F)) + 
2 \sum_{i=1}^{d'(F)} \sum_{a \neq b\,\, a,b\neq i\,\, a,b=1}^{c(F)} (2k_a -2) (2k_b -2) +
\end{equation}
$$
2 \sum_{i=1}^{d'(F)}\sum_{j\neq i, \,\, j=1}^{c(F)} (2k_j-2) (2k_j + 2+ n-c(F))
$$
when cutting sister leaves is allowed in a component with at least $3$ leaves,
while in the case where no cut of two leaves below the same vertex is allowed, 
the factor $2k_j + 2+ n-c(F)$ is replaced by $2k_j + 1+ n-c(F)$.
\end{enumerate}
\end{prop}

\proof (1) If $F=T$ has only one connected component then there is no outgoing EM arrow as EM merges two
components of the workspace. For $r\ge 2$ connected components the possible outgoing EM arrows correspond to all the possible
choices of two components of $F$ that are merged by EM, so there are $\binom{r}{2}$ arrows.

\smallskip

(2) Any nonplanar full binary tree $T$ with a nonempty set of edges can be written in a unique way as $T=\fM(T_1,T_2)$,
hence it is the image under EM of a workspace of the form $T_1 \sqcup T_2$. Suppose given a workspace 
$F=T_1\sqcup\cdots \sqcup T_r$, where $T_1, \ldots, T_{c(F)}$ have a non-empty set of edges, 
so that $T_i = \fM(T_{i,1}, T_{i,2})$ for $i=1,\ldots, c(F)$. We write $F= T_i \sqcup \hat F_i
=\fM(T_{i,1}, T_{i,2})\sqcup \hat F_i$ where $\hat F_i=\sqcup_{j\neq i} T_j$. 
Then $F$ is the target of $c(F)$ EM arrows with sources the workspaces $T_{i,1} \sqcup T_{i,2} \sqcup \hat F_i$, 
for $i=1,\ldots, c(F)$. 

\smallskip

(3) Each component $T_i$ of $F$ that has more than two edges is the source of a number of nontrivial (non-identity) 
IM arrows, one for each accessible term $T_v \subset T_i$ where the non-root vertex $v$ is not one of the two
vertices immediately below the root of $T_i$  (with the form of the coproduct we consider, IM applied to one of those
two accessible terms would just be the identity). If the component tree $T_i$ has $k_i$ leaves, there are $2k_i - 2$
non-root vertices hence $2k_i-4$ excluding the two below the root, so we obtain \eqref{IMedgesN}. 

\smallskip

(4) For $F$ to be the target of an IM arrow, one of the components $T_i=\fM(T_{i,1},T_{i,2})$ of $F$ 
with nontrivial set of edges must be obtained via IM from a tree of the form $T_{i,2}\lhd_e T_{I,1}$ or
$T_{i,1}\lhd_{e'} T_{I,2}$ where $\lhd_e$ is the insertion pre-Lie operation at an edge $e$. We consider
here as source of the arrows the workspaces where the $T_i$ component of $F$ is replaced by
one of these edge-insertions and all the other components of $F$ are the same. If $k_i=k_{i,1}+k_{i,2}=\# L(T_i)$
with $k_{i,1}=\# L(T_{i,1})$, $k_{i,2}=\# L(T_{i,2})$, there are $2k_{i,1}-2$ edges in $T_{i,1}$ where the insertion
$\lhd_e$ can take place and $2k_{i,2}-2$ edges in $T_{i,2}$, for a total of $2k_i-4$ possible edge insertions,
that is, $2k_i-4$ sources for IM arrows, for each component $T_i$ of $T$. This gives a total number of incoming
IM arrows as in \eqref{IMedgesN2}. 

\smallskip

(5) Minimal Sideward Merge operations select accessible terms or components that are single leaves and
combines them into a new component. The possible cases are: two single leaves that belong to the same
component that is not a cherry, two leaves belonging to two 
different components in $\pi_{0,E}(F)$, and one component of $F$ consisting of a single leaf and another
leaf from one of the components with non-empty set of edges. (The case of two components that are
both single leaves is already counted among the EM arrows and is not a case of SM.) Thus, the number
of outgoing minimal SM arrows is the sum of these three possibilities. 

\smallskip

The first case, 
adding over the $\tilde c(F)=c(F)-d'(F)$ components that are not cherries gives
$$ \sum_{i=d'(F)+1}^{c(F)} \binom{k_i}{2} \, , $$
while in the case where one does not allow any extraction of a pair of leaves in a cherry,
even if the cherry is part of a larger tree, one needs to subtract $d''(F)$.

\smallskip

The second case adds up to a contribution of
$$ \sum_{T_i\neq T_j \in \pi_{0,E}(F)} k_i \cdot k_j \, . $$

\smallskip 

The third case adds up to
$$ (r-c(F)) \cdot  \sum_{i=1}^{c(F)} k_i    $$
with the first factor giving the choice of a single leaf component and the second factor the choice
of a leaf in one of the components with non-empty set of edges. Thus the total number of
outgoing minimal SM arrows gives \eqref{minSMedgesN}.

\smallskip

(6) Since every minimal SM transformation generates a cherry component, if $F$ has $d'(F)=0$
there can be no incoming SM arrows. If $d'(F)\neq 0$ then each of the $d'(F)$ cherry component
can be generated by an SM transformation.
For a given cherry component of $F$, say $T_i=\{ \alpha, \beta \}$, the source workspace of an SM
arrow that creates $T_i$ will have one of the following possibilities for the leaves $\alpha$ and $\beta$:
\begin{enumerate}
\item one of $\alpha$ or $\beta$ is a single component and the other is inserted at an edge of a component
in $\pi_{0,E}(F)$,
\item one of $\alpha$ or $\beta$ is a single component and the other is joined to a single-leaf component of $F$ 
to form a cherry component,
\item $\alpha$ and $\beta$ are edge-inserted in two different components in $\pi_{0,E}(F)$,
\item one of $\alpha$ or $\beta$ is inserted in a component in $\pi_{0,E}(F)$ and the other is joined to a single leaf,
\item both $\alpha$ and $\beta$ are edge-inserted in the same component in 
$\pi_{0,E}(F)$ (without forming a cherry subtree, in the case where cutting the two edges 
below the same vertex inside a tree is not allowed) 
\item one of $\alpha$ or $\beta$ is joined to a single leaf to form a cherry and the other is inserted at an
edge of that cherry.
\end{enumerate} 
The first case accounts for
$$ 2 \sum_{i=1}^{d'(F)}\sum_{j\neq i, \,\, j=1}^{c(F)}  (2k_j-2) \, . $$
The second case accounts for
$$ 2 \cdot d'(F) \cdot  (n-c(F))\, . $$
The third case accounts for 
$$ 2 \sum_{i=1}^{d'(F)} \sum_{a \neq b\,\, a,b\neq i\,\, a,b=1}^{c(F)} (2k_a -2) (2k_b -2)\, . $$
The fourth case accounts to 
$$ 2 \cdot (n-c(F)) \sum_{i=1}^{d'(F)}\sum_{j\neq i, \,\, j=1}^{c(F)}  (2k_j-2) \, . $$
The fifth case accounts for 
$$ 2 \sum_{i=1}^{d'(F)} \sum_{j\neq i\, j =1}^{c(F)} (2k_j-2) (2k_j+1)  $$
where the factor $2k_j -2$ accounts for the choice of the edge of $T_j$ where the first insertion happens.
There are then $2k_j +1$ choices of an edges in $T_j \lhd_e \alpha$ where $\beta$ can be inserted.
In the case where we do not allow cutting of the two leaves of a cherry inside a tree component,
we only have $2k_j$ choices for the second edge, so the counting is replaced by
$$ 2 \sum_{i=1}^{d'(F)} \sum_{j\neq i\, j =1}^{c(F)} (2k_j-2) 2k_j \, . $$
The sixth case accounts for
$$ 4\cdot d'(F) \cdot (n-c(F)) \, . $$
This results in \eqref{minSMedgesN2}.
\endproof

\smallskip

As we will further discuss, one of the reasons why it is preferable, from the point of view of
the properties of the Merge dynamical system, to include the admissible cuts that 
remove both edges below a vertex, is that the counting of the in/out-degrees for the
SM arrows becomes then dependent only on data of the partition $\wp=p(F)$ rather than
on the specific topology of the forest $F$, while distinguishing between admissible
cuts that sever both edges below the same vertex or not depends on specific
information about $F$ and not just on the partition $\wp$. As we will discuss in \S \ref{DecompSec},
the fact that the in/out degrees only depend on $\wp$ greatly simplifies the dynamics.

\smallskip
\subsection{Sparsity}\label{sparseSec}

Proposition~\ref{EdgesGnA} also shows the {\em sparsity} property of these graphs.
Namely, as $n$ grows in size, the graph $\cG_{n,A}$ 
becomes increasingly sparse at a very fast rate. The density of a directed graph is defined as
\begin{equation}\label{Densegraph}
D(\cG) = \frac{ \# E(G)}{\# V(G) \cdot ( \# V(G) -1) } \, ,
\end{equation}
namely the ratio to the maximal number of edges. A family of graphs
$\cG_n$ is growing increasingly sparse if
$$ \lim_{n\to \infty} D(\cG_n) =0 \, . $$

\begin{cor}\label{sparse}
The density $D(\cG_{n,A})$ satisfies
\begin{equation}\label{Gdensity}
D(\cG_{n,A}) \leq \frac{P(n)}{\Lambda_n} \, ,
\end{equation}
with $P(n)$ a polynomial, at most cubic in $n$, and with $\Lambda_n$ as in \eqref{ForestsAnCount}.
\end{cor}

\proof The number of edges $\# E(\cG_{n,A})$ can be obtained from the
number of incoming and outgoing edges computed in Proposition~\ref{EdgesGnA} 
since, with $\deg$ denoting the valence of vertices,
$$ \# E(\cG_{n,A})=\frac{1}{2} \sum_{F \in V(\cG_{n,A})} \deg(F)  $$
$$ = \frac{1}{2} \sum_{F \in V(\cG_{n,A})} ( N^{\rm out}_{\rm EM}(F)  + N^{\rm in}_{\rm EM}(F) +
N^{\rm out}_{\rm IM}(F)  + N^{\rm in}_{\rm IM}(F) + N^{\rm out}_{\rm SM\, min}(F)  + N^{\rm in}_{\rm SM\, min}(F)) \, . $$
We can bound these terms by
$$  N^{\rm out}_{\rm EM}(F) \leq \binom{n}{2} =O(n^2)  , \ \ \  N^{\rm in}_{\rm EM}(F) \leq n $$
$$ N^{\rm out}_{\rm IM}(F) \leq 2 n (n-2) = O(n^2)   , \ \ \   N^{\rm in}_{\rm IM}(F) \leq 2 n (n-1) =O(n^2)  $$
$$ N^{\rm out}_{\rm SM\, min}(F) \leq  2 n(n-1) = O(n^2)  , \ \ \  
N^{\rm in}_{\rm SM\, min}(F) \leq    2n^2 (3 n+1) =O(n^3)   \, .   $$
We then obtain an overall polynomial growth estimate times the number of vertices 
$$ \# E(\cG_{n,A}) \leq P(n)\cdot \# V(\cG_{n,A})  \, . $$
Thus, the density satisfies
$$ D(\cG_{n,A}) \leq \frac{ P(n)}{\# V(\cG_{n,A})} \, . $$

The number of vertices $\# V(\cG_{n,A})$ is the number of forests $F\in \fF_{A,n}$. 
This can be computed as in Lemma~\ref{CountFlem}.
\endproof

\smallskip
\subsection{Strong connectedness} \label{strongconnGAnSec}

Recall that a directed graph $G$ is strongly connected if there is a a directed path of edges
between any pair of vertices in $G$. The graph $G$ is periodic 
if all the cycle-lengths in the graph have a common 
gcd greater than $1$ and is aperiodic if the gcd of the cycle-lengths is equal to $1$.

\begin{prop}\label{StrongConnAper}
The graphs $\cG_{n,A}$ are strongly connected and aperiodic for every choice of $(n,A)$.
\end{prop}

\proof The strongly connected property of the graphs  $\cG_{n,A}$,  
 was shown in \cite{MCB} and in \cite{MLH}, and relies on the presence of the minimal SM arrows.
To show aperiodicity, in the case $n=3$ one can see by direct inspection that there
are cycles of length $2$ and $3$, so the gcd of the cycle-lengths is $1$. For $n\geq 4$, 
notice that every graph $\cG_{n,A}$ contains a vertex $F$ that corresponds
to the workspace $F=T_n$ with a single component given by the comb tree
$$ T_n = \Tree[ $\alpha_n$ [ $\alpha_{n-1}$ [ $\cdots$ [ $\alpha_2$ $\alpha_1$ ] ] ] ] \, . $$ 
Among the possible IM transformations that act on this $F$ there is the one that extracts the
leaf $\alpha_1$ and results in a tree with the same comb topology but with a cyclic permutation
of the labels. One can then apply the same type of IM operation that extracts the $\alpha_2$ leaf
to the resulting tree and so on. 
This composition of Internal Merge operations gives a cycle of length $n-1$
$$ \Tree[ $\alpha_n$ [ $\alpha_{n-1}$ [ $\cdots$ [ $\alpha_2$ $\alpha_1$ ] ] ] ] \mapsto 
\Tree[ $\alpha_1$ [ $\alpha_n$ [ $\cdots$ [ $\alpha_3$ $\alpha_2$ ] ] ] ] \mapsto \cdots \mapsto  
 \Tree[ $\alpha_n$ [ $\alpha_{n-1}$ [ $\cdots$ [ $\alpha_2$ $\alpha_1$ ] ] ] ] \, . $$
The graph  $\cG_{n,A}$ also contains, for every $2< k < n$ a workspace of the form
$$ F = T_k \sqcup \alpha_{k+1} \sqcup \cdots \sqcup \alpha_n \, . $$
Internal Merge acts on these workspaces by acting on the $T_k$ component and we again
have a cycle of length $k$ as above. Thus  $\cG_{n,A}$ contains cycles of all lengths 
$k-1 =2,\ldots, n-1$, so again the gcd of the cycle-lengths is $1$. 
\endproof

\section{Random walks versus Hopf algebra Markov chains}\label{RWvsHMCsec}

There are two somewhat different ways in which one can  study the dynamics of the
Merge action on workspaces, given the Merge graph described above. One possible
approach consists of describing the dynamical systems of the Merge action as the
random walk on the Merge graph $\cG$ (which means the random walk on each
$\cG_{n,A}$ graph). The other is to study the Merge operator as a Hopf algebra
Markov chain (HAMC), as already discussed in \cite{MCB} and \cite{MLH}.

In this section we highlight how these two approaches differ, and why the second
is the right one to follow in our model. We do this by describing the 
main steps in constructing the associated dynamical system in both cases and
we highlight their main properties. In the
rest of the paper, we will focus only on the Hopf algebra Markov chain construction.

The key idea that describes the
difference is that, in the first case, one maintains the same basis of the combinatorial
Hopf algebra of workspaces (the forests $F$) and rescales the entries of the matrix 
representation of Merge in this basis (that is, the adjacency matrix of the Merge graph)
so that it is normalized to be a stochastic matrix, while in the second case one looks
for the appropriate basis in which the Merge operator itself is a Markov chain: the
right basis is obtained by rescaling the basis elements $F$ by the components of the
Perron--Frobenius eigenvector of the adjacency matrix of the Merge graph. So the
main difference is whether we study the Merge operator itself as a Markov chain,
as in the HAMC case, 
which reflects the linguistics setting where one considers the Markovian property
of Merge (see \cite{ChomskyGK}, \cite{ChomskyElements}), or we modify the
Merge operator by 
weights depending on the valencies of the
vertices of the Merge graphs (the number of possible Merge transformations
that can be applied to a given workspace), so as to obtain another operator that
is a Markov chain in the same combinatorial basis, as in the case of the
random walk operator. 

There are valid reasons to consider both possibilities: the usual claim in the
linguistics literature that the Merge operator itself is Markovian is rigorously
formalized by the HAMC approach described here, while the random walk approach
can be absorbed into the fact that a weighting of the Merge operations is
needed in any case, to describe cost functions and optimality properties, as we
will discuss more in detail in \S \ref{CostSec}, and the rescaling used in the
first approach can also be combined into an overall cost function (though
with a less transparent interpretation in terms of the linguistic model).

However, there is an additional property that makes the HAMC approach
preferable, namely the HAMC dynamical system satisfies an entropy 
optimization property, as we discuss in \S \ref{MaxEntSec}. In fact, as we
will discuss further in \S \ref{EntCostSec}, other forms of optimization with respect to
the Shannon entropy will play an important role in the dynamical properties
of the Merge action. 

For these reasons we will discuss both approaches here, but we will
only focus on the HAMC formulation. 

\subsection{Random walk on the Merge graph}\label{RWsec}

In the case of an undirected graph $G$, the random walk has a simple behavior.
It is the Markov chain with transition probabilities, for $v,v'\in V(G)$, 
\begin{equation}\label{undirRW}
\P(v,v') =\left\{ \begin{array}{ll} \frac{1}{\deg(v)} & \exists e\in E(G): \partial e =\{ v,v' \} \\
0 & \text{otherwise.}
\end{array} \right.
\end{equation}
The stationary distribution
\begin{equation}\label{statdistrRW}
 \pi(v')=\sum_{v\in V(G)} \pi(v) \P(v,v') 
\end{equation} 
is then just given by
$$ \pi(v)= \frac{\deg(v)}{2 \# E(G)} \, , $$
where we again use the fact that $\sum_{v\in V(G)}\deg(v)=2 \# E(G)$.

However, the random walk on a directed graph is more subtle. It is given by
the Markov chain with transition probabilities, for $v,v'\in V(G)$,
\begin{equation}\label{undirRW}
\P(v,v') =\left\{ \begin{array}{ll} \frac{1}{\deg^{\rm out}(v)} & \exists e\in E(G): v=s(e), v'=t(e)  \\
0 & \text{otherwise.}
\end{array} \right.
\end{equation}
where now only the outward directions are chosen at each vertex. In this case
one does not have a simple form of the stationary distribution in terms of valencies as in the undirected
case, since we would now have $\sum_v \deg^{\rm out}(v) \P(v,v') =\deg^{\rm in} (v')$. 
If the graph $G$ is strongly connected, however, it is
possible to obtain the stationary distribution through the Perron--Frobenius theorem.

The convergence to the stationary distribution then depends on the aperiodicity of the graph.
Periodicity with the gcd of the cycle-lengths equal to $k>1$ is equivalent to the existence 
of a morphism $G \to C_k$, where $C_k$ is a 
cyclically directed polygon graph of length $k$, and is in turn equivalent to the spectrum of the
transition matrix containing $k$-th roots of unity, and in turn implies that one does not have convergence by
iteration to the stationary distribution, but an oscillatory behavior. On the other hand, in the case of
aperiodicity, the Perron-Frobenius theorem ensures that all the other eigenvalues 
have smaller absolute value so the random walk always converges to a stationary distribution
given by the  Perron-Frobenius eigenvector \eqref{statdistrRW}. 

In the case of the Merge graph, this approach considers the random walk on the directed graphs $\cG_{n,A}$
described above. We know from Proposition~\ref{StrongConnAper} that the Merge graphs $\cG_{n,A}$ are both
strongly connected and aperiodic.
The original Merge operator is the linear operator on 
$\cV(\fF_{\cS\cO_0})$ represented in matrix form $\cK=(\cK_{F,F'})$ for 
the standard basis $\{ F \in \fF_{\cS\cO_0} \}$, on the invariant subspaces
$\cV(\fF_{A,n})$, by the adjacency matrix of the directed graph $\cG_{n,A}$. The random walk
operator is then the linear operator on the same vector space, represented in the same basis by 
rescaling the entries $\cK_{F,F'}\mapsto \deg^{out}(F)^{-1}\, \cK_{F,F'}$ with
$$ \deg^{out}(F)=N^{\rm out}_{\rm EM}(F)  +
N^{\rm out}_{\rm IM}(F)   + N^{\rm out}_{\rm SM\, min}(F)\, . $$ 
We use the notation $\cR\cW$ and 
\begin{equation}\label{RWoper}
 \cR\cW_{F,F'} := \deg^{out}(F)^{-1}\, \cK_{F,F'} 
\end{equation}  
to denote this random walk operator on the Merge graph, writing $\cR\cW^{(A,n)}$ and $\cR\cW_{F,F'}^{(A,n)}$
when we  want to explicitly keep track of the data $(A,n)$ of the graph $\cG_{n,A}$. 

From the point of view of the linguistic model that motivates our construction, the problem with the use of
the random walk operator $\cR\cW^{(A,n)}$ is that one does not have a natural theoretical justification for the rescaling 
of the entries of $\cK=(\cK_{F,F'})$ (the linguistic Merge operations) by the inverses $\deg^{out}(F)^{-1}$ 
the out-degrees: from the perspective of theoretical linguistics it is the Merge action itself that is Markovian. 
Thus, a better formulation is provided by the following approach. 

\subsection{Hopf algebra Markov chain} \label{HMCsec}

The other approach we mentioned above does not normalize the entries of the $\cK_{F,F'}$ matrix to make it stochastic,
rather it finds the appropriate 
change of basis in which the original linear operator is represented by a stochastic matrix. This can be
done, as discussed in \cite{DiPaRa}, \cite{Pang}. Since the graphs $\cG_{n,A}$ are strongly connected, we can apply
the Perron-Frobenius theorem that ensures the existence of right Perron-Frobenius eigenvector $\eta(F)=\eta_{n,A}(F)>0$,
\begin{equation}\label{rightPF}
\sum_{F'} \cK^{(A,n)}_{F,F'} \, \, \eta(F') = \lambda \, \eta(F)\, ,
\end{equation}
where $\lambda=\lambda_{n,A}$ is the Perron-Frobenius eigenvalue of the adjacency matrix $\cK^{(A,n)}$ of the
graph $\cG_{n,A}$. We can then rescale the basis elements by
$F \mapsto \eta(F)^{-1} F$. In this new basis $\{ \eta(F)^{-1} F \}$ the same Merge operator is represented
by the matrix $\hat \cK_{F,F'}^{(A,n)}$ of the form
\begin{equation}\label{hatKmat}
 \hat \cK_{F,F'}^{(A,n)} = \lambda^{-1} \frac{\eta(F')}{\eta(F)} \cK_{F,F'}^{(A,n)} \, . 
\end{equation}  
The matrix $\hat\cK^{(A,n)}=(\hat \cK^{(A,n)}_{F,F'} )$ obtained as in \eqref{hatKmat} is by construction stochastic, 
hence it defines the transition matrix of a Markov chain. This is the rigorous formulation of the ``Merge is
Markovian" property, with \eqref{hatKmat} the actual associated Hopf algebra Markov chain.

One can then study the stationary distribution \eqref{statdistrRW} with $\P(F,F')= \hat \cK_{F,F'}^{(A,n)}$,
and other properties of this Markov chain that describe the dynamics of Merge on workspaces.  
Proposition~\ref{StrongConnAper} already provides some important information about the Merge 
action on workspaces.

\smallskip

\begin{prop}\label{HAMCstatdistLRPF}
For a strongly connected graph $\cG$ with adjacency matrix $\cK$, with Perron-Frobenius
eigenvalue $\lambda$ and (right) eigenvector $\eta$, let $\hat\cK$ be the associated Markov chain 
defined as in \eqref{hatKmat}. The stationary distribution $\pi$ of $\hat\cK$,
\begin{equation}\label{pihatKstatdistr}
\sum_x \pi(x)\, \hat\cK_{x,y} = \pi(y) \, , 
\end{equation}
is given by
\begin{equation}\label{pihatKLRK}
\pi(x) = \psi(x) \eta(x) \, ,
\end{equation}
where $\psi, \eta$ are the left and right Perron-Frobenius eigenvectors of $\cK$.
 \end{prop}
 
\proof Since the graph $\cG$ is strongly connected, the Markov chain $\hat\cK$ is
irreducible hence its (normalized) stationary distribution is unique. Thus, it suffices
to see that $\psi(x) \eta(x)$ satisfies \eqref{pihatKstatdistr}. We have
$$ \sum_x  \psi(x) \eta(x) \lambda^{-1} \frac{\eta(y)}{\eta(x)} \cK_{x,y} =
\sum_x  \psi(x) \cK_{x,y} \lambda^{-1} \eta(y) = \psi(y) \eta(y) \, . $$
\endproof

We will prove in Proposition~\ref{projreg} a more refined version of this
simple observation, adapted to a projection between Markov chains.

\smallskip

\begin{rem}\label{ergodicitydef}{\rm
A Markov chain with transition matrix $\cK$ is {\em ergodic} if the stationary distribution is unique and 
any given initial distribution on the set of vertices of the underlying graph $G$
converges under iterations of $\cK$ to the stationary distribution. }
\end{rem}

\smallskip

Ergodicity, as formulated in Remark~\ref{ergodicitydef} above, follows
from strong connectedness and aperiodicity of $G$, hence Proposition~\ref{StrongConnAper}
has the following direct consequence.

\begin{prop}\label{ErgodicK}
The Merge Hopf algebra Markov chain $\hat\cK^{(A,n)}$ is ergodic, and so is the Merge random walk $\cR\cW^{(A,n)}$.
\end{prop}

\smallskip

This result carries interesting theoretical information about the Minimalist model of syntax.
It shows that syntax formation and transformation (what the Merge operator models) is an 
ergodic dynamical system. This means that the dynamics has a unique
stationary distribution that captures the long-term behavior of the dynamics (the long-term 
probability of being in one of the possible states, which are the workspaces, the syntactic structures).
Ergodic Markov chains are those for which the long-term behavior is {\em predictable}: 
ergodicity, for dynamical systems, is usually viewed as the property that
the time average is the same as the space average. In the case of Markov
chains, this can be more precisely expressed as in Remark~\ref{ergodicitydef}, as the property
that the long-term probability of being in a given state does not depend 
on the choice of the initial state. 
The values $\pi(F)$ of the stationary distribution also compute
the expected return times of the dynamics, which are given by $1/\pi(F)$. 

\smallskip

Note here, however, that these nice properties of the dynamics hold when we
include Sideward Merge as part of the dynamics (because the strong connectedness
property of the graphs $\cG_{n,A}$ require the SM arrows). We will discuss more
in detail in \S \ref{IMEMsec} what changes when the SM arrows are absent, and when the different types
of arrows, EM, IM, and SM are weighted differently by a suitable cost function.

\subsection{Hopf algebra Markov chain as maximal entropy random walk} \label{MaxEntSec}

We also point out here an important property that the HAMC formulation of Merge, as in \eqref{hatKmat}
satisfies, that the random walk operator $\cR\cW$ does not have: maximal entropy optimization.

Suppose given a directed graph $\cG$ with adjacency matrix $\cK_\cG$. Suppose that
$\cS$ is a Markov chain on the same underlying directed graph $\cG$, namely a weighted
version of the adjacency matrix $\cK_\cG$ that is a stochastic matrix. 

The entropy of a Markov chain  with transition probabilities $\cS(x,y)$ and with
stationary distribution $\pi(x)$ is given by
\begin{equation}\label{ShMarkov}
 {\rm Sh}(\cS) = - \sum_x \pi(x) \sum_y \cS(x,y) \log \cS(x,y) \, . 
\end{equation} 
This is also called the {\em entropy rate} of $\cS$.

The entropy rate ${\rm Sh}(\cS)$ is equal to the asymptotic behavior (in the length
of the path becoming large) of the average entropy of the probability distribution in the space of paths.
Namely, suppose given a directed path $\gamma = e_0 e_1 \ldots e_{\ell-1}$ in the graph $\cG$,
with source and target vertices of the edges $s(e_i)=x_i$ and $t(e_i)=x_{i+1}$. The
probability distribution on paths of length $\ell$ is given by
\begin{equation}\label{Pellgamma}
 \P_\ell (\gamma) = \pi(x_0)\, \cS(x_0,x_1) \cdots \cS(x_{\ell-1}, x_\ell) \, . 
\end{equation} 
The entropy rate  \eqref{ShMarkov} is the limit
\begin{equation}\label{ShSlimit}
 {\rm Sh}(\cS) = \lim_{\ell \to \infty} \frac{{\rm Sh}(\P_\ell)}{\ell} \, , 
\end{equation} 
where ${\rm Sh}(\P_\ell)$ is the Shannon entropy of the probability \eqref{Pellgamma}. 

When $\cS$ varies over all Markov chains with the same
underlying directed graph $\cG$, it is known (see \cite{Parry}, \cite{Ruelle}) that the maximal 
possible entropy rate ${\rm Sh}(\cS)$ is the topological entropy of the directed graph $\cG$,
which is given by
\begin{equation}\label{StopG}
S_{\rm top}(\cG):= \log \lambda_{\cG}\, , 
\end{equation}
where $\lambda_{\cG}$ is the Perron-Frobenius eigenvalue of the adjacency matrix $\cK_\cG$. 

The following fact is well known (see for instance \cite{Cover}, \cite{DelLib}, \cite{Parry}, \cite{Ruelle}) but we
recall it here for the reader's convenience.

\begin{prop}\label{MERWprop}
Suppose given a strongly connected directed graph $\cG$. Let $\cK_\cG$
be the adjacency matrix of $\cG$, with Perron-Frobenius eigenvalue $\lambda$ and
with $\eta$ and $\psi$ the right and left Perron-Frobenius eigenvectors, normalized
so that 
\begin{equation}\label{psietanorm}
\sum_{x,y} (\psi(x)\eta(y))^2 =1 \, .
\end{equation}
The maximal entropy random walk (MERW) is the Markov chain 
\begin{equation}\label{MERWK}
\cS_{\rm MERW}(x,y)= \hat\cK_\cG(x,y) =\frac{1}{\lambda} \frac{\eta(y)}{\eta(x)} \cK_\cG(x,y) \, .
\end{equation}
The $\cS_{\rm MERW}$ satisfies the following properties.
\begin{enumerate}
\item The probability distribution \eqref{Pellgamma}, with $\cS=\cS_{\rm MERW}$, 
is uniform on all paths with the same length and endpoints.
\item The entropy rate is maximal, 
$$ {\rm Sh}(\cS_{\rm MERW}) = S_{\rm top}(\cG) \, . $$
\end{enumerate}
\end{prop}

\proof One defines the MERW Markov chain as \eqref{MERWK}, and 
one can see that it determines a uniform distribution on the space of paths with fixed endpoints and length
by observing that by \eqref{Pellgamma} we have
\begin{equation}\label{PgammaMERW}
  \P_{{\rm MERW},\ell} (\gamma) = \frac{1}{\lambda^\ell}  \pi(x_0)  \frac{\eta(x_\ell)}{\eta(x_0)} \cK_\cG(x_0,x_1) \cdots \cK_\cG(x_{\ell-1},x_\ell) = 
 \frac{ \psi(x_0) \eta(x_\ell) }{\lambda^\ell}  \, ,  
\end{equation} 
 where in the second equality we used Proposition~\ref{HAMCstatdistLRPF} and the fact that, 
 for a directed edge $e_i$ in $\cG$ with $s(e_i)=x_i$ and $t(e_i)=x_{i+1}$, we have $\cK_\cG(x_i, x_{i+1})  =1$. 
 To see that MERW also maximizes the entropy rate, we use the limit \eqref{ShSlimit}. We have
 $$ {\rm Sh}(\P_{{\rm MERW},\ell}) = - \sum_\gamma  \P_{{\rm MERW},\ell} (\gamma) \log  \P_{{\rm MERW},\ell} (\gamma)\, , $$
 with the sum over all paths $\gamma$ of length $\ell$. This is
 $$ {\rm Sh}(\P_{{\rm MERW},\ell}) = - \sum_{x_0,x_\ell} \cK^\ell_\cG(x_0,x_\ell)  \frac{ \psi(x_0) \eta(x_\ell) }{\lambda^\ell} \log\left(\frac{ \psi(x_0) \eta(x_\ell) }{\lambda^\ell} \right)$$ $$ = - \sum_{x_0,x_\ell} \frac{\cK^\ell_\cG(x_0,x_\ell)}{\lambda^\ell} \, \psi(x_0) \eta(x_\ell) \log( \psi(x_0) \eta(x_\ell) ) +\sum_{x_0,x_\ell}
\frac{\cK^\ell_\cG(x_0,x_\ell)}{\lambda^\ell} \, \psi(x_0) \eta(x_\ell) \log(\lambda^\ell)  \, , $$
 where $\cK^\ell_\cG(x,y)=\#\{ \gamma \,|\, s(\gamma)=x, \, t(\gamma)=y\, , {\rm length}(\gamma)=\ell \}$.
Thus, we have
$$ {\rm Sh}(\cS_{\rm MERW}) =\lim_{\ell \to\infty} \frac{{\rm Sh}(\P_{{\rm MERW},\ell})}{\ell}  = \lim_{\ell \to\infty}
\sum_{x,y} \frac{\cK^\ell_\cG(x,y)}{\lambda^\ell} \, \psi(x) \eta(y) \log(\lambda)\, . $$
The Perron-Frobenius theorem shows that, when $\ell\to \infty$
$$ \frac{\cK^\ell_\cG(x,y)}{\lambda^\ell} \sim \psi(x) \eta(y)\, , $$
hence the above limit and the normalization \eqref{psietanorm} give
$$ {\rm Sh}(\cS_{\rm MERW}) =\sum_{x,y} ( \psi(x) \eta(y))^2 \,\, \log(\lambda) = \log(\lambda) = S_{\rm top}(\cG) \, . $$
\endproof

Note how the transition probabilities of MERW are non-local (unlike those of the random walk $\cR\cW$ that only
depend on the (out)degrees at vertices) and depend on the whole graph topology through the
Perron-Frobenius eigenvalue and eigenvector of the adjacency matrix. 

\smallskip
\subsection{Free energy optimization}\label{FreeEnSec}

We also recall,  for later use, a generalization of Proposition~\ref{MERWprop}, where instead of
the adjacency matrix one uses a weighted version and MERW induces a Boltzmann distribution on
paths and correspondingly optimizes the free energy instead of the entropy rate (see \cite{Ruelle}). 

\smallskip

Let $\cG$ be a strongly connected directed graph and let $\cK^{(\omega,\beta)}_\cG$ be a 
weighted adjacency matrix, namely
\begin{equation}\label{wAdjK}
\cK^{(\omega,\beta)}_\cG(x,y)=e^{-\beta \, \omega(x,y)}\,\, 
\cK_\cG(x,y) \, ,
\end{equation}
where $\omega(x,y)\in \R$ are real weights and $\cK_\cG$ is the adjacency matrix of $\cG$.
We define the energy function $\bE: E(\cG) \to \R_+$ on edges of $\cG$ as
$$ \bE(e)=\omega(x,y) \, ,  \ \ \ \text{ for } x=s(e)  \text{ and } y=t(e) \, . $$
Even though the weights are not necessarily positive, the energy is bounded below since the graph is finite,
hence it is non-negative after a global shift. 
Let $\cS$ be a Markov chain with underlying directed graph $\cG$.
The {\em free energy} $\bF(\cS)$ of the Markov chain is defined as 
\begin{equation}\label{freeendef}
 \bF(\cS)=\overline{\bE}(\cS) -\beta^{-1} {\rm Sh}(\cS)
\end{equation} 
where $\bar\bE(\cS)$ is the expected energy   
\begin{equation}\label{barE}
\overline{\bE}(\cS) := \sum_x \pi^\cS(x) \sum_y \cS(x,y) \, \omega(x,y) \, ,
\end{equation}
with $\pi^\cS$ the stationary distribution of $\cS$.
One looks for Markov chains that {\em minimizes} $\bF(\cS)$. Note that in the case with no weights, $\omega=0$,
this is equivalent to maximizing the entropy  ${\rm Sh}(\cS)$ so one recovers the previous case. 
Sometimes one rephrases this question equivalently in terms of a  
{\em maximization} of the free energy $\tilde\bF$ written as 
\begin{equation}\label{freeendefMax}
 \tilde\bF(\cS)={\rm Sh}(\cS)  - \beta \overline{\bE}(\cS) 
\end{equation} 
as in the formulation in  \cite{Ruelle}, while minimization of \eqref{freeendef} that we follow here 
is more suitable for direct comparison with the thermodynamics semirings used in \cite{MarBer}. 

\smallskip

\begin{defn}\label{BoltzDef}
The {\em Boltzmann distribution} on the set of edges of $\cG$ is given by
\begin{equation}\label{PeBoltz}
\bP_{(\omega,\beta)}(e) :=\frac{ \psi^{(\omega,\beta)}(s(e))\eta^{(\omega,\beta)}(t(e)) e^{-\beta \bE(e)}}{Z_\omega(\beta)}  \ \ \ \text{ with } Z_\omega(\beta)= \sum_e \psi^{(\omega,\beta)}(s(e))\eta^{(\omega,\beta)}(t(e))  \,\, e^{-\beta \bE(e)} \, ,
\end{equation}
with  $\psi^{(\omega,\beta)}$ and $\eta^{(\omega,\beta)}$
the left and right Perron-Frobenius eigenvectors of $\cK^{(\omega,\beta)}_\cG$, with
Perron-Frobenius eigenvalue $\lambda_{(\omega,\beta)}$. We use the normalization
\begin{equation}\label{piNorm1}
\sum_x  \psi^{(\omega,\beta)}(x)\eta^{(\omega,\beta)}(x) =1 \, .
\end{equation}
\end{defn}

\smallskip

\begin{lem}\label{lemZlambda}
The partition function of the Boltzmann distribution \eqref{PeBoltz} satisfies
\begin{equation}\label{Zlambda}
 Z_\omega(\beta)= \lambda_{\omega,\beta} 
\end{equation}
\end{lem}

\proof We have
$$ \sum_e \psi^{(\omega,\beta)}(s(e))\eta^{(\omega,\beta)}(t(e))  \,\, e^{-\beta \bE(e)}=\sum_{x,y} 
\psi^{(\omega,\beta)}(x)\eta^{(\omega,\beta)}(y) e^{-\beta \omega(x,y)}\cK_\cG(x,y) = $$ $$
\lambda_{\omega,\beta}  \sum_x \psi^{(\omega,\beta)}(x)\eta^{(\omega,\beta)}(x) = \lambda_{\omega,\beta}  $$
with the normalization \eqref{piNorm1}.
\endproof

\smallskip

We denote by $\hat\cK^{(\omega,\beta)}_\cG$  the Markov chain associated to the
weighted adjacency matrix by
\begin{equation}\label{hatKbeta}
\hat\cK^{(\omega,\beta)}_\cG(x,y) := \frac{1}{\lambda_{(\omega,\beta)}} \frac{\eta^{(\omega,\beta)}(y)}{\eta^{(\omega,\beta)}(x)}
\,\, \cK^{(\omega,\beta)}_\cG(x,y) \, . 
\end{equation}
This has the following property.

\begin{lem}\label{RWfreeEn}
For the Markov chain $\hat\cK^{(\omega,\beta)}_\cG$, the 
expected energy satisfies
\begin{equation}\label{BolzEbar}
 \bF(\hat\cK^{(\omega,\beta)}_\cG)= \overline{\bE}_{\bP_{(\omega,\beta)}} -\beta^{-1} {\rm Sh}(\bP_{(\omega,\beta)})\, ,
\end{equation}
with $\bP_{(\omega,\beta)}$ the Boltzmann distribution \eqref{PeBoltz} and 
$$ \overline{\bE}_{\bP_{(\omega,\beta)}} =\sum_e \bP_{(\omega,\beta)}(e) \, \bE(e) \, . $$
 \end{lem}

\proof We have
$$ \bF(\hat\cK^{(\omega,\beta)}_\cG)= \sum_x \pi^{(\omega,\beta)}(x) \sum_y \hat\cK^{(\omega,\beta)}_\cG(x,y)( \omega(x,y) +\beta^{-1} \log (\hat\cK^{(\omega,\beta)}_\cG(x,y)) ) $$
$$ = \sum_x \pi^{(\omega,\beta)}(x) \sum_y \hat\cK^{(\omega,\beta)}_\cG(x,y)( \omega(x,y) + \beta^{-1} \log(e^{-\beta \omega(x,y)}) + \beta^{-1}  \log (\hat\cK_\cG(x,y)) ) $$
$$ = \beta^{-1}  \sum_x \pi^{(\omega,\beta)}(x) \sum_y  e^{-\beta \omega(x,y)} \hat\cK_\cG(x,y)  \log (\hat\cK_\cG(x,y)) )\, . $$

By  Lemma~\ref{lemZlambda} and Proposition~\ref{HAMCstatdistLRPF} applied to $\hat\cK^{(\omega,\beta)}_\cG$, we also have
 $$ \overline{\bE}_{\bP_{(\omega,\beta)}}  = \sum_{x,y} \frac{\psi^{(\omega,\beta)}(x) \eta^{(\omega,\beta)}(y)}{\lambda_{(\omega,\beta)}}  e^{-\beta \omega(x,y)} \, \cK_\cG(x,y)\, \omega(x,y) = \sum_x \pi^{(\omega,\beta)}(x) \sum_y e^{-\beta \omega(x,y)}\, \hat\cK_\cG(x,y)\, \omega(x,y) $$
 and
 $$ -\beta^{-1} {\rm Sh}(\bP_{(\omega,\beta)}) =
  \beta^{-1} \sum_x \pi^{(\omega,\beta)}(x) \sum_y  \hat\cK^{(\omega,\beta)}_\cG (x,y) \log (\hat\cK^{(\omega,\beta)}_\cG(x,y)) $$
$$ =    \beta^{-1} \sum_x \pi^{(\omega,\beta)}(x) \sum_y \hat\cK^{(\omega,\beta)}_\cG(x,y) \log(e^{-\beta \omega(x,y)}) +
  \beta^{-1} \sum_x \pi^{(\omega,\beta)}(x) \sum_y e^{-\beta \omega(x,y)} \hat\cK_\cG(x,y) \log(\hat\cK_\cG(x,y))
 $$ 
 $$ =  - \sum_x \pi^{(\omega,\beta)}(x) \sum_y e^{-\beta \omega(x,y)}\, \hat\cK_\cG(x,y)\, \omega(x,y) 
 +  \beta^{-1} \sum_x \pi^{(\omega,\beta)}(x) \sum_y e^{-\beta \omega(x,y)} \hat\cK_\cG(x,y) \log(\hat\cK_\cG(x,y)) \, , $$
 so that
 $$ \overline{\bE}_{\bP_{(\omega,\beta)}} -\beta^{-1} {\rm Sh}(\bP_{(\omega,\beta)}) = 
 \beta^{-1} \sum_x \pi^{(\omega,\beta)}(x) \sum_y e^{-\beta \omega(x,y)} \hat\cK_\cG(x,y) \log(\hat\cK_\cG(x,y))\, . $$
\endproof

We also have the following properties of the Boltzmann distribution.

\begin{lem}\label{BoltzFESh}
The Boltzmann distribution $\bP_{(\omega,\beta)}$ of \eqref{PeBoltz}
satisfies
\begin{equation}\label{BoltzOpt1}
\overline{\bE}_{\bP_{(\omega,\beta)}} = - \frac{d}{d\beta} \log Z_\omega(\beta) \, , 
\end{equation}
\begin{equation}\label{BoltzOpt2}
{\rm Sh}(\bP_{(\omega,\beta)}) = \log Z_\omega(\beta) - \beta \frac{d}{d\beta} \log Z_\omega(\beta) 
+ {\rm Sh}(\pi^{(\omega,\beta)}) \, , 
\end{equation}
\begin{equation}\label{BoltzOpt3}
\bF_{\bP_{(\omega,\beta)}}  = \overline{\bE}_{\bP_{(\omega,\beta)}} -\beta^{-1} {\rm Sh}(\bP_{(\omega,\beta)})= -\beta^{-1}({\rm Sh}(\pi^{(\omega,\beta)})+ \log Z_\omega(\beta)) \, .
\end{equation}
\end{lem}

\proof
We obtain \eqref{BoltzOpt1} as the identity
$$ \overline{\bE}_{\bP_{(\omega,\beta)}} = \sum_e \bP_{(\omega,\beta)}(e)\, \bE(e)=
\frac{\sum_e \psi^{(\omega,\beta)}(s(e)) \eta^{(\omega,\beta)}(t(e)) e^{-\beta \bE(e)} \bE(e)}{Z_\omega(\beta)} $$ $$ =
\frac{ -\frac{d}{d\beta} \sum_e \psi^{(\omega,\beta)}(s(e)) \eta^{(\omega,\beta)}(t(e)) e^{-\beta \bE(e)}}{Z_\omega(\beta)} =
-\frac{d}{d\beta} \log Z_\omega(\beta) \, . $$
We also have 
$$ {\rm Sh}(\bP_{(\omega,\beta)}) =-\sum_e \frac{\psi^{(\omega,\beta)}(s(e)) \eta^{(\omega,\beta)}(t(e)) e^{-\beta \omega(e)}}{Z_\omega(\beta)} \log( \frac{\psi^{(\omega,\beta)}(s(e)) \eta^{(\omega,\beta)}(t(e)) e^{-\beta \omega(e)}}{Z_\omega(\beta)} ) $$
$$ = \beta\, \sum_e \frac{\psi^{(\omega,\beta)}(s(e)) \eta^{(\omega,\beta)}(t(e)) e^{-\beta \omega(e)}}{Z_\omega(\beta)} 
\omega(e) $$
$$ + \sum_e \frac{\psi^{(\omega,\beta)}(s(e)) \eta^{(\omega,\beta)}(t(e)) e^{-\beta \omega(e)}}{Z_\omega(\beta)} \log Z_\omega(\beta) $$
$$ - \sum_e \frac{\psi^{(\omega,\beta)}(s(e)) \eta^{(\omega,\beta)}(t(e)) e^{-\beta \omega(e)}}{Z_\omega(\beta)} \log( \psi^{(\omega,\beta)}(s(e)) \eta^{(\omega,\beta)}(t(e)) ) $$
$$ = \beta \overline{\bE}_{\bP_{(\omega,\beta)}} + \log Z_\omega(\beta) - \sum_{x,y} \frac{\psi^{(\omega,\beta)}(x)\eta^{(\omega,\beta)}(y) \cK^{(\omega,\beta)}_\cG(x,y)}{Z_\omega(\beta)} \log( \psi^{(\omega,\beta)}(x)\eta^{(\omega,\beta)}(y) ) $$
$$ = \log Z_\omega(\beta) - \beta \frac{d}{d\beta} \log Z_\omega(\beta)  - \sum_{x,y} \frac{\psi^{(\omega,\beta)}(x)\eta^{(\omega,\beta)}(y) \cK^{(\omega,\beta)}_\cG(x,y)}{\lambda_{(\omega,\beta)}} \log( \psi^{(\omega,\beta)}(x)) $$ $$  - \sum_{x,y} \frac{\psi^{(\omega,\beta)}(x)\eta^{(\omega,\beta)}(y) \cK^{(\omega,\beta)}_\cG(x,y)}{\lambda_{(\omega,\beta)}} \log(\eta^{(\omega,\beta)}(y)) = \log Z_\omega(\beta) - \beta \frac{d}{d\beta} \log Z_\omega(\beta) $$ $$ 
 - \sum_x \psi^{(\omega,\beta)}(x)\eta^{(\omega,\beta)}(x) \log( \psi^{(\omega,\beta)}(x)) 
 - \sum_y \psi^{(\omega,\beta)}(y)\eta^{(\omega,\beta)}(y) \log(\eta^{(\omega,\beta)}(y)) $$
 $$ = \log Z_\omega(\beta) - \beta \frac{d}{d\beta} \log Z_\omega(\beta) -  
 \sum_x \psi^{(\omega,\beta)}(x)\eta^{(\omega,\beta)}(x) \log( \psi^{(\omega,\beta)}(x)\eta^{(\omega,\beta)}(x)) \, , $$
 which gives \eqref{BoltzOpt2}, from which \eqref{BoltzOpt3} then follows. 
\endproof

We also have the following generalization of Proposition~\ref{MERWprop}.

\begin{prop}\label{MERWweight}
Suppose given a strongly connected directed graph $\cG$ with weighted adjacency matrix $\cK^{(\omega,\beta)}_\cG$ as
above, and
with Perron-Frobenius eigenvalue $\lambda_{(\omega,\beta)}$ and left and right 
eigenvectors $\psi^{(\omega,\beta)}$ and $\eta^{(\omega,\beta)}$. 
Then the weighted MERW on $\cG$ is the Markov chain \eqref{hatKbeta},
\begin{equation}\label{MERWKw}
\cS^{(\omega,\beta)}_{\rm MERW}(x,y)= \hat\cK^{(\omega,\beta)}_\cG(x,y) =\frac{1}{\lambda_{(\omega,\beta)}} 
\frac{\eta^{(\omega,\beta)}(y)}{\eta^{(\omega,\beta)}(x)} \cK^{(\omega,\beta)}_\cG(x,y) \, .
\end{equation}
The weighted MERW has the following properties.
\begin{enumerate}
\item The probability on the space of paths induced by $\hat\cK^{(\omega,\beta)}_\cG$ is the 
Boltzmann distribution 
\begin{equation}\label{BoltzmannMERW}
\P_{(\omega,\beta),\ell}(\gamma)=  \frac{1}{Z_{\omega,\ell}(\beta)}    \psi^{(\omega,\beta)} (x_0)\eta^{(\omega,\beta)}(x_\ell) 
 \exp(-\beta\, \bE(\gamma))\, ,
\end{equation}
for a path $\gamma=e_0 e_1\ldots e_\ell$ with $s(e_i)=x_i$ and $t(e_i)=x_{i+1}$, with energy 
$$ \bE(\gamma):= \sum_{i=0}^\ell \omega(x_i,x_{i+1}) \, , $$
with the partition function
$$  Z_{\omega,\ell}(\beta) = \sum_\gamma \psi^{(\omega,\beta)} (x_0)\eta^{(\omega,\beta)}(x_\ell) \exp(-\beta\, \bE(\gamma)) = \lambda_{(\omega,\beta)}^\ell  \, , $$
with the sum taken over paths $\gamma$ of length $\ell$.
\item The free energy $\bF_{\bP_{(\omega,\beta),\ell}}$ of the Boltzmann distribution on paths of length $\ell$ satisfies 
\begin{equation}\label{FreeEn}
\lim_{\ell\to \infty} \frac{ \bF_{\bP_{(\omega,\beta),\ell}} }{\ell} = \bF(\cS^{(\omega,\beta)}_{\rm MERW}) +
\beta^{-1} {\rm Sh}(\pi^{(\omega,\beta)})  \, , 
\end{equation}
where
 \begin{equation}\label{FreeEn2}
\bF_{\bP_{(\omega,\beta),\ell}} := -\beta^{-1}  {\rm Sh}(\P_{(\beta,\omega),\ell}) + \sum_\gamma \P_{(\beta,\omega),\ell}(\gamma)\, \bE(\gamma) \, , 
\end{equation}
and the limit \eqref{FreeEn} gives the optimal value 
\begin{equation}\label{maxFreeEn}
 \lim_{\ell\to \infty} \frac{ \bF_{\bP_{(\omega,\beta),\ell}} }{\ell} =- \beta^{-1} \log \lambda_{(\omega,\beta)} \, . 
\end{equation}
\end{enumerate}
\end{prop} 

\proof The argument follows the same steps as in Proposition~\ref{MERWprop}, but replacing the
adjacency matrix $\cK_\cG$ with the weighted version $\cK^{(\omega,\beta)}_\cG(x,y)
=\exp(-\beta \, \omega(x,y)) \cK_\cG(x,y)$. The distribution induced by $\hat\cK^{(\omega,\beta)}_\cG$ 
on paths is given by
$$ \P_{(\omega,\beta),\ell}(\gamma)=  \pi^{(\omega,\beta)}(x_0) 
\hat\cK^{(\omega,\beta)}_\cG(x_0,x_1) \cdots \hat\cK^{(\omega,\beta)}_\cG(x_{\ell-1},x_\ell) $$
$$ = \frac{1}{\lambda^\ell_{(\omega,\beta)}} \psi^{(\omega,\beta)}(x_0) \eta^{(\omega,\beta)}(x_\ell) 
\cK^{(\omega,\beta)}_\cG(x_0,x_1) \cdots \cK^{(\omega,\beta)}_\cG(x_{\ell-1},x_\ell) =\frac{1}{\lambda^\ell_{(\omega,\beta)}} \psi^{(\omega,\beta)}(x_0) \eta^{(\omega,\beta)}(x_\ell) e^{-\beta \bE(\gamma)}  \, . $$
 To show that this agrees with \eqref{BoltzmannMERW} it suffices to show that
 $Z_{\omega,\ell}(\beta)=\lambda_{(\omega,\beta)}^\ell$. This follows as in Lemma~\ref{lemZlambda} with
$$ Z_{\omega,\ell}(\beta)= \sum_{x,y} \psi^{(\omega,\beta)}(x) \eta^{(\omega,\beta)}(y) (\cK^{(\omega,\beta)}_\cG)^\ell(x,y) 
= \lambda_{(\omega,\beta)}^\ell \sum_x \psi^{(\omega,\beta)}(x) \eta^{(\omega,\beta)}(x) = \lambda_{(\omega,\beta)}^\ell  \, . $$
 The free energy of the distribution $\P_{(\omega,\beta),\ell}$ on the set of paths of length $\ell$ is given by
 $$ \bF_{\bP_{(\omega,\beta),\ell}} = -\beta^{-1}  {\rm Sh}(\P_{(\beta,\omega),\ell}) + \sum_\gamma \P_{(\beta,\omega),\ell}(\gamma)\, \bE(\gamma) $$
 where
 $$ -\beta^{-1}  {\rm Sh}(\P_{(\beta,\omega),\ell}) =
  \beta^{-1} \sum_{x,y} \psi^{(\omega,\beta)}(x) \eta^{(\omega,\beta)}(y)  \frac{(\cK^{(\omega,\beta)}_\cG)^\ell (x,y)}{\lambda_{(\omega,\beta)}^\ell} \, \log( \frac{\psi^{(\omega,\beta)}(x) \eta^{(\omega,\beta)}(y) e^{-\beta \bE(\gamma)}}{\lambda_{(\omega,\beta)}^\ell}) $$
 using as in Proposition~\ref{MERWprop} the fact that along the edges $e_i$ of the path $\gamma$
 each $\cK_\cG(x_i,x_{i+1})=1$, so that we obtain 
 $$ -\beta^{-1}  {\rm Sh}(\P_{(\beta,\omega),\ell}) = \beta^{-1} \sum_{x,y} \psi^{(\omega,\beta)}(x) \eta^{(\omega,\beta)}(y)
 \frac{(\cK^{(\omega,\beta)}_\cG)^\ell (x,y)}{\lambda_{(\omega,\beta)}^\ell} \, ( -\beta \bE(\gamma) + 
 \log ( \frac{\psi^{(\omega,\beta)}(x) \eta^{(\omega,\beta)}(y) }{\lambda_{(\omega,\beta)}^\ell } )  )   \, ,  $$
 while
 $$  \sum_\gamma \P_{(\beta,\omega),\ell}(\gamma)\, \bE(\gamma) =
  \sum_{x,y} \psi^{(\omega,\beta)}(x) \eta^{(\omega,\beta)}(y)
   \frac{(\cK^{(\omega,\beta)}_\cG)^\ell (x,y)}{\lambda_{(\omega,\beta)}^\ell} \, \bE(\gamma) \, . 
 $$
 This gives
 $$ \bF_{\bP_{(\omega,\beta),\ell}} = \beta^{-1} \sum_{x,y} \psi^{(\omega,\beta)}(x) \eta^{(\omega,\beta)}(y)
 \frac{(\cK^{(\omega,\beta)}_\cG)^\ell (x,y)}{\lambda_{(\omega,\beta)}^\ell} 
  \log ( \frac{\psi^{(\omega,\beta)}(x) \eta^{(\omega,\beta)}(y) }{\lambda_{(\omega,\beta)}^\ell } ) \, . $$
  $$ =  \beta^{-1} ( \sum_x \psi^{(\omega,\beta)}(x) \eta^{(\omega,\beta)}(x) \log (\psi^{(\omega,\beta)}(x)) 
  + \sum_y \psi^{(\omega,\beta)}(y) \eta^{(\omega,\beta)}(y) \log ( \eta^{(\omega,\beta)}(y)) 
  -  \beta^{-1}  \log( \lambda_{(\omega,\beta)}^\ell) $$
  $$ = -  \beta^{-1}  ({\rm Sh}(\pi^{(\omega,\beta)})   + \ell \log \lambda_{(\omega,\beta)}) \, , $$
so that the limit of $\bF_{\bP_{(\omega,\beta),\ell}}/\ell$ gives the optimal value 
$\min_\cS \bF(\cS)=-\beta^{-1} \log(\lambda_{(\omega,\beta)})$.
\endproof

\subsection{Linguistic implications of entropy optimization}\label{LingEmtropySec} 

In the case of the Merge dynamics on $\cG_{n,A}$, we consider the matrix
representation $\hat\cK^{(A,n)}$ of the Merge operations in the appropriate
basis in which it is a Hopf algebra Markov chain, as explained above. As we have
seen in this section, the form $\hat\cK^{(A,n)}$ of the Merge dynamics is
exactly the MERW on the Merge graph $\cG_{n,A}$.

\smallskip

Thus, this shows a particularly interesting property of the action of Merge
on workspaces, modeled as a Hopf algebra Markov chain: Merge is {\em information
maximizing}, in the sense that it maximizes the entropy rate ${\rm Sh}(\cS)$ over all
the possible Markov chains on the same directed graph. 

\smallskip

As we showed in this section, 
the Hopf algebra Markov chain $\hat\cK^{(A,n)}$ induces the {\em
uniform distribution} on the set of directed paths in the graph $\cG_{n,A}$.
As in Proposition~\ref{MERWprop}, we know that the probability
distribution on paths $\gamma$ of length $\ell$ between workspaces $F$ and $F'$ is
equal to
$$ \frac{\psi(F) \eta(F')}{\lambda^\ell}\, . $$
What this means is that, for example, if we start a Merge derivation with given
lexical material (the initial workspace $F$) and we obtain a resulting
workspace $F'$ (say, a fully formed sentence $F'=T$ consisting of a single syntactic
objects), all viable derivations of the same length (directed path of edges in 
the graph $\cG_{n,A}$ from $F$ to $F'$) have the same probability in the
Markov chain. This is a stronger Markovian property than just the usual
observation that ``Merge is Markovian". It means that the history of the
derivation does not matter and only the starting and final point and the
length determine the likelihood of that derivation, but not in any way the 
specific history. 

\smallskip

Indeed, the fact that such a stronger Markovian property should hold for Merge 
was already observed by Chomsky in \cite{ChomskyGK}, p.20, where he talks about
``derivations strictly Markovian in a strong sense, beyond the normal 
Markovian property". This property is described as the fact that in the
process of  structure formation in syntax the current state of the
derivation does not contain the history of the derivation. We interpret
this here as the property that the probability of a derivation (a path) is
only dependent on its result and not on the path itself (the history).

\smallskip

We will discuss the effects of Minimal Yield and other optimality constraints
in \S \ref{CostSec}, but we want to first point out here that the uniform
distribution on paths, that the dynamics $\hat\cK^{(A,n)}$ satisfies, already
ensures that such a stricter Markovian property holds, according to which the 
history of the derivation has no effect on the dynamics (except for its length, 
which is inevitable, since minimal length derivations ought to be favored over 
longer derivations, for optimality and for convergence). 

\smallskip

There is another related optimization property involving the Shannon entropy,
that will provide an important additional piece of information when evaluating
the cost functions that are usually considered in linguistics (Minimal Search and
Resource Restrictions) and their effect on the dynamics. We will show that it is only
when one also incorporates an information optimizing property that the cost
functions have the desired effect of dampening the Sideward Merge contributions
to the dynamics and ensure the convergence of the structure formation process
of syntax to the most connected workspaces (single trees) via External Merge. 
We will return to discuss this in \S \ref{EntCostSec}.

\smallskip

The optimization of free energy instead of entropy that we discussed in this
section will be useful when we discuss the Merge dynamics with
the optimization of cost functions, where the entropy-maximizing Markov 
chain $\hat\cK^{(A,n)}$ induces a free-energy minimizing Markov chain.

\section{Decomposing the Hopf algebra Markov chain} \label{DecompSec}

We proceed to analyze the dynamics of the Merge Hopf algebra Markov chain by
separately analyzing different parts of the Merge graph and different Merge operations.
In the linguistic model, External Merge is responsible for structure formation (including
theta role assignments, to which we will return in \S \ref{ColorSec}), while Internal Merge 
is responsible for transformation via movement. Sideward Merge is invoked to justify 
some special linguistic phenomena (head-to-head movement, for example) as
discussed in \cite{MLH}. Its existence and necessity as part of the Minimalism is
controversial. We will discuss in detail its role in the dynamics of Merge and
what changes if only EM and IM are considered, \S \ref{IMEMsec}.

The way we proceed in analyzing different parts of the dynamics is by separating out
the set $V(\cG_{n,A})$ of vertices $F$ of the graph $\cG_{n,A}$ into subsets $V_{\wp,n,A}$ 
labelled by a partition $\wp\in \cP'(n)$ of the integer $n$, with $\cP'(n)$ the subset of partitions as in \eqref{Pprimen}. 
We then consider, for each $\wp\in \cP'(n)$, the subgraphs $\cG_{n,A,\wp}$ on 
the set of vertices $V_{\wp,n,A}$ with only the edges given by IM arrows. 
We study, separately, the dynamics within
each graph $\cG_{n,A,\wp}$. We also consider a graph $\cG_{\cP'(n)}$ that has as
vertices the partitions $\wp \in \cP'(n)$ and a directed edge $\wp \to \wp'$ whenever
there are directed edges in $\cG_{n,A}$ between some vertex in $V_{\wp,n,A}$ 
and some vertex in $V_{\wp',n,A}$. The dynamics on $\cG_{\cP'(n)}$ will provide 
the remaining information needed to understand how the dynamical systems on the graphs
$\cG_{n,A,\wp}$ are intertwined in the full Merge dynamics over the graph $\cG_{n,A}$. We will
see that this strategy separates out in a natural way the different linguistic forms of Merge and
explains the effects on the dynamics of their interaction. 

\begin{prop}\label{GpnA}
Let $\fT_{\cS\cO_0,k}$ denote the set of nonplanar full binary rooted trees with $k$ leaves
decorated by elements of $\cS\cO_0$.
Let $\cP(n)$ be the set of all partitions and $\cP'(n)=\cP(n)\smallsetminus \{ 1,1,\ldots, 1 \}$ as in \eqref{Pprimen}. 
Consider the decomposition \eqref{VGpnA} of the set of vertices
$V(\cG_{n,A}) = \sqcup_{\wp\in \cP'(n)} V_{\wp,n,A}$.
Let $\cG_{n,A,\wp}$ denote the subgraph of $\cG_{n,A}$
with vertex set $V_{\wp,n,A}$ and only the edges consisting of IM arrows. 
The EM arrows map vertices of $\cG_{n,A,\wp}$ to vertices of $\cG_{n,A,\wp'}$ where $\wp'$ is a
a partition of the form $\wp'=\{ k_i+k_j \}\cup \{ k_\ell \}_{\ell\neq i,j}$, for some $i\neq j\in \{1,\ldots, r\}$.
If the partition $\wp$ contains a $2$ and a $1$, then there are some minimal SM arrows that map the set
$V_{\wp,n,A}$ to itself. All other minimal SM arrows map vertices of $\cG_{n,A,\wp}$ 
to vertices of $\cG_{n,A,\wp'}$ where 
$\wp'\neq \wp$ is a
a partition of the form $\wp'=\{ 2, k_i-1, k_j -1 \}\cup \{ k_\ell \}_{\ell\neq i,j}$ or of the form
 $\wp'=\{ 2, k_i-2 \}\cup \{ k_\ell \}_{\ell\neq i}$, for some $i\neq j \in \{1,\ldots, r\}$ or $i \in \{1,\ldots, r\}$, respectively. 
 Thus, in all of these cases $\cG_{n,A,\wp}$ is the induced subgraph on the set of vertices $V_{\wp,n,A}$.
\end{prop}

\proof Given the decomposition of the set of vertices as in Lemma~\ref{VertPart}, we obtain a map
$p: \fF_{A,n} \to \cP'(n)$, that assigns to each $F$ the corresponding partition of the $n$ leaves
into the leaves of each connected component. The range of this map consists of $\cP'(n)$,
where the partition $n=1+1+\cdots +1$ is not counted because, as previously discussed, 
we have removed the forest consisting of only
leaves with no edges from the dynamics, since it is always a transient state. The set of vertices $V_{\wp,n,A}$ is the preimage
$V_{\wp,n,A}=p^{-1}(\wp)$ under this map. At a vertex $F$ of the graph $\cG_{n,A}$, each IM arrow acts on 
one of connected component $T_i$ of the forest $F$ (acting as the identity on the other components), 
preserving the number of leaves $k_i=\# L(T_i)$ of each component. Thus, the IM transformations preserve the
partition, namely we have $p\circ {\rm IM} = p$, so that the sets $V_{\wp,n,A}$ are preserved under IM, or in other
words the IM arrows are edges of the graphs $\cG_{n,A,\wp}$. If the connected components of the workspace 
do not include a union $\fM(\alpha_i, \alpha_j) \sqcup \alpha_k$ of a cherry tree and a single leaf (the partition $\wp$ 
does not contain a $2$ and a $1$), the graph $\cG_{n,A,\wp}$ agrees with
the induced subgraph of $\cG_{n,A}$ on the same set of vertices $V_{\wp,n,A}$. Indeed, in such cases
only the IM arrows can preserve $\wp$. EM arrows always decrease the number of connected components
by one by grafting two connected components to a common root,
modifying the partition $\wp=\{ k_1, \cdots, k_r\}$ of $n$ to a partition of $n$ where two terms 
$k_i,k_j$ are replaced by a term $k_i+k_j$.  
The only minimal SM arrows that can preserve $\wp$ are those that act on two components
given by a cherry and a single leaf as either of
$$ \Tree[ $\alpha_i$ $\alpha_j$ ] \sqcup \alpha_k \mapsto \Tree[ $\alpha_i$ $\alpha_k$ ] \sqcup \alpha_j $$
$$ \Tree[ $\alpha_i$ $\alpha_j$ ] \sqcup \alpha_k \mapsto \Tree[ $\alpha_j$ $\alpha_k$ ] \sqcup \alpha_i \, . $$
All the minimal SM arrows produce a new cherry tree, so a term $2$ is added to the partition $\wp$, so unless
a previously single leaf component is eliminated by merging with another leaf into the new cherry, the number
of connected compoments would decrease, with either two terms $k_i,k_j$ are 
decreased to $k_i-1, k_j-1$ (including the case where either $k_i$ or $k_j$, but not both, is equal to $1$)
or a single term $k_i\geq 3$ is decreased to $k_i-2$. 
Thus, only in the cases where the components of the forests $F\in V_{\wp,n,A}$ include a cherry and a single leaf,
the IM graph $\cG_{n,A,\wp}$ is smaller than the induced subgraph with vertex set $V_{\wp,n,A}$. 
On the other hand, all the IM arrows will occur as edges of some $\cG_{n,A,\wp}$. 
\endproof

This decomposition in particular means that, if we want to study the Internal Merge dynamics, we can
consider a single graph $\cG_{n,A,\wp}$ at a time. 

\section{The Internal Merge dynamics}\label{IMdynSec}

We discuss here the long term behavior of the Internal Merge dynamics and convergence to
a stationary distribution.

\begin{defn}\label{IMgraph}
Let $\cG_{n,A}^{\rm IM}$ be the graph with vertex set $\fF_{A,n}$ and with edges the IM
arrows (namely the graph $\cG_{n,A}$ with the EM and SM arrows removed). 
We write $\cG^{\rm IM}_{n,A,\wp}$  for the induced subgraph of  $\cG_{n,A}^{\rm IM}$ 
with vertex set $V_{\wp, n, A}$.
\end{defn}

\smallskip

\subsection{Connected components}\label{connIMsec}

We first consider the Internal Merge dynamics on workspaces consisting of
a single tree.

\begin{prop}\label{connIM}
The graph $\cG_{n,A,\{ n \}}^{\rm IM}$ with $\wp=\{ n \}$, with vertices the
workspaces $F=T\in \fF_{A,n}$ consisting of a single tree, 
is connected.
\end{prop}

\proof The graph has $(2n-3)!!$ vertices corresponding to all the possible
tree topologies of $T$ (as a non-planar tree with labelled leaves).
We need to show that all of them are in the same connected component.
A full binary rooted tree $T$ is a {\em comb tree} iff it 
contains only one cherry. In that case, if IM extracts any one of the leaves, it produces
another comb tree with a permutation of the leaves labels,
and in fact all permutations of the labels occur in this way. So all the
comb trees are in the same connected component. Thus, 
to show that $\cG_{n,A,\{ n \}}^{\rm IM}$ is connected, it suffices to show that, 
for any $T$ with $n$-leaves labelled by
the elements of $A$, that is not already a comb tree, there is a path of IM arrows,
each extracting a single leaf, that connects $T$ to some comb tree on $n$-leaves, 
with some labeling by the elements of $A$. 
We can do this inductively, starting with $n=4$ (for $n=3$ there are only comb trees),
where one IM on one of the leaves of the non-comb topology produces a comb tree.
Suppose that for $m<n$ leaves the statement holds. 
If $T$ is not a comb tree then there is more than one cherry. Apply an IM that extracts one
of the two leaves, call it $\alpha$, in one of the cherries. The resulting tree has that leaf $\alpha$ 
attached directly to the root, with the vertex $v$ that was the root of $T$ as its sister vertex,
which is now the root vertex of $T/^d \alpha$. The tree $T/^d \alpha$ has $n-1$ leaves,
so there is a path of IM arrows extracting one leaf at a time that connect $T/^d \alpha$
to a comb. But then performing the same path of IM operations on the same sequence
of leaves, on the tree $\fM(\alpha, T/^d \alpha)$ also produces a comb tree, because
the comb sequence of leaves extracted and merged at the root continue a comb sequence
already starting with $\alpha$.
\endproof

\smallskip

Note that this argument proves connectedness but not strong connectedness. In fact
the graph $\cG_{n,A,\{ n \}}^{\rm IM}$ is also strongly connected: this will follow from the
more general argument given in Proposition~\ref{dynIMunif}.

\smallskip

\begin{prop}\label{dynIMunif1}
For a partition $\wp\in \cP'(n)$ of the form $\wp=\{ k_1, \ldots, k_1, \ldots, k_r, \ldots k_r \}$ with $n=a_1 k_1+\cdots + a_r k_r$,
the graph $\cG^{\rm IM}_{n,A,\wp}$ consists of separate connected components. The number of
connected components is equal to the generalized multinomial coefficient of \eqref{genmultinomial},
$$ b_0(\cG^{\rm IM}_{n,A,\wp}) =\Upsilon_{\wp,n}= \frac{1}{a_1 ! \cdots a_r !} \binom{n}{\underbrace{k_1,\ldots, k_1}_{a_1\text{-times}} \ldots, 
\underbrace{k_r, \ldots, k_r}_{a_r\text{-times}}}  \, . $$
We denote by $C_\wp$ the set of connected components of the graph $\cG^{\rm IM}_{n,A,\wp}$ and
we write $\cG^{\rm IM}_{n,A,\wp,\sigma}$ for the component $\sigma\in C_\wp$. The set $C_\wp$
can be obtained as a quotient of the conjugacy class ${\rm Conj}_\wp$ in the symmetric group $S_n$ of cycle type $\wp$, 
with the relation
$$ \Upsilon_{\wp,n} \cdot \prod_{i=1}^r ((k_i-1)!)^{a_i} = \# {\rm Conj}_\wp =\frac{n!}{\prod_{i=1}^r a_i ! \,\,  k_i^{a_j}} \, . $$
Each connected component $\cG^{\rm IM}_{n,A,\wp,\sigma}$ contains a number of vertices equal to
$$ \prod_{i=1}^r ((2k_i-3) !!)^{a_i} \,\,  . $$
\end{prop}

\proof Internal Merge acting on a workspace $F=T_1\sqcup\cdots\sqcup T_r$
acts separately on one connected component  $T_i$ of the workspace at a time,
so that the action on workspaces is completely determined by the action on trees.
In particular, the action preserves the splitting $L(F)=\sqcup_i L(T_i)$
by preserving each set $L(T_i)=L(\fM(T_v, T_i/T_v))$. This means that
a forest $F\in V_{\wp,n,A}$ will always be mapped by IM arrows to other
forests in $V_{\wp,n,A}$ with the same $\wp$. The partition of the set of labels 
over the subsets $L(T_i)$ of leaves is also preserved by IM (unlike SM
and EM): suppose 
a given $F\in V_{\wp,n,A}$ for $\wp=\{ k_1, \ldots, k_r \}$ has an assigned 
partition of the labels set $A=\{ \alpha_1, \ldots, \alpha_n \}$ into subsets of sizes $k_i$ assigned to the
sets of leaves $L(T_i)$ of the components $T_i$ of $F$. Then Internal Merge will preserve this partition of the
labels and cannot mix leaf labels assigned to different components of the workspace. Thus, these partitions
of the set $A$ are invariant along each directed path of Internal Merge arrows on the set $V_{\wp,n,A}$.
This means that each such orbit, determined by a label assignment, is a different connected 
component of the graph $\cG_{n,A,\wp}$. 

If the partition $\wp$ is of the form $n=k_1+\cdots+ k_r$ with all the $k_i\leq 2$,
there are no Internal Merge arrows (we exclude the IM arrows that extract one 
of the accessible terms at the vertices immediately below the root since those 
are just the identity, so both the single leaves and the cherry trees
do not have IM dynamics). So in this case each vertex of $\cG^{\rm IM}_{n,A,\wp}$ 
is a connected component. 

When the partition $\wp$ has at least one of the $k_i$ satisfies $k_i\geq 3$, there are nontrivial IM
arrows in $\cG^{\rm IM}_{n,A,\wp}$. For any partition of the labels set $A$ over the components
of the workspaces, the $\prod_i (2k_i -3)!!$ different tree topologies are in the same connected
component, by Proposition~\ref{connIM}.

Thus, the connected components of $\cG^{\rm IM}_{n,A,\wp}$ are in bijective correspondence with
the partitions of the set $A$ of labels over components of the workspace, 
so that there are $\Upsilon_{\wp,n}$ components. 

Consider permutations $\sigma\in S_n$ that have the cycle
structure determined by the partition $\wp$, namely a number $a_i$ of $k_i$-cycles.
The cycle type $\wp$ determines a conjugacy class ${\rm Conj}_\wp$  in the symmetric group $S_n$.
When counting elements in the conjugacy class ${\rm Conj}_\wp$ of cycle type $\wp$,
a factor equal to the multinomial coefficient $\mu_{\wp,n}$ as in \eqref{multinomial} counts the choices of 
elements that to go in each cycle. Each of these choices generates $k_i !/k_i =(k_i-1)!$ assignments of labels
counted modulo cyclic permutations. This counting is taken up to permutations of the $a_i$ cycles of size $k_i$, resulting
in $\# {\rm Conj}_\wp=\Upsilon_{\wp,n} \cdot \prod_i ((k_i-1)!)^{a_i}= \#C_\wp \cdot \prod_i ((k_i-1)!)^{a_i}$. 

In each component $\cG^{\rm IM}_{n,A,\wp, \sigma}$, one has a choice of all
the possible tree topologies of the trees $T_i$ in the workspace $F=T_1 \sqcup \cdots \sqcup T_r$,
which gives a product of $((2k_i -3)!!)^{a_i}$ for each component $T_i$, which determines the
number of vertices in each connected component.
\endproof

\subsection{Strong connectedness of components}\label{sconnIMsec}

We can then prove a stronger result about the Internal Merge dynamics on 
each connected component $\cG^{\rm IM}_{n,A,\wp,\sigma}$, $\sigma\in C_\wp$.

\begin{prop}\label{dynIMconv}
Let $\wp\in \cP'(n)$ be a partition $\wp=\{ k_1,\ldots, k_r \}$ with at least one of the $k_i\geq 3$, so that
there are non-trivial IM arrows. Then 
each connected component  $\cG^{\rm IM}_{n,A,\wp, \sigma}$ of the graph $\cG^{\rm IM}_{n,A,\wp}$ 
is strongly connected and aperiodic, hence 
the Internal Merge Hopf algebra Markov chain $\hat\cK^{(A,n,\wp,\sigma)}$ is ergodic on each connected
component.
\end{prop}

\proof As we mentioned in the case of Proposition~\ref{ErgodicK}, 
ergodicity follows from strong connectedness and aperiodicity. To
show that the Hopf algebra Markov chain $\hat\cK^{(A,n,\wp,\sigma)}$ of Internal Merge is ergodic
we need a refinement of Proposition~\ref{StrongConnAper} to show that not only the graphs
$\cG_{n,A}$ are strongly connected and aperiodic but also each connected component 
$\cG^{\rm IM}_{n,A,\wp,\sigma}$ of the subgraphs $\cG^{\rm IM}_{n,A,\wp}$ has the same properties.

For a given component $\cG^{\rm IM}_{n,A,\wp, \sigma}$, we look at the restriction $\hat\cK^{(A,n,\wp,\sigma)}$
of the dynamics of $\hat\cK^{(A,n,\wp)}$ to that component. 
By Proposition~\ref{EdgesGnA} we know that there are, at each $F\in V_{\wp,n,A,\sigma}$
as many outgoing as incoming IM edges,
\begin{equation}\label{IMinoutd}
 N_{\rm IM}^{\rm out}(F) = \sum_{i=1}^{c(F)} (2k_i-4) = N_{\rm IM}^{\rm in}(F) \, , 
\end{equation} 
where $F=T_1 \sqcup \cdots \sqcup T_r$ with $\wp$ the partition $n=k_1+\cdots+ k_r$
for $k_i=\# L(T_i)$, and where we assume for simplicity of notation that $T_1,\ldots, T_{c(F)}$ 
are the components with nonempty set of edges, with $c(F)\leq r$. The numbers $k_i$ as 
well as $c(F)=\# \{ i\,|\, k_i >1 \}$ are determined by the partition $\wp$ so they are the same
for all vertices $F\in V_{\wp,n,A}$ for fixed $\wp$.

The strong connectedness property then 
follows from the fact that each vertex of $\cG_{n,A,\wp,\sigma}$ has $\deg^{\rm out}(F)=\deg^{\rm in}(F)$, which
implies that the directed graph $\cG_{n,A,\wp,\sigma}$ is Eulerian, namely there is a directed
cycle that visits each edge of $\cG_{n,A,\wp,\sigma}$ exactly once. This then implies strong connectedness
as, starting at any vertex $F$ one can reach any other vertex $F'$ along the Eulerian circuit. We can see
that $\cG_{n,A,\wp,\sigma}$ must also be aperiodic by refining the argument given in 
Proposition~\ref{StrongConnAper} that identifies explicit cycles of Internal Merge operations. 
It suffices to consider the case of $\wp=\{ n \}$ namely the workspaces $F=T$ consisting of a single
connected component $T\in \fT_{\cS\cO_0,n}$. The graph $\cG_{n,A,\wp=\{ n \}}$ has a
single connected component. If we show that these graphs
are aperiodic, then all the other graphs $\cG_{n,A,\wp,\sigma}$ will also be aperiodic because 
for every $k_i\geq 3$ in the partition $\wp=\{ k_1, \ldots, k_r\}$ they will contain copies of the cycles 
we exhibit in $\cG_{n,A, \{ k_i \}}$. 
We consider separately the cases $n=3$ and $n=4$ leaves and the general case $n\geq 5$ leaves.
We show that in all of these cases we always have an IM cycle of order $2$ and an IM cycle of order $3$
so that the gcd of the cycle-lengths is equal to $1$.
For $n=3$, by extracting one of the other leaf of the cherry, we obtain cycles 
$$ \Tree[ $\alpha$ [ $\beta$ $\gamma$ ] ] \mapsto \Tree[ $\beta$  [ $\alpha$ $\gamma$ ] ] \mapsto \Tree[ $\alpha$ [ $\beta$ $\gamma$ ] ]  $$
$$ \Tree[ $\alpha$ [ $\beta$ $\gamma$ ] ] \mapsto \Tree[ $\beta$  [ $\alpha$ $\gamma$ ] ] \mapsto 
\Tree[ $\gamma$ [ $\beta$ $\alpha$ ] ] 
\mapsto  \Tree[ $\alpha$ [ $\beta$ $\gamma$ ] ]  \, . $$
In the case with $n=4$ we similarly have a cycle of length $3$
$$ \Tree[ $\alpha_1$ [ $\alpha_2$ [ $\alpha_3$ $\alpha_4$ ] ] ] \mapsto \Tree[ $\alpha_3$ [ $\alpha_1$ [ $\alpha_2$ $\alpha_4$ ] ] ]
\mapsto \Tree[ $\alpha_2$ [ $\alpha_3$ [ $\alpha_1$ $\alpha_4$ ] ] ]\mapsto \Tree[ $\alpha_1$ [ $\alpha_2$ [ $\alpha_3$ $\alpha_4$ ] ]  ] $$
and also a cycle of length $2$ (by first extracting the $\{ \alpha_3, \alpha_4 \}$ cherry and then the single leaf $\alpha_1$)
$$ \Tree[ $\alpha_1$ [ $\alpha_2$ [ $\alpha_3$ $\alpha_4$ ] ] ] \mapsto \Tree[ [ $\alpha_3$ $\alpha_4$ ] [ $\alpha_1$ $\alpha_2$ ] ]  \mapsto  \Tree[ $\alpha_1$ [ $\alpha_2$ [ $\alpha_3$ $\alpha_4$ ] ] ]  \, . $$
Thus $\cG_{2,A \{ 2 \}}$ and $\cG_{3,A \{ 3 \}}$ are aperiodic. Now suppose that $n\geq 5$. Consider trees of the form
$$ T_{r,k,n-rk}=\Tree[ $T_1$ [ $T_2$ [ $\cdots$ [ $T_r$ $T_{r+1}$ ]] ] ] =  \Tree[ \qroof{$k$-leaves}. [ \qroof{$k$-leaves}. [ $\cdots$ [ \qroof{$k$-leaves}. \qroof{$\ell$-leaves}. ] ] ] ] $$
with subtrees $T_i$ for $i=1,\ldots, r$ with $k$-leaves with identical tree topologies, 
and an additional subtree $T_{r+1}$ with $\ell$ leaves, for $n=rk+\ell$ and $\ell\geq 1$, $r\geq 2$ and $k\geq 1$. 
The trees considered in Proposition~\ref{StrongConnAper} are the special case $T_{n-1,1,1}$. Then the iteration
of IM operations that extract the accessible terms $T_r$ then $T_{r-1}$ and so on until $T_1$ gives a cycle of length $r$.
For $n\geq 5$ these cycles always include a cycle of length $r=2$ and a cycle of length $3$ so again the gcd of the cycle 
lengths is equal to one, so aperiodicity always holds.
\endproof

\subsection{Uniform distribution}\label{unifdistrIMsec} 
The same property \eqref{IMinoutd} of the in/out-degrees of the Internal Merge dynamics
also give us the form of the stationary distribution on the $\cG^{\rm IM}_{n,A,\wp,\sigma}$ graphs.

\begin{prop}\label{dynIMunif}
For partitions with at least one $k_i\geq 3$, the connected component 
$\cG^{\rm IM}_{n,A,\wp, \sigma}$ of the graph $\cG^{\rm IM}_{n,A,\wp}$ 
are strongly connected and for any $\sigma\in C_\wp$,
the part $\hat\cK^{(A,n,\wp,\sigma)}$ of the Hopf algebra Markov chain given by the action 
of Internal Merge on $\cG^{\rm IM}_{n,A,\wp, \sigma}$
has the uniform distribution as stationary distribution.  When all $k_i\leq 2$  
there are no IM arrows and each vertex is a component.
\end{prop}

\proof As shown in Propositions~\ref{EdgesGnA} and \ref{dynIMconv}, 
at each vertex $F\in V_{\wp,n,A,\sigma}$ there are as many outgoing as incoming IM edges,
as in \eqref{IMinoutd}, $\deg^{\rm out}(F)=\deg^{\rm in}(F)$.   
Note that the value in \eqref{IMinoutd}
is constant for all vertices $F\in V_{\wp,n,A}$, and depends only on the partition $\wp$,
$$ N_{\rm IM}^{\rm out}(F) = N_{\rm IM}^{\rm in}(F)  =: d_\wp \, . $$
Let $\cK^{(A,n,\wp,\sigma)}=(\cK^{(A,n,\wp,\sigma)}_{F,F'})$ denote the adjacency matrix of the graph 
$\cG^{|rm IM}_{n,A,\wp,\sigma}$.

For a directed graph where all the vertices have the same in and out degree, all of them equal to $d_\wp$,
the Perron-Frobenius theorem, applied to the adjacency matrix, shows that the Perron-Frobenius eigenvalue satisfies
$$ \min_F \sum_{F'} \cK^{(A,n,\wp,\sigma)}_{F,F'} \leq \lambda \leq \max_F \sum_{F'} \cK^{(A,n,\wp,\sigma)}_{F,F'} \, , $$
where $\sum_{F'} \cK^{(A,n,\wp,\sigma)}_{F,F'} =\deg^{\rm out}(F)$, so if $\deg^{\rm out}(F)=d_\wp$ for all 
$F\in V_{\wp,n,A,\sigma}$,
we have that $\lambda=d_\wp$ is the Perron-Frobenius eigenvector and correspondingly the uniform
distribution $\eta(F)=1$ is the Perron-Frobenius eigenvalue 
$$ \sum_{F'} \cK^{(A,n,\wp,\sigma)}_{F,F'} \eta(F')= d_\wp \, \eta(F)\, . $$
Thus, the part of the Hopf algebra Markov chain given by Internal Merge simply has transition matrix
$$ \hat\cK^{(A,n,\wp,\sigma)}_{F,F'}= d_\wp^{-1}\cdot  \cK^{(A,n,\wp)}_{F,F'}\, , $$
with no rescaling of the basis, since $\eta(F)=1$ for all $F\in V_{\wp,n,A}$. This Markov chain
has stationary distribution
$$ \sum_F \pi(F) \hat\cK^{(A,n,\wp,\sigma)}_{F,F'} =\pi(F') $$
given by the uniform distribution $\pi(F)=1$ (normalized to $\pi(F)=1/\# V_{n,A,\wp}$), since
$$ d_\wp^{-1} \sum_F \cK^{(A,n,\wp,\sigma)}_{F,F'} = \frac{\deg^-(F')}{d_\wp} =1\, , $$
since $\deg^+(F)=\deg^-(F)=d_\wp$ for all $F\in V_{\wp,n,A,\sigma}$. 
\endproof

\smallskip

\smallskip

One can rephrase Propositions~\ref{dynIMunif} and \ref{dynIMconv} in the following way.

\begin{cor}\label{KIMconverge}
The uniform distribution is the unique stationary distribution of the Internal Merge 
Hopf algebra Markov chains $\hat\cK^{(A,n,\wp,\sigma)}$, on each connected
component $\cG_{n,A,\wp,\sigma}$ of $\cG_{n,A,\wp}$, 
and any initial distribution on the set of workspaces $V_{\wp,n,A, \sigma}$ converges
under iterations of $\hat\cK^{(A,n,\wp, \sigma)}$ to the uniform distribution. 
\end{cor}

\smallskip

This shows that Internal Merge by itself is an ergodic dynamical system with uniform stationary distribution.
This fact singles out Internal Merge as having especially nice dynamical properties.

\section{Internal and External Merge dynamics}\label{IMEMsec}

In this and the next section we outline the opposite roles that the minimal Sideward Merge and the External Merge 
arrows plays in the dynamics by showing what happens if one either removes all the SM arrows, while keeping 
IM and EM arrows or (in the next section) if one removes the EM arrows and keeps IM and SM arrows. In both
cases the strong connectedness property breaks down and the vertices of the graph split into transient and
recurrent states, with the stationary distribution supported only on the recurrent states. The nature of the
recurrent states is complementary in the two cases: predictably, the dynamics moves to the most connected
structure when no SM is present and to the least connected ones when no EM is present. 

\smallskip

The dynamical system involving only External and Internal Merge is the one that is usually considered in the
linguistics literature, where it is argued (see \cite{ChomskyUCLA}, \cite{ChomskyGK}, \cite{ChomskyElements})
that the Sideward Merge part of the dynamics is excluded on the basis of optimization with respect to certain
cost functions, that we will discuss in \S \ref{CostSec}. However, from the mathematical perspective, as well as
in comparison with similar settings in physics models, one can argue (see \cite{MLH}) that cost
optimization is not a hard structural constraint that eliminates SM, but rather a soft constraint that makes it
more unlikely, or contributing less significantly, the more it deviates from optimality. Using this viewpoint, it
was shown in \cite{MLH} that certain linguistic phenomena can be explained by a minimally optimality
violating form of SM, which is exactly the minimal SM we are considering in the graphs $\cG_{n,A}$ here.

\smallskip

In this section, we present the Merge dynamics in the way it would be according to the original
proposal of \cite{ChomskyUCLA}, \cite{ChomskyGK}, \cite{ChomskyElements}, where the SM
arrows of the graph $\cG_{n,A}$ are completely eliminated. We show that, while in this case the
dynamics does indeed converge through EM to the connected structures (workspaces consisting
of a single tree) with the remaining action of IM on trees as the only recurrent part of the dynamics, 
the strong connectedness property is broken, eliminating the desirable dynamical property or
ergodicity and the independence of the dynamics on the initial position, as well as the entropy
optimizing property we discussed earlier that also relies on strong connectedness.   
We will return in \S \ref{CostSec} and \S \ref{EntCostSec} to discuss how to use the type of cost
functions proposed in \cite{ChomskyUCLA}, \cite{ChomskyGK}, \cite{ChomskyElements}
and in \cite{MCB} and the entropy optimization property to correct the dynamics by 
incorporating cost-weights that suppress the undesirable effects of Sideward Merge without
losing the desirable strong connectedness and ergodicity and entropy optimization
properties of the dynamics.

\smallskip

We consider here the subgraph $\cG^{\rm IM-EM}_{n,A} \subset \cG_{n,A}$
with the same set of vertices $V(\cG^{\rm IM-EM}_{n,A})=V(\cG_{n,A})=\fF_{A,n}$, 
that only includes as edges the IM and the EM arrows but not the SM arrows of $\cG_{n,A}$.
It is then clear that the strong connectedness property fails for this subgraph $\cG^{\rm IM-EM}_{n,A}$.
We refer to the resulting Markov chain as $\hat\cK^{{\rm IM-EM},(A,n)}$.

\subsection{Reducibility and transient states}\label{redIMEMsec}
A transient state in a Markov chain is a state that has a nonzero probability of never 
being returned to. 

\begin{lem}\label{transient}
Every vertex $F\in V(\cG^{\rm IM-EM}_{n,A})$ given by a workspace with more than one connected component
(that is, with image $\wp=p(F)\in \cP'(n)$ with $\wp\neq \{ n \}$) is a transient state of
$\hat\cK^{{\rm IM-EM},(A,n)}$.
\end{lem}

\proof Since EM always decreases the number of connected components
of the workspace by one, it is clear that a state $F=T_1\sqcup \cdots \sqcup T_r\in \fF_{A,n}$ 
with $r>1$, namely a vertex in $\cG^{\rm IM-EM}_{n,A,\wp=\{ k_1, \ldots, k_r\}}$ has a non-zero probability of 
leaving the graph $\cG^{\rm IM-EM}_{n,A,\wp=\{ k_1, \ldots, k_r\}}$ via one of the EM arrows and entering
one of the graphs with $\cG^{\rm IM-EM}_{n,A,\wp'=\{ k_i+k_j, k_1, \ldots, \hat k_i, \ldots, \hat k_j, \ldots  k_r\}}$,
where the notation $\hat k_i$ means that the entry $k_i$ is removed.  Since without the SM arrows
every arrow in the graph $\cG^{\rm IM-EM}_{n,A}$ either maintains the same number of connected
components of the workspace (IM arrows) or decreases it by one (EM arrows), there is no way
to return from any vertex of $\cG^{\rm IM-EM}_{n,A,\wp'=\{ k_i+k_j, k_1, \ldots, \hat k_i, \ldots, \hat k_j, \ldots  k_r\}}$
to a vertex of $\cG^{\rm IM-EM}_{n,A,\wp=\{ k_1, \ldots, k_r\}}$, hence the vertices of 
$\cG^{\rm IM-EM}_{n,A,\wp=\{ k_1, \ldots, k_r\}}$
are transient. In the case of vertices of $\cG^{\rm IM-EM}_{n,A,\wp=\{ n \}}$ there are no outgoing EM arrows so
the vertex of this subgraph remain (via IM arrows) within the same $\cG^{\rm IM-EM}_{n,A,\wp=\{ n \}}$ and these
are the only recurrent states of the Markov chain $\hat\cK^{{\rm IM-EM},(A,n)}$.
\endproof

\begin{cor}\label{KIEreduce}
The Hopf algebra Markov chain $\hat\cK^{{\rm IM-EM},(A,n)}$ is reducible. Equivalently, the graph
$\cG^{\rm IM-EM}_{n,A}$ is not strongly connected.
\end{cor}

\proof
A finite irreducible Markov chain cannot have any transient state, hence Lemma~\ref{transient}
implies that $\hat\cK^{{\rm IM-EM},(A,n)}$ is reducible.
\endproof

\subsection{Communication classes}\label{commclIMEMsec}

A reducible Markov chain can be decomposed into {\em communication classes}: 
this is the same as the decomposition of the
underlying directed graph into {\em strongly connected components}. 

\smallskip

Given a directed
graph $G$, one says that two vertices $v,v' \in V(G)$ are strongly connected, or communicating,
if either $v=v'$ or if $v\neq v'$ and there is a directed path of edges from $v$ to $v'$ and a directed
path of edges from $v'$ to $v$. Strong connectivity so defined is an equivalence relation on
$V(G)$ and the equivalence classes are the strongly connected components of $G$ (or communication classes).

\smallskip

The {\em condensation graph} $\bar G$ is the directed graph obtained from $G$ by contracting each
strongly connected component of $G$ to a single vertex. The graph $\bar G$ is an {\em acyclic directed graph}.
In general $\bar G$ so obtained is a multigraph, with possibly multiple parallel edges between the same pair 
of vertices. We can further identify any set of parallel edges in $\bar G$ to a weighted single edge. 

\smallskip

A {\em closed communication class} is a strongly connected component with no arrows pointing to
other strongly connected components, namely such that its strong connectedness equivalence class
is a sink of the condensation directed acyclic graph.

\smallskip

Consider again the graphs $\cG^{\rm IM}_{n,A,\wp,\sigma}$, with $\sigma\in C_\wp$,
that are the  connected components of the IM graph $\cG^{\rm IM}_{n,A,\wp}$.

The following description of communication classes for $\hat\cK^{{\rm IM-EM},(A,n)}$ follows directly from
the description of the IM dynamics on the graphs $\cG^{\rm IM}_{n,A,\wp, \sigma}$ and the EM arrows.

\begin{prop}\label{SCCdecomp}
The graphs $\cG^{\rm IM}_{n,A,\wp, \sigma}$, for $\wp \in \cP'(n)$ and 
$\sigma\in C_\wp$, are the strongly connected components 
of the graph $\cG^{\rm IM-EM}_{n,A}$ (communication classes of the Markov chain $\hat\cK^{{\rm IM-EM},(A,n)}$).
The condensation graph $\bar \cG^{\rm IM-EM}_{n,A}$ is the graph $\tilde\cG^{\rm EM}_{\cP'(n)}$ 
with set of vertices 
$$ V(\tilde\cG^{\rm EM}_{\cP'(n)})=\{ (\wp,\sigma)\,|\, \wp \in \cP'(n)\, , \, \, \sigma\in C_\wp \} $$
and one directed edge from $(\wp,\sigma)$ to $(\wp',\sigma')$ for $$ \wp=\{ k_1, \ldots, k_r \} \ \  \text{ and }  \ \
\wp'=\{ k_i+k_j, k_1, \ldots, \hat k_i, \ldots, \hat k_j, \ldots  k_r\} $$ (with the notation as in Lemma~\ref{transient})
for some $k_i, k_j\in \wp$, and with $\sigma'$ the image of $\sigma$ under the induced map $C_\wp \to
C_{\wp'}$. 
The graph $\cG^{\rm IM}_{n,A,\wp}$ with $\wp=\{ n \}$, where vertices are workspaces consisting of 
a single tree with $n$ leaves, is the only closed communication class. 
\end{prop}

\smallskip

The induced map $C_\wp \to C_{\wp'}$ is simply determined by the fact that each assignment of
labels of $A$ to the boxes of $D(\wp)$ determines a corresponding assignment to the boxes of $D(\wp')$.

\subsection{Stationary distribution on the closed communication class}\label{stadistrIMEMsec}
This implies that, when we do not consider any Sideward Merge as part of the dynamics, the limiting
distribution for $\hat\cK^{{\rm IM-EM},(A,n)}$ is only supported on the connected structures. This describes
Merge as a process of structure formation where the dynamics evolves via External Merge toward 
connected structures and away from workspaces with multiple connected components, and once 
connected structures are reached the only dynamics left consists of the transformations via movement 
implemented by Internal Merge. This indeed agrees with the usual description of the Merge dynamics
given by Chomsky in \cite{ChomskyUCLA}, \cite{ChomskyGK}, \cite{ChomskyElements}, 
in terms of just EM and IM.

\begin{prop}\label{EMIMstationary}
The stationary distribution for the reducible Markov chain $\hat\cK^{{\rm IM-EM},(A,n)}$ is
the uniform distribution supported on the vertices of the graph $\cG^{\rm IM}_{n,A,\wp=\{ n \}}$ (workspaces
consisting of a single tree with $n$ leaves).
\end{prop}

\proof If a stationary distribution $\pi$ exists, it has the property that 
\begin{equation}\label{pilimPN}
 \pi(F) = \lim_{N \to \infty} \frac{1}{N} \sum_{i=1}^N (\hat\cK^{{\rm IM-EM},(A,n)})^N_{F,F'} \, \,  , 
\end{equation} 
for any initial state $F'$ in the same communication class. In the reducible case, $\pi(F)=0$ whenever $F$ is
a transient state of $\hat\cK^{{\rm IM-EM},(A,n)}$.  
Thus, $\pi$ can only be supported on the recurrent states, that in this case are the vertices
of the subgraph $\cG^{\rm IM}_{n,A,\wp=\{ n \}}$. Indeed, the stationary distribution for a reducible Markov chain is
always supported on the closed communication classes. In this case there is only one
such class, consisting of the vertices of $\cG^{\rm IM}_{n,A,\wp=\{ n \}}$. Moreover, on those vertices, $\pi$ should agree
with the stationary distribution of the Internal Merge dynamics on the graph $\cG^{\rm IM}_{n,A,\wp=\{ n \}}$
because choosing both $F,F'$ in \eqref{pilimPN} to be vertices in $\cG^{\rm IM}_{n,A,\wp=\{ n \}}$ leads to all
the nontrivial entries in the $(\hat\cK^{{\rm IM-EM},(A,n)})^N_{F,F'}$ being IM arrows so the limit has to
be the stationary distribution of the IM dynamics on $\cG^{\rm IM}_{n,A,\wp=\{ n \}}$, which is the
uniform distribution.  
\endproof

\subsection{The IM--EM structure formation model: linguistic aspects}\label{modelIMEMsec}

The dynamics described in this section and summarized in the result of Proposition~\ref{EMIMstationary}
gives, from the linguistic perspective, a very satisfactory view of structure formation in syntax,
obtained via free symmetric Merge, which shows convergence to the connected structures 
(fully formed sentences) and Internal Merge
dynamics acting on the connected structures in an ergodic way, without privileging any
particular subset of available structures (uniform stationary distribution). It indeed reflects what was described in Chomsky's
model of \cite{ChomskyUCLA} and \cite{ChomskyGK}, where Sideward Merge is excluded
on the ground of its failure of optimality. 

However, as we pointed out at the beginning of this section, this has the problem that elimination
of the Sideward Merge arrows looks artificial, in the sense that these Merge operations are 
clearly part of the natural Hopf algebra Markov chain definition of Merge, as compositions
$$ \fM_{S,S'}=\sqcup \circ (\cB \otimes {\rm id}) \circ \delta_{S,S'} \circ \Delta \, ,  $$
where non-primitive terms of the coproduct are used, so there is no {\em algebraic}
argument for their elimination, and cost counting only makes these Merge
operations subdominant with respect to EM and IM, but still present. 

As discussed in \cite{MCB}, one could replace the coproduct  $\Delta$ with its primitive part $\Delta^P$,
that only forms partitions of workspaces without extracting accessible terms, 
and then this same formula would provide only External Merge without Sideward Merge.
On the other hand this would not be consistent with maintaining the extraction of accessible
terms in Internal Merge and still being able to see Internal Merge and External Merge as cases
of the same operation (with an $\fM_{S,1}$ operation involved in IM). Thus, just removing the
SM arrows from the Hopf algebra Markov chain $\hat \cK$ of Merge is not directly justified by the
underlying algebraic structure. 

In Chomsky's formulation of Minimalism \cite{ChomskyUCLA}, \cite{ChomskyGK}, the selection of 
External Merge and Internal Merge over Sideward Merge is formulated as optimality 
(minimization of some natural cost function), where EM and IM are optimal
while SM is not. There are different types of optimality that are discussed, in this context, in 
the linguistics literature on Merge: Minimal Search and forms of Resource Restriction like Minimal Yield.
These are analyzed in terms of the mathematical model in \cite{MCB}. A detailed discussion of
these different cost functions and Sideward Merge is also further expanded in \cite{MLH}, where it
is shown more precisely what roles Minimal Search and Resource Restrictions have in
constraining possible forms of Sideward Merge. 

The way to deal with different weights by cost functions in our mathematical setting
is to replace the adjacency matrix of the Merge graph $\cG_{n,A}$ with a weighted
version, where the weights incorporate the cost functions, in a way similar to the
case we analyzed in \S \ref{FreeEnSec}, and then argue that this weighted form of
the Merge dynamics better approximates the EM--IM dynamics described in this section,
with the effects of SM dampened by the cost function. While this idea seems natural
and is easily confirmed in the simple example of $n=3$ leaves computed in \cite{MLH},
the general case of arbitrary $n$ is more subtle. We will return to discuss cost functions in \S \ref{CostSec}
and the effects on the Merge dynamics in \S \ref{CostSec} and \S \ref{EntCostSec}. 

\section{Internal and Sideward Merge dynamics} \label{IMSMsec}

We now consider, for comparison, what would happen if we only had IM and SM arrows, without EM.
This is of course a purely hypothetical scenario that does not correspond to any viable linguistic model,
since it clearly would not lead to structure formation in any way. It is useful, however, to discuss this
briefly, in preparation for the direct comparison of EM and SM that we introduce in the following section.

\smallskip

The first observation we can make regarding the Sideward Merge part of the Merge dynamical
system is that it plays a role similar to External Merge but in the ``opposite direction" namely
moving the dynamics toward more disconnected structures. Thus, if we keep the Internal Merge
part of the dynamics unchanged, and we consider only the SM arrows without the EM arrows,
we obtain an overall picture similar to what we have seen with IM and EM, but in the reverse
direction. We assume, as in the full Merge dynamics $\hat\cK$ introduced in \S \ref{GraphSec},
that only the minimal Sideward Merge arrows are included (those that only extract two
single leaves). 

\begin{prop}\label{IMSMgraph}
Let $\cG^{\rm IM-SM}_{n,A}\subset \cG_{n,A}$ be the subgraph on the same set of vertices, with
edges given only by the IM and the minimal SM arrows (without any EM arrows). 
The associated
Markov chain $\hat\cK^{\rm IM-SM, (A,n)}$ is reducible.
The strongly connected components (communication classes) are given by the following cases.
\begin{itemize}
\item For partitions $\wp\in \cP'(n)$ containing no pair $k_i=2$, $k_j=1$, 
the strongly connected components are graphs $\cG^{\rm IM}_{n,A,\wp,\sigma}$, with $\sigma\in C_\wp$; 
\item For partitions $\wp\in \cP'(n)$ containing some $\{ 2, 1 \}$ pairs, the strongly connected components are 
graphs $\cG^{\rm IM-SM}_{n,A,\wp,\tilde\sigma}$ with
$\tilde\sigma\in C_{\tilde\wp}$ for $\tilde\wp$ the partition (in some $\cP'(m)$ with $m< n$)
obtained from $\wp$ by removing all the $\{ 2, 1 \}$ pairs.
\end{itemize}
The only closed communication class is given by the induced subgraph on $\cG^{\rm IM-SM}_{n,A,\wp=\{ 2, 1, \ldots, 1 \}}$. 
This graph has no IM arrows and only minimal SM arrows.
The stationary distribution is the uniform distribution supported on the vertices of the subgraph 
$\cG^{\rm IM-SM}_{n,A,\wp=\{ 2, 1, \ldots, 1 \}}$. 
The condensation graph $\bar \cG^{\rm IM-SM}_{n,A}$ is the graph $\tilde\cG^{\rm SM}_{\cP'(n)}$
with vertex set 
$$ V(\tilde\cG^{\rm SM}_{\cP'(n)})=\{ (\wp,\sigma) \,|\, \wp\in \cP_{2,1}'(n)^c \, \sigma\in C_\wp \}
\cup \{ (\wp,\sigma) \,|\, \wp\in \cP_{2,1}'(n) \, \sigma\in C_{\tilde\wp} \} $$
where $\cP_{2,1}'(n)$ is the set of partitions with $\{ 2,1 \}$ pairs and $ \cP_{2,1}'(n)^c$ its complement,
and SM arrows $(\wp,\sigma) \to (\wp', \sigma')$ for $\wp=\{ k_1, \ldots, k_r \}$ and
$\wp'=\{ 2, k_1, \ldots, k_i-2, \ldots, k_r \}$ or $\wp'=\{ 2, k_1, \ldots, k_i-1, \ldots, k_j-1, \ldots k_r \}$
with $\sigma'$ the image of $\sigma$ under the induced map $C_\wp \to C_{\wp'}$.
This condensation graph has source vertices the pairs $(\wp,\sigma)$ with partitions $\wp\in \cP'(n)$ where all $k_i\neq 2$. 
\end{prop}

\proof Any workspace $F\in \cG^{\rm IM-SM}_{n,A,\wp}$ with $\wp\neq \{ 2, 1, \ldots, 1 \}$ is a transient state.
Indeed, any workspace $F$ that contains either more than one cherry component or at least one component
with $k_i\geq 3$ leaves has at least one minimal SM arrow with target $F'$ in a different $\cG^{\rm IM-SM}_{n,A,\wp'}$.
In the case with some $k_i\geq 3$, there are some SM arrows where 
the partition $\wp'$ replaces a $k_i$ with a pair $\{ 2, k_i-2 \}$, while 
in the case with more than one cherry, $\wp'$ replaces a pair $\{ k_i=2, k_j=2 \}$
with a triple $\{ 2, 1, 1 \}$, increasing the number of connected components, and 
there are no SM (nor IM) arrows that can go in the reverse direction, decreasing the number of components. 
On the other hand in the induced subgraph 
$\cG^{\rm IM-SM}_{n,A,\wp=\{ 2, 1, \ldots, 1 \}}$ there are no IM arrows and the only arrows in $\cG^{\rm IM-SM}_{n,A}$
that start at a vertex of $\cG^{\rm IM-SM}_{n,A,\wp=\{ 2, 1, \ldots, 1 \}}$  are SM arrows of the form
\begin{equation}\label{SMperms}
 \Tree[ $\alpha_i$ $\alpha_j$ ] \sqcup \alpha_k \bigsqcup_{\ell \neq i,j,k} \alpha_\ell \mapsto 
\Tree[ $\alpha_i$ $\alpha_k$ ] \sqcup \alpha_j \bigsqcup_{\ell \neq i,j,k} \alpha_\ell \, . 
\end{equation}
These arrow preserve $\cG^{\rm IM-SM}_{n,A,\wp=\{ 2, 1, \ldots, 1 \}}$ so this subgraph consists of recurrent states.

For $\wp$ with no $2$'s and $1$'s among the $k_i$'s (these have $k_i\geq 3$) the subgraph 
$\cG^{\rm IM-SM}_{n,A,\wp,\sigma}$ of $\cG^{\rm IM-SM}_{n,A,\wp}$ on the same vertices as
$\cG^{\rm IM}_{n,A,\wp,\sigma}$ have no minimal SM arrows mapping the 
vertex set to itself, so these graphs agree with the IM graphs $\cG^{\rm IM}_{n,A,\wp,\sigma}$ that 
are strongly connected and aperiodic. 

When the partition contains some $2$'s and some $1$'s,
on the other hand, $\cG^{\rm IM-SM}_{n,A,\wp}$ contains minimal SM edges that  map a subgraph 
$\cG^{\rm IM-SM}_{n,A,\wp,\sigma}$ to some other $\cG^{\rm IM-SM}_{n,A,\wp,\sigma'}$ through permuting some of the leaf 
labels as in \eqref{SMperms}. Thus, some different  $\cG^{\rm IM-SM}_{n,A,\wp,\sigma}$ become part of the same
strongly connected component in this case, and the resulting strong components depend only on the assignments
of labels on the rows of the diagram $D(\tilde\wp)$. In fact all the labels in cherries and the single
leaves components of these workspaces can be permuted by repeated application of arrows like \eqref{SMperms}.
The condensation graph $\bar \cG^{\rm IM-SM}_{n,A}$ has vertices the pairs $(\wp,\sigma)$ that
list the strongly connected components of $\cG^{\rm IM-SM}_{n,A}$. Among these, the ones that are sources
are those that have no incoming SM arrows, hence where the workspace has no cherry component. 
\endproof

As we observed above, for a partition $\wp$ that contains at least a pair $k_i=2$, $k_j=1$, there is,
for each such pair, an SM transformation mapping $\cG_{n,A,\wp}$ to itself, that acts on the 
components $T_i =\fM( \alpha_a, \alpha_b )$ and $T_j=\alpha_c$ as
$$  \Tree[ $\alpha_a$ $\alpha_b$ ] \sqcup \alpha_c \mapsto \Tree[ $\alpha_a$ $\alpha_c$ ] \sqcup \alpha_b $$
and as the identity on the other components of the workspace $F=T_1\sqcup \cdots \sqcup T_r$. 
These are the only kinds of SM arrows that preserve $\wp$. We denote them here as ${\rm SM}_{2,1}$.
For any $F\in V(\cG_{n,A,\wp})$ we have
\begin{equation}\label{inoutSM21}
N_{{\rm SM}_{2,1}}^{\rm out}(F) = N_{{\rm SM}_{2,1}}^{\rm in}(F) =\#\{ \text{ pairs } \{ 2, 1 \} \text{ in } \wp \} \, . 
\end{equation}

\section{External and Sideward Merge dynamics via Young diagrams}

Our analysis so far has shown that, in order to understand the properties of the full Merge
dynamics on $\cG_{n,A}$, which includes all the EM, IM, and minimal SM arrows, we need
to understand how the SM and EM arrows interact with each other and drive the behavior
of the dynamics. As before, we restrict to only the minimal SM arrows, since these are enough
SM arrows to make the Hopf algebra Markov chain ergodic and the results of \cite{MLH} also
show that these suffice to account for any linguistic phenomena for which SM may be a preferable
explanation. 

As we have previously discussed in \S \ref{inoutdegSec} (and also in \cite{MLH}),
one can also choose to only allow those minimal SM arrows that do not
cut both edges of a cherry. This choice is
natural because it gives somewhat better algebraic properties for the deletion
coproduct $\Delta^d$ (as will be discussed elsewhere). On the other hand, 
we will see that including these SM arrows makes the dynamical properties 
of the Merge Hopf algebra Markov chain nicer. Thus, given that this is the
main focus of our investigation, we will always 
make here  the assumption that we include among the minimal SM arrows
also those that extract the two leaves of a cherry. We can limit these
to the case where the cherry is a subtree of a larger tree, as when
the cherry is a component, the resulting SM would be just the identity
so we can ignore it. 

Since we already know the behavior of the IM dynamics, the first step in further 
investigating the properties of the Hopf algebra Markov chain $\hat\cK$ with IM, EM, and 
minimal SM arrows is to study graphs $\cG^{\rm EM-SM}_{\cP'(n)}$ with only EM and minimal SM arrows. 
Note that the graph $\cG^{\rm EM-SM}_{\cP'(n)}$ (unlike the $\tilde\cG^{\rm EM}_{\cP'(n)}$
and $\tilde\cG^{\rm SM}_{\cP'(n)}$ discussed above) is no longer the condensation 
graph of $\cG_{n,A}$ because $\cG_{n,A}$ is itself strongly connected so it has 
only one strongly connected component and no condensation graph. 

We analyze the EM--SM dynamics in steps. First we introduce the graph $\cG^{\rm EM-SM}_{\cP'(n)}$
that has vertex set $\cP'(n)$ and arrows that describe the effect of EM and minimal SM on
partitions of leaves in the workspace. This graph does not include the  tree topologies and
does not account for assignments of labels in $A=\{ \alpha_1, \ldots, \alpha_n \}$
at the leaves. 

As a second step, we introduce a way to keep track of the labels $\alpha_i$
at the leaves, by replacing the Young diagrams describing the partitions in $\cP'(n)$ with
Young tableaux on the alphabet $A$. Finally, we discuss how to account for the different 
tree topologies via appropriate multiplicities. 

We will discuss in \S \ref{FullDynSec} how the full Merge dynamics on
$\cG_{n,A}$ relates to the dynamics on $\cG^{\rm EM-SM}_{\cP'(n)}$ that 
we discuss in this section, through an intermediate graph $\cG_{\cP'(n)}$
that also incorporates the IM part of the dynamics.

\subsection{The Young diagrams graph} \label{YoungDsec}

Young diagrams are usually organized as a lattice either with respect to the dominance order on $\cP(n)$ 
 or by inclusion of diagrams (also called covering relation) on $\cP=\sqcup_n \cP(n)$ (Young lattice). These
 order relations and the resulting graphs can be interpreted as discrete dynamical systems, which have
 applications to models like the Sand Pile Model (a dynamical system exhibiting self-organized criticality). 
 
 \smallskip
 
 We consider here a different kind of discrete dynamical system on Young diagrams (not
 coming from a lattice structure), which is a Markov
 chain $\hat\cK_{\cP'(n)}$ induced on the set of partitions $\cP'(n)$
  by the Merge action $\hat\cK^{(A,n)}$ on $\cG_{n,A}$ through a projection map,
  resulting in a graph $\cG^{\rm EM-SM}_{\cP'(n)}$, which we consider here with only EM and SM
  arrows, and we extend in \S \ref{FullDynSec} by additional IM self-loops. 
  
 \begin{defn}\label{projmap}
 Let $p: \fF_{A,n} \to \cP'(n)$ denote the projection map that
 assigns to a forest $F=T_1 \sqcup \cdots \sqcup T_r$ in $\fF_{A,n}$
 the partition $\wp=p(F)$, $n=k_1+\cdots + k_r$, defined by the partitioning 
 of the $n$ leaves into the different components of the workspace, 
 $$ L(F) =L(T_1) \sqcup \cdots \sqcup L(T_r) $$
 into the $k_i=L(T_i)$ leaves of each component $T_i$. 
 \end{defn} 
  
We consider here a graph that has partitions as vertices and edges
corresponding to the effect on partitions of the EM and SM arrows on
workspaces.

   \begin{defn}\label{graphGPn}
   Let $\cG^{\rm EM-SM}_{\cP'(n)}$ denote the graph that has vertex set given by the set $\cP'(n)$ of partitions 
   of $n$, excluding the partition $\wp_1=\{ 1,1, \ldots, 1 \}$, representing the partitions $\wp$
   as their Young diagrams $D(\wp)$. The edges 
   are induced by the EM  and SM arrows of $\cG_{n,A}$: there is an EM or SM
   edge $\wp \to \wp'$ in $\cG^{\rm EM-SM}_{\cP'(n)}$ whenever there is an EM or SM arrow $F\to F'$
   for some $F\in p^{-1}(\wp)$ and some $F'\in p^{-1}(\wp')$ in the graph $\cG_{n,A}$,
   under the projection $p: \fF_{A,n} \to \cP'(n)$. The EM arrows in $\cG^{\rm EM-SM}_{\cP'(n)}$ combine
 together two of the rows of the diagram into a single row, while the SM arrows extract either
 two boxes from the same row or two boxes from different rows and form with them a new row 
 of length $2$. 
   \end{defn}

 For example, in the case $n=5$ the graph $\cG^{\rm EM-SM}_{\cP(5)}$ takes the form   {\scriptsize
  $$ \begin{tikzcd}
  \ydiagram{2,1,1,1} \arrow[bend left=10, r, "{\rm EM}"] \arrow[bend left=30, rr, "{\rm EM}"] \arrow[loop above, distance=5em, "{\rm SM}"] & \ydiagram{2,2,1} \arrow[bend left=10, l, "{\rm SM}"]  \arrow[bend left=30, rr, "{\rm EM}"]  \arrow[bend right=40, rrr, "{\rm EM}"]  \arrow[loop above, distance=5em, "{\rm SM}"] & \ydiagram{3,1,1} \arrow[bend right=30, rr, "{\rm EM}"] \arrow[r, "{\rm EM}"]  \arrow[l, "{\rm SM}"]  \arrow[bend left=50, ll, "{\rm SM}"]  &
  \ydiagram{3,2}  \arrow[bend left=30, rr, "{\rm EM}"]  \arrow[bend left=50, ll, "{\rm SM}"] & \ydiagram{4,1} \arrow[r, "{\rm EM}"] \arrow[l, "{\rm SM}"] \arrow[bend right=40, lll, "{\rm SM}"] & \ydiagram{5} \arrow[bend left=30, ll, "{\rm SM}"]
  \end{tikzcd} $$ }

 We now investigate some general properties of the graph $\cG^{\rm EM-SM}_{\cP'(n)}$ and the Markov
 chain $\hat\cK^{\rm EM-SM}_{\cP'(n)}$.

 \subsection{Strong connectivity}\label{strconnSMEMsec}
 Strong connectivity in fact follows from strong connectivity of the
 $\cG_{n,A}$ graph via projection, but we spell out a direct proof anyway (which reflects
 the proof for $\cG_{n,A}$ given in \cite{MLH}), because it will be useful later.
 
 \begin{lem}\label{strongconnGPn}
 The graph $\cG^{\rm EM-SM}_{\cP'(n)}$ for all $n\geq 4$, is strongly connected and aperiodic. 
 Thus, the Markov chain $\hat\cK^{\rm EM-SM}_{\cP'(n)}$ is ergodic: it has a unique stationary
 distribution (non-uniform, in general). Any initial distribution converges to the 
 stationary distribution under iterations of $\hat\cK^{\rm EM-SM}_{\cP'(n)}$.
 \end{lem}
 
 \proof To show strong connectivity it suffices to show that, for every $\wp,\wp' \in V(\cG^{\rm EM-SM}_{\cP'(n)})=\cP'(n)$, there is
 a directed path of minimal SM arrows from $\wp$ to the partition $\{ 2, 1, \ldots, 1 \}$ and that
 there is a directed path of EM arrows from $\{ 2, 1, \ldots, 1 \}$ to any other $\wp'$.  For the first one, 
 suppose that $\wp=\{ k_1, \ldots, k_r \}$. For any $k_i\geq 3$ there is a sequence of SM
 arrows that extract two boxes from a row of length $k_i$ and gives $\{ 2, k_1, \ldots, k_i-2, \ldots, k_r \}$,
 until each such $k_i$ is reduced to $2$ or $1$. Applying this sequence of SM transformations
 results in a partition that has only $2$'s and $1$'s. Any such partition can be reduced to 
 the partition with a single $2$ and all $1$'s by minimal SM arrows by repeatedly extracting
 one box from two rows of length $2$ resulting in one new row of length $2$ and two rows of
 length $1$ (and all the other rows unchanged), this reduces by one the number of rows of length
 $2$ and increases by two the number of rows of length $1$, until one reaches $\{ 2, 1, \ldots, 1 \}$. 
 To obtain the EM path from $\{ 2, 1, \ldots, 1 \}$ to another vertex $\wp'$, note that $\wp'$ 
 contains at least one row of length $k_i\geq 3$ or/and additional rows of length $2$.
 Each of these rows of length $k_i \geq 2$ can be obtained by successive EM arrows that start
 with combining two rows of $\{ 2, 1, \ldots, 1 \}$ and keep adding one box at a time to
 rows, using the remaining rows of length $1$ of $\{ 2, 1, \ldots, 1 \}$.
 The aperiodicity is immediate because of the self-loop
 SM arrows at $\{ 2, 1, \ldots, 1 \}$ and at other partitions consisting of $2$'s and $1$'s. 
 Ergodicity then follows from strong connectedness and aperiodicity. 
\endproof

Note that even if we restrict to the subgraph without the SM self-loops  we
still do have aperiodicity, for all $n\geq 4$, because any three $\wp, \wp', \wp''$ that contain the diagrams 
shown below have both a length $2$ and a length $3$ cycle coming from this subgraph of 
$\cG^{\rm EM-SM}_{\cP(4)}$: {\scriptsize
$$ \begin{tikzcd} \ydiagram{2,1,1} \arrow[bend left=30, rr, "{\rm EM}"] & \ydiagram{2,2} \arrow[l, "{\rm SM}"]  & \ydiagram{3,1} \arrow[l, "{\rm SM}"] 
\arrow[bend left=30, ll, "{\rm SM}"] 
 \end{tikzcd} $$ }
 
 \begin{rem}\label{strongkconn}{\rm 
 While strong connectedness holds for all the graphs $\cG^{\rm EM-SM}_{\cP'(n)}$, direct inspection of the
 example of $n=5$ above shows that it is not a ``robust" property of the graph, in the sense
 that it is sufficient to remove a single edge (the SM arrow between the $\{ 2,2 , 1\}$ and the $\{ 2,1,1,1 \}$
 partitions) to lose the strong connectivity property. Indeed, without that arrow, the graph would split
 out into two strongly connected components (the first recurrent, the second transient): {\scriptsize
 $$ \left\{ \ydiagram{2,2,1}\, , \,\ydiagram{3,2}\, , \, \ydiagram{4,1}  \, , \, \ydiagram{5} \right\} \ \ \text{ and } \ \ \left\{ \ydiagram{2,1,1,1}\, , \, \ydiagram{3,1,1} 
 \right\} \, .
 $$ } 
 This ``robustness" property can be measured in terms of {\em strong $k$-edge connectedness}. A graph
 is $k$-edge connected if it remains connected after the removal of up to $k-1$ edges and it is strongly
 $k$-edge connected if it remains strongly connected after removal of up to $k-1$ edges. The graph $\cG^{\rm EM-SM}_{\cP(5)}$
 is strongly $1$-edge connected (namely it is strongly connected) but it is not strongly $2$-edge connected.  }
 \end{rem}
 
 Indeed, the strong edge connectedness never improves when $n$ grows. 
 
 \begin{prop}\label{strongedgeconn}
 For all $n\geq 3$ the graphs $\cG^{\rm EM-SM}_{\cP'(n)}$ are strongly $1$-edge-connected but never strongly $2$-edge-connected.
 \end{prop}
 
 \proof Strong $1$-edge-connectedness (strong connectedness) is shown in Lemma~\ref{strongconnGPn}.
 The case $n=3$ just consists of two vertices with an arrow in each direction, so clearly removing one
 edge causes the strong connectedness property to fail. For any $n\geq 4$ consider the two paritions
 $\wp=\{ 2, 1, \ldots, 1 \}$ and $\wp'=\{ 3, 1, \ldots, 1 \}$. The partition $\wp'$ cannot be the target of any minimal SM
 arrow because it does not contain any $k_i=2$, so it can only be reached via an EM arrow. There is only
 one such arrow, which is coming from the vertex labelled by $\wp$. Thus, removing a single edge
 given by this EM arrow $\wp\longrightarrow \wp'$ will cause $\wp'$ to become a  transient
 state, in fact a source vertex, so the remaining graph is no longer strongly connected.
 \endproof
 
 This shows that ``fragility" of strong connectedness in the $\cG^{\rm EM-SM}_{\cP'(n)}$ networks is coming from the fact that
 some vertices cannot be directly reached via minimal SM arrows. Those vertices are only reachable
 via EM arrows. This makes certain EM arrows particularly crucial to the strong connectivity of the graph.
 Note that the argument of Proposition~\ref{strongedgeconn} identifies a crucial kind of EM arrow 
 as being, in the early stages of structure formation, the arrows that go from a cherry and a single
 leaf to a $3$-leaves structure
 $$ {\rm EM} : \,\,\, \alpha \sqcup \Tree[ $\beta$ $\gamma$ ] \longmapsto \Tree[ $\alpha$ [ $\beta$ $\gamma$ ] ] \, . $$
 
 One can argue, from the linguistics perspective, that these EM arrows do indeed play a special role, as
 they are the ones that make it possible to transition from simple binary associations $\{ \alpha, \beta \}$
 of two elementary structures, to actual hierarchical structures like $\{\{ \alpha, \beta \}, \gamma \}$. 
 
 This property also shows the fundamental asymmetry between EM and SM, with EM (and especially these
 particular EM arrows) playing a more crucial role in determining the properties of the dynamics. 
 
 Note, however, that if one includes all SM arrows, not just the minimal ones, then there would be SM arrows
 with target a partition that does not contain any $k_i=2$, such as, for example
 $$ \ydiagram{2,2,1} \stackrel{\text{non-min SM}}{\longrightarrow} \ydiagram{3,1,1} $$
 corresponding to non-minimal SM transformations like
 $$ \text{non-min SM} : \Tree[ $\alpha_1$ $\alpha_2$ ] \sqcup \Tree[ $\alpha_3$ $\alpha_4$ ] \sqcup \alpha_5 \longmapsto
  \Tree[ $\alpha_2$ [ $\alpha_3$ $\alpha_4$ ] ] \sqcup \alpha_1 \sqcup \alpha_5 \, . $$
  Thus, maintaining minimality of SM (as discussed in \cite{MLH}, see also \S \ref{CostSec} below) 
  causes the ``fragility" of strong connectedness and the special role of certain EM arrows. 
 
 \subsection{Young tableaux diagram}\label{YtablSec}
 
 As an intermediate step, before returning to discuss the full Merge action on workspaces
 and the graph $\cG^{\rm EM-SM}_{n,A}$, we refine 
 the graph $\cG^{\rm EM-SM}_{\cP'(n)}$ by including the assignment of labels $\alpha_i$ at the leaves,
 partitioned according to the partition $\wp\in \cP'(n)$. Thus, for $A=\{ \alpha_1, \ldots, \alpha_n \}$, 
 we consider the graph $\cG^{\rm EM-SM}_{\cP'(n), A}$ where the vertices are pairs $(\wp,\alpha)$ with
 $\alpha: A \to Y(\wp)$ a map that bijectively assigns the elements $\alpha_i \in A$ to the
 boxes of the Young diagram $Y(\wp)$ of the partition $\wp\in \cP'(n)$, considered up to the action of the
 Young group $S_\wp$ of permutations within the rows, and up to permuting rows of equal length, so
 that the total number of such assignments $\alpha: A \to Y(\wp)$ is given by the generalized
 multinomial coefficient $\Upsilon_{\wp,n}$ of \eqref{genmultinomial}.
 
   Again we do not include
 the partition $\{ 1, \ldots, 1 \}$ in $V(\cG^{\rm EM-SM}_{\cP'(n), A})$. The maps are then the same maps EM and SM 
 arrows as in $\cG^{\rm EM-SM}_{\cP'(n)}$ on the partitions, but now one also keeps track of the labeling, so that,
 for example an arrow like {\scriptsize
 $$ \ydiagram{3,1,1} \stackrel{\rm SM}{\longrightarrow} \ydiagram{2,2,1} $$ }
 splits into two arrows {\scriptsize
 $$ 
  \begin{ytableau}
  \alpha_1 & \alpha_2 \\
  \alpha_3  & \alpha_4 \\
  \alpha_5 
  \end{ytableau}
 \stackrel{\rm SM}{\longleftarrow}
  \begin{ytableau}
 \alpha_1 & \alpha_2 & \alpha_3 \\ \alpha_4 \\ \alpha_5 
 \end{ytableau}
 \stackrel{\rm SM}{\longrightarrow}
 \begin{ytableau}
  \alpha_1 & \alpha_2 \\
  \alpha_3  & \alpha_5 \\
  \alpha_4 
  \end{ytableau}
 $$ }
 and similarly for the other arrows of $\cG^{\rm EM-SM}_{\cP'(n)}$ that have different possible preimages
 under the projection map $\cG^{\rm EM-SM}_{\cP'(n), A} \to \cG^{\rm EM-SM}_{\cP'(n)}$ that forgets the labels $\alpha_i\in A$.
 
 \begin{rem}\label{noloopsSM}{\rm
 Note that the self-looping SM maps present in the graph $\cG^{\rm EM-SM}_{\cP'(n)}$, at vertices given by
 partitions $\wp$ containing $2$'s and $1$'s, are no longer self-loops in $\cG^{\rm EM-SM}_{\cP'(n), A}$,
 since now they map a vertex $(\wp,\alpha)$ to a different vertex $(\wp,\alpha')$ with $\alpha\neq \alpha'$, as 
 for example {\scriptsize
 $$ \begin{ytableau}
  \alpha_1 & \alpha_2 \\
  \alpha_3   
  \end{ytableau}
  \stackrel{\rm SM}{\longrightarrow}
  \begin{ytableau}
  \alpha_1 & \alpha_3 \\
  \alpha_2   
   \end{ytableau}
$$ }
 }\end{rem}

 \begin{lem}\label{GPnAstrconn}
 The graph $\cG^{\rm EM-SM}_{\cP'(n), A}$ is strongly connected and aperiodic, hence the
 Markov chain $\hat\cK^{\rm EM-SM}_{\cP'(n), A}$ is ergodic. 
 \end{lem}
 
 \proof
 Given the strong connectedness of $\cG^{\rm EM-SM}_{\cP'(n)}$, we can show that $\cG^{\rm EM-SM}_{\cP'(n), A}$ is also
 strongly connected by showing that any permutation of the labels on a fixed partition $\wp$
 can be achieved by a concatenation of arrows in $\cG^{\rm EM-SM}_{\cP'(n), A}$. Consider the following
 case as illustration of the general method. Suppose in the following Young tableau we want
 to permute the labels $\alpha_3$ and $\alpha_5$. This can be done with the following chain
 of SM and EM arrows:
 {\scriptsize
 $$
 \begin{ytableau}
 \alpha_1 & \alpha_2 & \alpha_3 & \alpha_4 \\
 \alpha_5 & \alpha_6
  \end{ytableau}
 \stackrel{\rm SM}{\longrightarrow} 
  \begin{ytableau}
 \alpha_1 & \alpha_2 & \alpha_4 \\
 \alpha_5 & \alpha_3 \\
 \alpha_6
  \end{ytableau}
   \stackrel{\rm SM}{\longrightarrow} 
   \begin{ytableau}
 \alpha_1 & \alpha_2 \\
  \alpha_5 &  \alpha_4 \\
 \alpha_3  \\
 \alpha_6
  \end{ytableau} 
  \stackrel{\rm EM}{\longrightarrow} 
     \begin{ytableau}
 \alpha_1 & \alpha_2 &  \alpha_5 &  \alpha_4 \\
 \alpha_3  \\
 \alpha_6
  \end{ytableau} 
  \stackrel{\rm EM}{\longrightarrow} 
    \begin{ytableau}
 \alpha_1 & \alpha_2 &  \alpha_5 &  \alpha_4 \\
 \alpha_3  &  \alpha_6
  \end{ytableau} 
  $$}
  This same procedure can be used to permute any two labels in different rows of the diagram.
  Since the assignments of labels to the diagrams is only up to permutations within rows, compositions
  of these operations suffice to achieve all possible permutations. 
 The same argument used for $\cG^{\rm EM-SM}_{\cP'(n)}$ shows aperiodicity for $\cG^{\rm EM-SM}_{\cP'(n), A}$. 
 \endproof
 
 Thus, we obtain the following result.
 
 \begin{cor}\label{stationaryGPA}
 The Markov chain $\hat\cK^{\rm EM-SM}_{\cP'(n), A}$ induced on $\cG^{\rm EM-SM}_{\cP'(n), A}$ by the Merge Hopf algebra
 markov chain $\hat\cK^{(A,n)}$ on $\cG_{n,A}$ has a unique asymptotic distribution
 $$ \pi_{\cP'(n), A}\,\, \hat\cK^{\rm EM-SM}_{\cP'(n), A} = \pi_{\cP'(n), A}\, , $$
 with $\pi_{\cP'(n), A} (\wp,\alpha) >0$ for all $(\wp,\alpha)\in V(\cG^{\rm EM-SM}_{\cP'(n), A})$ 
 and any initial distribution converges to $\pi_{\cP'(n), A}$ under iterations of $\hat\cK^{\rm EM-SM}_{\cP'(n), A}$.
  \end{cor}

Note that, unlike Internal Merge dynamics, the stationary distribution $\pi_{\cP'(n), A}$ on $\cG^{\rm EM-SM}_{\cP'(n), A}$
for External and Sideward Merge dynamics is usually not uniform.

\section{Full Merge dynamics with minimal Sideward Merge} \label{FullDynSec}

We now combine all the results that we have obtained by separately investigating the Internal Merge
dynamics on the graphs $\cG_{A,n,\wp,\sigma}$ and the combined Sideward Merge and External Merge 
induced dynamics on the graphs $\cG^{\rm EM-SM}_{\cP'(n)}$, to understand the full Merge dynamics
on $\cG_{n,A}$. To this purpose, we will consider here an extension of the graph $\cG^{\rm EM-SM}_{\cP'(n)}$,
to a graph $\cG_{\cP'(n)}$ that also includes the self-loops that are images of the IM arrows of $\cG_{n,A}$. 
To study the relation between the dynamics on this graph $\cG_{\cP'(n)}$ and on the full Merge graph
$\cG_{n,A}$ we first prove in \S \ref{ProjSec} a general result on projections of Markov chains. 

\subsection{Markov Chain Projection}\label{ProjSec}

We show here that when a finite Markov chain has a particular symmetry with respect
to a projection of the space of states to a smaller state, the stationary distribution can be
computed on a smaller graph. 

\begin{defn}\label{projreg}
Let $X$ be a finite set and $p: X \to Y$ a surjective map. Let $G_X$ be 
a strongly connected aperiodic directed graph with 
$V(G_X)=X$, and let $\cK^X$ be the adjacency matrix of $G_X$.  
We say that the graph $G_X$ is $p$-symmetric if it has the property that
the number of edges from a vertex $x$ to vertices in a fiber $p^{-1}(y')$ is constant
over $x\in p^{-1}(y)$ 
\begin{equation}\label{Rpsymm}
 \sum_{x'\in p^{-1}(y')} \cK^X_{x,x'} = \cK^{Y,R}_{y,y'} \, , \ \ \ \forall x \in p^{-1}(y) 
\end{equation}
and the number of edges to a vertex $x'$ from vertices in a
fiber $p^{-1}(y)$ is constant over $x'\in p^{-1}(y')$,
\begin{equation}\label{Lpsymm}
 \sum_{x\in p^{-1}(y)} \cK^X_{x,x'} = \cK^{Y,L}_{y,y'} \, , \ \ \ \forall x' \in p^{-1}(y') \, ,
\end{equation} 
with
\begin{equation}\label{LRKXY}
\# p^{-1}(y) \cdot \cK^{Y,R}_{y,y'} = \# p^{-1}(y') \cdot \cK^{Y,L}_{y,y'}\, .
\end{equation}
counting the total number of directed edges in $G_X$ from the fiber $p^{-1}(y)$
to the fiber $p^{-1}(y')$.
\end{defn}

The number computed by \eqref{LRKXY} can be interpreted geometrically
as follows. Let $G_Y$ be the graph with vertex set $Y$ and with a directed edge 
$y\to y'$ whenever there are some $F\in p^{-1}(y)$ and $F'\in p^{-1}(y')$ such
that there is a directed edge $F \to F'$ in $G_X$ between them. This way the
map $p: X \to Y$ extends to a map $p: G_X \to G_Y$. We refer to $G_Y$ as
the graph induced by $G_X$ and $p$.
Then both sides of \eqref{LRKXY} compute the number of edges of $G_X$ 
that are above an edge $y\to y'$ in $G_Y$ with respect to the projection $p$.

\begin{prop}\label{propPFpi}
Let $X$ be a finite set and $p: X \to Y$ a surjective map and $G_X$ a
strongly connected aperiodic $p$-symmetric graph on the vertex set $X$.
Then the induced graph $G_Y$ is strongly connected. Assume that
it is also aperiodic. 
Let $c_y=\# p^{-1}(y)$ denote the size of the fiber over a point $y\in Y$.
Let $\lambda$ and $\eta_X$ be the Perron-Frobenius eigenvalue
and (right) eigenvector of the adjacency matrix $\cK^X$ of $G_X$,
\begin{equation}\label{rightPFetaX}
\sum_{x'\in X} \cK^X_{x,x'}\, \eta_X(x') = \lambda\, \eta_X(x) \, ,
\end{equation}
and let $\hat\cK^X$ be the associated Markov chain transition matrix
$$ \hat\cK^X_{x,x'} := \lambda^{-1} \frac{\eta_X(x')}{\eta_X(x)} \cK^X_{x,x'} \, . $$
\begin{enumerate}
\item The right Perron-Frobenius  eigenfunction $\eta_X(x)=\eta_Y(y)$ of $\cK^X$ is constant on fibers $p^{-1}(y)$
for $y\in Y$ and its values $\eta_Y(y)$ are given by the right Perron-Frobenius eigenfunction of $\cK^{Y,R}$
\begin{equation}\label{PFYeq}
 \sum_{y'\in Y} \cK^{Y,R}_{y,y'} \,\, \eta_Y(y') = \lambda \, \eta_Y(y) \, . 
\end{equation} 
\item The left Perron-Frobenius  eigenfunction $\xi_X(x)=\xi_Y(y)$ of $\cK^X$ is also constant on fibers $p^{-1}(y)$
and given by the left Perron-Frobenius eigenfunction of $\cK^{Y,L}$
\begin{equation}\label{PFYeq}
 \sum_{y\in Y} \xi_Y(y)\,\, \cK^{Y,L}_{y,y'}  = \lambda \, \xi_Y(y) \, . 
\end{equation} 
\item The stationary distribution $\pi_X(x)$ of $\hat\cK^X$,
\begin{equation}\label{stathatKX}
\sum_{x\in X} \pi_X(x)\, \hat\cK^X_{x,x'} =\pi_X(x') 
\end{equation}
is also constant, $\pi_X(x)=\pi_Y(y)$, on fibers $x\in p^{-1}(y)$, with $\pi_Y(y)$ the stationary distribution
\begin{equation}\label{stathatKY}
\sum_{y\in y} \pi_Y(y)\, \,\hat\cK^Y_{y,y'} =\pi_Y(y') \, ,
\end{equation}
of the Markov chain on $Y$ with transition matrix
$$ \hat\cK^Y_{y,y'}= \lambda^{-1}\, \frac{\eta_Y(y')\cdot c_y}{\eta_Y(y)\cdot c_{y'}} \,\cK^{Y,R}_{y,y'} \, . $$
\item Equivalently, the function
\begin{equation}\label{psiY}
 \psi_Y(y):= \frac{\pi(y) \, c_y}{\eta_Y(y)} 
 \end{equation}
 is the (left) Perron-Frobenius eigenvector of $\cK^{Y,R}$, satisfying
 \begin{equation}\label{KYRleftPF}
 \sum_y \psi_Y(y) \, \cK^{Y,R}_{y,y'} = \lambda \, \psi_Y(y') \, , 
 \end{equation}
 \item The stationary distribution of the Markov chain $\hat\cK^X$ is given by
 \begin{equation}\label{statRLPF}
 \pi_X(x) = \frac{\eta_Y(p(x)) \cdot \psi_Y(p(x))}{c_{p(x)}} \, ,
 \end{equation}
 in terms of the left and right PF eigenvectors $\psi_Y$ and $\eta_Y$ of $\cK^{Y,R}$ and the sizes $c_y$ of
 the fibers of the projection $p: X \to Y$.
 \item the left Perron-Frobenius eigenvectors $\xi_Y$ of $\cK^{Y,L}$ and $\psi_Y$ of $\cK^{Y,R}$ are related
 (up to normalization) by
  \begin{equation}\label{LPFofKLR}
  \xi_Y(y) = \frac{\psi_Y(y)}{c_y}\, . 
  \end{equation}
 \end{enumerate}
\end{prop}

\proof The graph $G_Y$ is strongly connected since $G_X$ is. Indeed, given $y,y'\in Y$, 
two arbitrary vertices $F\in p^{-1}(y)$ and $F'\in p^{-1}(y')$ are connected by paths in $G_X$ in both directions
and subsequent edges of these paths in $G_X$ project to subsequent edges in $G_Y$.

(1) The graph $G_Y$ is also aperiodic by assumption, so the matrix $\cK^{Y,R}$ has 
a unique (up to scaling) right Perron-Frobenius eigenfunction with $\eta_Y(y)>0$ for all $y\in Y$, 
satisfying \eqref{PFYeq} with $\lambda>0$ the PF eigenvalue of $Y$. We then have
$$ \sum_{x'\in X} \cK^X_{x,x'}\, \eta_Y(y') = \sum_{y'}  \sum_{x'\in p^{-1}(y')} \cK^X_{x,x'} \, \eta_Y(y') =
\sum_{y'} \cK^{Y,R}_{y,y'} \, \eta_Y(y') = \lambda \eta_Y(y) \, , $$
so that $\eta_X(x):=\eta_Y(p(x))$ is a solution to \eqref{rightPFetaX} with the same $\lambda$.
The Perron-Frobenius theorem for $\cK^X$ ensures that the left/right PF eigenvectors are the
only eigenvectors of $\cK^X$ with all positive real entries. Thus, $\eta_X$ has to be the right
Perron-Frobenius eigenvector of $\cK^X$, which is unique up to scale. 

(2) The case of the left Perron-Frobenius eigenfunction of $\cK^X$ is analogous with
$$ \sum_{x\in X} \xi_X(x) \cK^X_{x,x'}  = \sum_y  \sum_{x\in p^{-1}(y)} \xi_Y(y) \cK^X_{x,x'} =
\sum_y \xi_Y(y) \cK^{Y,L}_{y,y'}  = \lambda \xi_Y(y)\, . $$

(3) In a similar way,  given
that $\eta_X$ is constant on preimages of $p$, if $\pi_Y(y)$ is the stationary distribution of
$\hat\cK^Y$ satisfying \eqref{stathatKY}, then $\pi_X(x):=\pi_Y(p(x))$ satisfies, for $x'\in p^{-1}(y')$,
$$ \sum_x \pi_X(x)\, \hat\cK_{x,x'}= \sum_y \pi_Y(y) \, \lambda^{-1} \frac{\eta_Y(y')}{\eta_Y(y)} 
\sum_{x\in p^{-1}(y)} \cK^X_{x,x'} = \sum_y \pi_Y(y) \, \lambda^{-1} \frac{\eta_Y(y')}{\eta_Y(y)} 
\cK^{Y,L}_{y,y'} $$
$$ = \sum_y \pi_Y(y) \, \lambda^{-1} \frac{\eta_Y(y')\cdot c_y}{\eta_Y(y) \cdot c_{y'}}  \cK^{Y,R}_{y,y'}  = \sum_y \pi_Y(y) \, \hat\cK^Y_{y,y'} = \pi_Y(y') = \pi_X(x') \, . $$
Since the graph $G_X$ is strongly connected and aperiodic, the stationary distribution $\pi_X$ 
is unique, so it has to be of this form. 

(4) We can equivalently write the equation
$$ \sum_{y\in y} \pi_Y(y)\, \,\lambda^{-1}\, \frac{\eta_Y(y')\cdot c_y}{\eta_Y(y)\cdot c_{y'}} \,\cK^{Y,R}_{y,y'} =\pi_Y(y') $$
in the form
$$ \sum_{y\in y} \frac{\pi_Y(y) c_y}{\eta_Y(y)} \,\cK^{Y,R}_{y,y'} = \lambda \, \frac{\pi_Y(y')  c_{y'}}{\eta_Y(y')} \, , $$
which says that  $\psi(y)$ as in \eqref{psiY} is the left Perron-Frobenius eigenvector equation for $\cK^{Y,R}$ 
with PF eigenvalue $\lambda$. 

(5) This follows directly from (4) and the form of the stationary distribution $\pi(x)$ as in \eqref{statRLPF}. 

(6) Since $\psi_Y$ is the left PF eigenvector of $\cK^{Y,R}$, we have
$$ \sum_y \frac{\psi_Y(y)}{c_y}\,\, \frac{c_y}{c_{y'}} \cK^{Y,R}_{y,y'} = \lambda \frac{\psi_Y(y')}{c_{y'}}\, , $$
which by \eqref{LRKXY} and the uniqueness (up to normalization) 
of the left PF eigenvector of $\cK^{Y,L}$ gives the relation \eqref{LPFofKLR}.

\endproof

\subsection{Stationary distribution of the Merge dynamics on workspaces} \label{statdistrMergeSec}
This result applies to the case of the Hopf algebra Markov chain of Merge we are considering here.

\begin{prop}\label{psymmMergeDyn}
Let $X=V(\cG_{n,A})=\fF_{A,n} \subset \fF_{\cS\cO_0}$ be the set of nonplanar full binary 
forests with $n$ leaves with set of labels $A=\{ \alpha_1. \ldots, \alpha_n \}$, and with non-empty
set of edges. Let $Y=\cP'(n)$
be the set of partitions of $n$ (not including the partition $\{ 1,\ldots, 1\}$).
Let $p: \fF_{A,n} \to \cP'(n)$ be the map that assigns to a forest $F$
with $n$ leaves the partition $\wp=p(F)$ of $n=k_1+\cdots+k_r$ into the numbers of leaves
of each tree component $F=T_1 \sqcup \cdots \sqcup T_r$ of the forest, $k_i=\# L(T_i)$.
The graph $\cG_{n,A}$ is $p$-symmetric with respect to this map. The projection 
$p: \fF_{A,n} \to \cP'(n)$ induces a projection of graphs $\cG_{n,A} \to \cG_{\cP'(n)}$,
where the graph $\cG_{\cP'(n)}$ has vertex set $\cP'(n)$ and EM, SM arrows as in Definition~\ref{graphGPn}
and an additional IM self-loop at a vertex $\wp$ whenever there are IM arrows $F\to F'$ between
some vertices $F,F'\in p^{-1}(\wp)=\cG_{n,A,\wp}$.
\end{prop}

\proof The key property of the graph $\cG_{n,A}$ that is responsible for the $p$-symmetry
with respect to the map $p: \fF_{A,n} \to \cP'(n)$ is the fact that the number of edges from
(or to) a vertex $F$ in $\cG_{n,A}$ depends only on $p(F)$ (assuming that we also include
Sideward Merge arrows that extract the two leaves of a cherry when that cherry is part
of a larger tree). This follows from Proposition~\ref{EdgesGnA}. 
This property, in fact, can be further refined to show that the counting of the edges out (or in)
a vertex $F$ on $\cG_{n,A}$ can be split into separate counting for the subsets of edges
with target (respectively, source) vertex in each fiber $p^{-1}(\wp')$, and the counting over
each of these sets also depends only on $\wp=p(F)$ and not on $F$ itself. To see this,
notice that for each type (EM/IM/SM) of arrows, the counting of the outgoing (or incoming)
edges of that type depends on counting choices of appropriate boxes in the partition $\wp=p(F)$.
Any such choice of boxes, along with the type of operation (EM/IM/SM), suffice to determine 
the shape of the resulting partition $\wp'$ associated to the target vertex of the arrow. Thus
each subset of arrows with source $F$ that have the other end in a fiber $p^{-1}(\wp')$ is
determined by knowing the partition $\wp$. This implies that
both \eqref{Lpsymm} and \eqref{Rpsymm} hold. By construction, the map
$\cG_{n,A} \to \cG_{\cP'(n)}$ has fibers $p^{-1}(\wp)=\cG_{n,A,\wp}$, hence
the EM and SM arrows of $\cG_{n,A}$ give rise to edges of $\cG_{\cP'(n)}$ as in
Definition~\ref{graphGPn} and the IM arrows, which are always along
components of fibers, give rise to looping edges in $\cG_{\cP'(n)}$. 
\endproof

Thus we obtain the form of the stationary distribution for the full Merge dynamics on the 
set of workspaces, given by the Hopf algebra Markov chain $\hat\cK^{A,n}$ on the
graph $\cG_{n,A}$.

\begin{thm}\label{MergeDynStatDistr}
The Merge Hopf algebra Markov chain $\hat\cK^{A,n}$ on the graph
$\cG_{n,A}$ has stationary distribution that is constant on the subsets
$V_{n,A,\wp}= p^{-1}(\wp)$, that are the fibers of the projection 
$p: \cG_{n,A} \to \cG_{\cP'(n)}$. Up to an overall multiplicative constant,
the stationary distribution $\pi(F)$ of the Merge Hopf algebra Markov chain $\hat\cK^{A,n}$ 
is obtained from the left and right Perron-Frobenius eigenvectors, $\eta_{\cP'(n)}, \psi_{\cP'(n)}$ of the
adjacency matrix $\cK^{\cP'(n),R}$ of the graph $\cG_{\cP'(n)}$ (with the IM looping edges included)
and the multiplicities of the fibers $\# V_{n,A,\wp}=\Lambda_{\wp,n}$ as in \eqref{Lgenmultinomial}, in the form
\begin{equation}\label{piMergeDyn}
\pi(F) = \frac{\eta_{\cP'(n)}(\wp) \, \, \psi_{\cP'(n)}(\wp)}{\Lambda_{\wp,n}} \ \ \ \text{ for } \wp=p(F) \, . 
\end{equation}
\end{thm}

\proof The form \eqref{piMergeDyn} of the distribution follows from Proposition~\ref{propPFpi}
and the fact that the distribution is uniform on the $V_{n,A,\wp}$ follows from \eqref{IMinoutd} and 
\eqref{inoutSM21}, as in Proposition~\ref{dynIMunif}.
\endproof

\bigskip
\bigskip

\subsection{The example of $n=4$}\label{ExSec}

We show here an explicit example, which is simple enough that everything can be computed
explicitly, but sufficient to illustrate all the steps of the procedure discussed above. The case
with $n=3$ leaves is too simple (and already fully discussed in \cite{MLH}) so we look at the
case of $n=4$ leaves. We first describe the geometry of the projection $p: \cG_{4,A} \to \cG_{\cP(4)}$.

\begin{lem}\label{ex4YD}
In the case $n=4$, the projection $p: \cG_{4,A} \to \cG_{\cP(4)}$ has multiplicities above each edge
given by 
{\scriptsize
$$ \begin{tikzcd} 
\ydiagram{2,1,1} \arrow[loop above, distance=5em, "24"]  \arrow[bend left=15, r, "6"]  
\arrow[bend left=30, rr, "12"] & 
\ydiagram{2,2}  \arrow[bend left=15, l, "12"]   \arrow[bend right=40, rr, "3"]  & 
\ydiagram{3,1}  \arrow[l, "36"] \arrow[bend left=40, ll, "36"]   \arrow[r, "12"]  
\arrow[loop above, distance=5em, "24"]  & 
\ydiagram{4} \arrow[bend left=50, ll, "45"] \arrow[loop above, distance=5em, "60"] 
 \end{tikzcd} $$ }
 while the multiplicities above the vertices are $6$ for the first vertex, $3$ for the second,
 $12$ for the third and $15$ for the fourth. 
\end{lem}

\proof The graph $\cG_{\cP(4)}$ with the EM and SM arrows, and with the self-looping IM arrows, 
is of the form
{\scriptsize
$$ \begin{tikzcd} 
\ydiagram{2,1,1} \arrow[loop above, distance=5em, "{\rm SM}"]  \arrow[bend left=15, r, "{\rm EM}"]  
\arrow[bend left=30, rr, "{\rm EM}"] & 
\ydiagram{2,2}  \arrow[bend left=15, l, "SM"]   \arrow[bend right=40, rr, "{\rm EM}"]  & 
\ydiagram{3,1}  \arrow[l, "{\rm SM}"] \arrow[bend left=40, ll, "{\rm SM}"]   \arrow[r, "{\rm EM}"]  
\arrow[loop above, distance=5em, "{\rm IM}"]  & 
\ydiagram{4} \arrow[bend left=50, ll, "{\rm SM}"] \arrow[loop above, distance=5em, "{\rm IM}"] 
 \end{tikzcd} $$ }
 To compute the weights of the arrows we need to compute the number of edges of
 $\cG_{n,A}$ that are above each of the edges in this graph through the projections
 $\cG_{n,A}\to \cG_{\cP'(n),A}\to \cG_{\cP'(n)}$. Since the number of incoming and
 outgoing edges at each vertex $F$ in $\cG_{n,A}$ only depends on the partition $\wp=p(F)$
 (as we are here also allowing SM arrows that cut two leaves of a cherry that is inside a larger tree),
 the number of edges above an edge $\wp\to \wp'$ can be computed as the product of the 
 number of vertices above $\wp$ in $\cG_{n,A}$ times the number of edges from one of
 the vertices above $\wp$ to the fiber above $\wp'$. Thus, we obtain the following:
 \begin{itemize}
\item There are $6$ vertices of $\cG_{4,A}$ above the vertex $\wp=\{ 2,1,1 \}$ of $\cG_{\cP(4)}$, namely
the choices of labels at the leaves
$$ \Tree[ $\alpha_i$ $\alpha_j$ ] \sqcup \alpha_k \sqcup \alpha_\ell $$
namely all choices of two out of $\{ \alpha_1, \alpha_2, \alpha_3, \alpha_4 \}$ to form the cherry leaves.
\item There are $3$ vertices of $\cG_{4,A}$ above the vertex $\wp=\{ 2,2 \}$ of $\cG_{\cP(4)}$, again given
by all choices of labels in
$$ \Tree[ $\alpha_i$ $\alpha_j$ ] \sqcup \Tree[ $\alpha_k$ $\alpha_\ell$ ] $$
which can be seen as the $6$ choices of two labels to form the first cherry, taken up to the symmetry between
the two cherries that cuts the counting down to $3$.
\item There are $12$ vertices of $\cG_{4,A}$ above the vertex $\wp=\{ 3,1 \}$ of $\cG_{\cP(4)}$ given by
the assignments of labels to
$$ \Tree[ $\alpha_i$ [ $\alpha_j$ $\alpha_k$ ] ] \sqcup \alpha_\ell $$
so $4$ choices of the label for the isolated vertex, and $3$ remaining choices of the label at the leaf attached
to the root in the tree.
\item There are $15$ vertices of $\cG_{4,A}$ above the vertex $\wp=\{ 4 \}$ of $\cG_{\cP(4)}$ given by
the sum over the two tree topologies of the different possible labelings at the leaves,
$$ \Tree[ [ $\alpha_i$  $\alpha_j$ ] [ $\alpha_k$ $\alpha_\ell$ ] ]  \ \ \text{ and } \ \  \Tree[  $\alpha_i$  [ $\alpha_j$ [ $\alpha_k$ $\alpha_\ell$ ] ] ] \, . $$
In the first case there are $3$ possible assignments of leaves (as in the case of $\wp=\{ 2,2 \}$) while in the second
case there are $12$ ($4$ choices of label at the vertex attached to the root and $3$ remaining at the second vertex).
\end{itemize}
We then just need to count the number of edges from one vertex above $\wp$ with target above $\wp'$
as this number is constant at all vertices above $\wp$. These are counted by our general counting of indegrees
and outdegrees for the $\cG_{n,A}$.
\begin{itemize}
\item There are $4$ Sideward Merge arrows from $F$ above $\wp=\{ 2,1,1 \}$ to the same fiber, hence $24 =6\cdot 4$
edges above $\{ 2,1,1 \} \to \{ 2,1,1 \}$.
\item There is $1$ External Merge arrow from $F$ above $\wp=\{ 2,1,1 \}$ to vertices above $\wp'=\{ 2,2 \}$, hence
$6$ edges above $\{ 2,1,1 \} \to \{ 2,2 \}$.
\item There are $2$ External Merge arrows from $F$ above $\wp=\{ 2,1,1 \}$ to vertices above $\wp'=\{ 3,1 \}$,
hence $12$ edges above $\{ 2,1,1 \} \to \{ 3,1 \}$.
\item There are $4$ Sideward Merge arrows from $F$ above $\wp=\{ 2,2 \}$ to vertices above $\wp'=\{ 2,1,1 \}$, hence
$12$ edges above $\{ 2,2 \}\to \{ 2,1,1 \}$.
\item There is $1$ External Merge arrow from $F$ above $\wp=\{ 2,2 \}$ to vertices above $\wp'=\{ 4 \}$, hence $3$ edges
above $\{ 2,2 \}\to \{ 4 \}$. 
\item There are $2$ Internal Merge arrows from a vertex $F$ above $\wp=\{ 3,1 \}$ to the same fiber, hence $24=12\cdot 2$
edges above $\{ 3,1 \} \to \{ 3,1 \}$.
\item There are $3$ Sideward Merge arrows from $F$ above $\wp=\{ 3,1 \}$ to vertices above $\wp'=\{ 2, 2\}$, hence
$36=12\cdot 3$ edges above $\{ 3,1 \} \to \{ 2,2 \}$.
\item There are $3$ Sideward Merge arrows from $F$ above $\wp=\{ 3,1 \}$ to vertices above $\wp'=\{ 2,1,1 \}$, hence
$36=12\cdot 3$ edges above $\{ 3,1 \} \to \{ 2,1,1 \}$.
\item There is $1$ External Merge arrows from $F$ above $\wp=\{ 3,1 \}$ to vertices above $\wp'=\{ 4 \}$, hence $12$ edges 
above $\{ 3,1 \} \to \{ 4 \}$. 
\item There are $3$ Sideward Merge arrows from $F$ above $\wp=\{ 4 \}$ to vertices  above $\wp'=\{ 2,2 \}$,
hence $45=15\cdot 3$ edges above $\{ 4 \} \to \{ 2,2 \}$.
\item There are $4$ Internal Merge arrows from a vertex $F$ above $\wp=\{ 4 \}$ to the same fiber, hence $60=15\cdot 4$ edges
above $\{ 4 \} \to \{ 4 \}$.
\end{itemize}
\endproof

We then compute the matrix $\cK^{\cP(4),R}$ of Proposition~\ref{propPFpi}.

\begin{lem}\label{KP4R}
For $n=4$, the matrix $\cK^{\cP(4),R}$ satisfying
$$ \cK^{\cP(4),R}_{\wp,\wp'}=\sum_{F'\in p^{-1}(\wp')} \cK^{(A,4)}_{F,F'} \, , $$
with $\cK^{(A,4)}_{F,F'}$ the adjacency matrix of $\cG_{4,A}$, is given by
\begin{equation}\label{matKP4R}
\cK^{\cP(4),R} = \begin{pmatrix}   
4   & 1 & 2 & 0   \\ 
4   & 0  & 0 &  1  \\ 
3   & 3 & 2 & 1   \\ 
 0  & 3 & 0 &  4
\end{pmatrix}
\end{equation}
and corresponds to the weighted graph (with vertices listed from left to right)
{\scriptsize
$$ \begin{tikzcd} 
\ydiagram{2,1,1} \arrow[loop above, distance=5em, "4"]  \arrow[bend left=15, r, "1"]  
\arrow[bend left=30, rr, "2"] & 
\ydiagram{2,2}  \arrow[bend left=15, l, "4"]   \arrow[bend right=40, rr, "1"]  & 
\ydiagram{3,1}  \arrow[l, "3"] \arrow[bend left=40, ll, "3"]   \arrow[r, "1"]  
\arrow[loop above, distance=5em, "2"]  & 
\ydiagram{4} \arrow[bend left=50, ll, "3"] \arrow[loop above, distance=5em, "4"] 
 \end{tikzcd} $$ }
The stationary distribution for $\hat\cK^{A,4}$ on $\cG_{4,A}$ takes
the form 
\begin{equation}\label{statdistr4}
 \pi_{4,A}(F)\sim \left\{ \begin{array}{llll}
   0.5267     & F\in p^{-1}(\{ 2,1,1,1 \}) \\
   0.3109     & F\in p^{-1}(\{ 2,2 \}) \\
   0.1217     & F\in p^{-1}(\{ 3,1 \}) \\
   0.0406    & F\in p^{-1}(\{ 4 \})\, .
\end{array}\right. 
\end{equation}
\end{lem}

\proof The matrix $\cK^{\cP(4),R}$ of \eqref{matKP4R} has 
Perron--Frobenius eigenvalue $\lambda\sim 6.9656$, right
Perron--Frobenius eigenvector with $\eta_{\cP(4)}(\{ 2,1,1,1 \})\sim 0.5565$,
$\eta_{\cP(4)}(\{ 2,2 \})\sim 0.3739$, $\eta_{\cP(4)}(\{ 3,1 \})\sim 0.6383$,
$\eta_{\cP(4)}(\{ 4 \})\sim 0.3782$, and left Perron--Frobenius eigenvector
$\psi_{\cP(4)}(\{ 2,1,1,1 \})\sim 0.8345$,
$\psi_{\cP(4)}(\{ 2,2 \})\sim 0.3666$, $\psi_{\cP(4)}(\{ 3,1 \})\sim 0.3361$,
$\psi_{\cP(4)}(\{ 4 \})\sim 0.2369$. 
The size $c_\wp=\# p^{-1}(\wp)$ of the fibers is computed in Lemma~\ref{ex4YD}
above, so using the expression \eqref{statRLPF} for the stationary distribution,
and normalizing the result so that it is a probability distribution,
we obtain \eqref{statdistr4}.
\endproof

\smallskip

\subsection{The Sideward Merge problem}\label{SMproblemSec}

One sees easily, already from the explicit example of $n=4$, the general problem produced by the
presence of the Sideward Merge arrows, even while restricting only to the minimal
ones: they tend to push the dynamics towards the more disconnected structures, entirely
preventing the process of structure formation. We see this reflected in the fact
that the fully connected structures (the trees in the preimage of
$\wp=\{ 4 \}$) have a very low weight in the stationary distribution, compared to
the most disconnected structures (in the preimage of $\wp=\{ 2,1,1 \}$) that
receive the highest weight. 

\smallskip

Notice how the behavior, already in the case $n=4$, is different from the case
$n=3$ analyzed in \cite{MLH}. In the case $n=3$, when EM, SM, IM arrows
are considered, the dynamics has a uniform stationary distribution where
all workspaces carry the same weight. This already shows a problem, created
by Sideward Merge: the dynamics does not converge to the fully formed
sentences (the workspaces consisting of a single tree). This problem is 
corrected by the introduction of cost functions, see \S \ref{CostSec}.

\smallskip

Indeed, this behavior clearly does not reflect the reality of linguistic
structure formation and it means that the different types of arrows, External
Merge, Internal Merge, and Sideward Merge, must have different weights in
the graph $\cG_{n,A}$, according to some cost function. As already discussed in
\cite{MLH}, Sideward Merge is indeed non-optimal with respect to some measures
of optimality used in the linguistics literature. Also there are only a few linguistic
phenomena that may require Sideward Merge as an explanation, so it should
occur very infrequently in Merge derivations, contrary to what one sees
in the dynamics of the Markov chain $\hat\cK^{A,n}$ we have been
discussing here, where Sideward Merge arrows end up playing a dominant
role in the dynamics. This means that a realistic model, as already discussed
in \cite{MLH}, will have to replace $\hat\cK^{A,n}$ with a weighted version,
where different costs are accounted for. We will be discussing this in
Section~\ref{CostSec}.

 \section{Merge dynamics with the contraction coproduct} \label{ContractSec}
 
 Before introducing cost functions and analyzing their effect on
 the different parts (EM, IM, SM) of the dynamics, we 
 comment in this section on possible variants of the dynamics
 discussed here. It is shown in \cite{MCB} that there
 are different possible forms of the coproduct on workspaces
 that is used to define the Merge operations. All the discussion
 and results in the previous sections are based on the use of
 what is called in \cite{MCB} the ``deletion" coproduct, where
 the quotient $T/^d F_{\underline{v}}$, for $F_{\underline{v}}=T_{v_1}\sqcup\cdots \sqcup T_{v_k}$
 a collection of disjoint accessible terms in $T$, is
 defined as the maximal full binary rooted tree obtained
 (via edge contractions) from the non-full binary tree
 (including non-branching vertices) resulting from the
 admissible cut that removes the edges above each of the vertices $v_i$
 (the roots of the accessible terms $T_{v_i}$). The essential property of
 this coproduct $\Delta^d$ that we used in analyzing the dynamics is that the
 number of leaves adds over this decomposition:
 $$ \# L(T)= \# L(F_{\underline{v}}) + \# L(T/^d F_{\underline{v}}) \, . $$
 In terms of the linguistics model, this coproduct is suitable for the
 interface of the core Merge dynamics with Externalization,
 where the trace of movement is not expressed. 
 
 The other coproduct used in \cite{MCB} is the ``contraction"
 coproduct, where the quotient $T/^c F_{\underline{v}}$ is obtained
 by shrinking each component $T_{v_i}$ of 
 $F_{\underline{v}}=T_{v_1}\sqcup\cdots \sqcup T_{v_k} \subset T$
 to its vertex $v_i$. Each of these vertices $v_i$ becomes in this way 
 a new leaf of $T/^c F_{\underline{v}}$, labelled by an element
 of $\cS\cO_0$ that maintains a ``trace" of the movement of the
 extracted accessible terms. One denotes such label by 
 the symbol $\text{\sout{$T_{v_i}$}}$, to mark that it corresponds
 to the cancellation of the lower copy of $T_{v_i}$ under movement. 
 This form of coproduct $\Delta^c$ now no longer satisfies the
 additivity of the number of leaves, as now we have
 $$ \# L(T) + k = \# L(F_{\underline{v}}) + \# L(T/^c F_{\underline{v}}) \, , $$
 for $F_{\underline{v}}=T_{v_1}\sqcup\cdots \sqcup T_{v_k}$.

 \subsection{Labels of traces of movement}\label{tracelabelSec}

 Before discussing further the Merge operation constructed
 using this coproduct (as discussed in \S 1.3 of \cite{MCB}), 
 it is useful to clarify the meaning of the labeling system
 $\text{\sout{$T_v$}}$ used to mark the new leaf $v$ of the quotient
 $T/^c T_v$, when the subtree $T_v \subset T$ is contracted
 to its root vertex $v$. This section clarifies and expands in more
 detail the brief discussion in \S 1.3 of \cite{MCB}. 
 
 The form of the quotient $T/^c T_v$ requires an explanation: since the
 label of the leaf resulting from the contraction of $T_v$ is usually written as
 $\text{\sout{$T_v$}}$, this does not look like it can be 
 compatible with having a fixed finite set of leaf labels.
 Indeed, if the label $\text{\sout{$T_v$}}$ would really contain
 the full information about $T_v$, together with a ``cancellation symbol"  to mark it as a ``trace"
 of movement, it would give rise to a countably infinite set of labels, since
 $T_v$ can be an arbitrary syntactic object (since $T$ itself can be arbitrary), and the set of all
 possible syntactic objects is countably infinite. However, the label $\text{\sout{$T_v$}}$ in fact 
 does not have to retain the full structure of $T_v$ as a syntactic object (as a full
 nonplanar binary rooted tree). After all, the extracted term $T_v$ is still present, on
 the other side of the coproduct, so that information is not lost. All that is needed is
 that  $\text{\sout{$T_v$}}$ remembers a feature of $T_v$ sufficient for interpretation: namely, that it
 remembers the role that the moved element $T_v$ played in its original position. 
 We show here that this can be achieved through a labeling algorithm that uses a head
 function (in the sense of \cite{MCB}, that we recall below) 
 and the structure of $T_v$ to assign a label to the non-leaf
 vertices, including the root $v$. 
 
 To begin with, we make a refined choice of the initial set $\cS\cO_0$ of
 lexical items and syntactic features. This is not needed in the formulation
 we gave in the previous sections,  as it is only relevant when the coproduct
 $\Delta^c$ is involved. 
 
 \begin{defn}\label{cancSO0def}
 We replace $\cS\cO_0$ with a set of the form
  $\cS\cO_0\sqcup\text{\sout{ $\cS\cO_0$ }}$, where $\text{\sout{ $\cS\cO_0$ }}$
 stands for the set of symbols of the form
 \begin{equation}\label{SO0canc}
 \text{\sout{ $\cS\cO_0$ }} := \{ \text{\sout{ $\alpha$ }} \,|\, \alpha \in \cS\cO_0 \} \, .
 \end{equation}
 We refer to elements of $\text{\sout{ $\cS\cO_0$ }}$ as ``trace labels".
 \end{defn}
 
 We also recall from \S 1.13 of \cite{MCB} the definition of a {\em head function}.
 
 \begin{defn}\label{headdef}
A head function on a full non-planar
binary rooted tree $T$ is a function $h_T: V^o(T) \to L(T)$, from the set of non-leaf
vertices of $T$ to the set of leaves, with the property that $h_T(v)\in L(T_v)$ and
that if $T_w \subset T_v$ and $h_T(v)\in L(T_w)$ then $h_T(v)=h_T(w)$.
 \end{defn}
 
 We can now define more precisely the meaning of the labels $\text{\sout{$T_v$}}$
 as follows.

  \begin{defn}\label{tracelabeldef}
 Let $T\in \fT_{\cS\cO_0}$ be a syntactic object endowed with a head function $h_T$,
 and let $T_v \subset T$ be an accessible term. If $T_v=\ell$ is a single leaf, decorated
 by an element $\alpha_\ell \in \cS\cO_0$, then we simply define $\text{\sout{$T_v$}}$ as
 the label $\text{\sout{$\alpha_\ell$}} \in \text{\sout{ $\cS\cO_0$ }}$. If $T_v$ is a subtree
 containing more than a single leaf, 
 the label $\text{\sout{$T_v$}}$ is defined as the element $\text{\sout{$\alpha_{h_T(v)}$}}$
 in $\cS\cO_0$,  where $\ell=h_T(v)$ is the value of the head function $h_T$ 
 at the vertex $v\in V^o(T)$ and $\alpha_{h(T_v)}$ is the labeling of this leaf in $\cS\cO_0$. 
 \end{defn}
 
 We need to verify that this definition of the labeling $\text{\sout{$T_v$}}$ is 
 compatible with the coassociativity of the coproduct $\Delta^c$.
 
 \begin{lem}\label{cancellabelcoprod}
The coproduct $\Delta^c$ with the leaves-labels $\text{\sout{$T_{v_i}$}}$ of the 
contracted accessible terms $T_{v_i}$ given by $\text{\sout{$\alpha_{h(T_v)}$}}$
with $h=h_T$ a head function is coassociative.
\end{lem}

\proof 
As shown in \cite{MCB} a head function partitions the vertices of $T$ into
disjoint paths $\gamma_\ell$ for $\ell\in L(T)$, where for all $\gamma_\ell$
consisting of more than the single vertex $\ell$, the other end-vertex $v_\ell$
is called the ``maximal projection" of $\ell$. For each vertex $v$ in $\gamma_\ell$
one has $h_T(v)=\ell$. The leaf $\ell$ with the property that
its maximal projection is the root of $T$ is called the head of $T$, and denoted by
$\ell = h(T)$. The partitioning of vertices into paths $\gamma_\ell$ determines
a labeling of all vertices of $T$, where vertices on the path $\gamma_\ell$
acquire the same label $\alpha_\ell$ originally assigned to the leaf $\ell$. 
Thus, the labeling $\text{\sout{$\alpha_{h(T_v)}$}}$ assigns to the
root vertex $v$ of $T_v$ the label $\alpha_{h_T(v)}$ with the additional
cancellation symbol $\text{\sout{$\alpha_{h_T(v)}$}}$ that indicates the new 
leaf $v$ is a ``trace" of movement. What we need to check is that this
assignment of labels maintains the coassociativity of the coproduct $\Delta^c$.
The coassociativity $(\Delta^c \otimes {\rm id}) \circ \Delta^c= ({\rm id}\otimes \Delta^c)\circ \Delta^c$
is an identity involving two successive admissible cuts $C,C'$, where
the first cut $C$ corresponds to the first application of the coproduct $\Delta^c$
and the second cut $C'$ can take place in either the left or the right side
($\Delta^c \otimes {\rm id}$ or ${\rm id}\otimes \Delta^c$) 
of the output of the first $\Delta^c$. Without leaf-labels, the coassociativity
of this form of coproduct is well known, so the issue is only whether the labeling
is consistent on the two sides of the coassociativity identity. Since the
coproduct is multiplicative, we only need to check this on a tree $T$, rather
than on an arbitrary forest $F$. For an admissible cut $C$ we use the
notation $\pi_C(T)=F_{\underline{v}}$ for the forest extracted by the cut,
which is in the left side of the corresponding terms of the coproduct, and
we write $\rho^c_C(T)=T/^c F_{\underline{v}}$ for the quotient term on
the right side, so the term of the coproduct determined by the cut $C$ is
$\pi_C(T) \otimes \rho^c_C(T)= F_{\underline{v}} \otimes T/^c F_{\underline{v}}$,
and the full coproduct is $\Delta^c(T)=\sum_C \pi_C(T) \otimes \rho^c_C(T)$,
including the primitive part $T\otimes 1 + 1 \otimes T$ in the form of the two
trivial cuts (above the root, below the leaves). 
The coassociativity relation compares terms of the form
$$ \pi_{C'}(\pi_C(T))\otimes \rho^c_{C'}(\pi_C(T)) \otimes \rho^c_C(T) $$
in $(\Delta^c \otimes {\rm id})(T)$ with terms of the form
$$ \pi_{C'}(T) \otimes \pi_C(\rho^c_{C'}(T)) \otimes \rho^c_C(\rho^c_{C'}(T)) $$
in $({\rm id}\otimes \Delta^c)(T)$, where in this second term $C'$ is
the first cut and $C$ the second. The terms
$\pi_{C'}(\pi_C(T))=\pi_{C'}(T)$ (since the cut $C'$ is placed below $C$)
do not involve any trace labeling so nothing needs to be checked.
The middle term $\rho^c_{C'}(\pi_C(T))=\pi_C(\rho^c_{C'}(T))$ has
trace labels at the vertices $w_j$ that are the roots of the
trees in $\pi_{C'}(T)=\sqcup_j T_{w_j}$. These are labelled as 
$\text{\sout{$\alpha_{h_T(w_j)}$}}$ in $\rho^c_{C'}(T))$, and maintain
the same label in $\pi_C(\rho^c_{C'}(T))$, while they are labelled as
$\text{\sout{$\alpha_{h_{T_{v_i}}(w_j)}$}}$
in $\rho^c_{C'}(\pi_C(T))$, where $\pi_C(T)=\sqcup_i T_{v_i}$.
But the properties of the head function ensure that
$h_{T_{v_i}}(w_j)=h_T(w_j)$
whenever $T_{w_j}\subset T_{v_i}\subset T$, so the
labels match. Similarly, in the case of the terms
$\rho^c_C(T)=\rho^c_C(\rho^c_{C'}(T))$, the trace labels are
at the vertices $v_i$ that were the root vertices of the $T_{v_i}$
with $\pi_C(T)=\sqcup_i T_{v_i}$. In both $\rho^c_C(T)$ and
$\rho^c_C(\rho^c_{C'}(T))$ these are labelled by 
these are labelled by $\text{\sout{$\alpha_{h_T(v_i)}$}}$
so again there is no labeling mismatch.
\endproof

This labeling method for the traces of movement require that all
the components $T_i$ of all the workspaces $F=\sqcup_i T_i$ 
have an assigned head function $h_{T_i}$. On a tree $T_i$ with
$k_i=\# L(T_i)$ leaves there are $2^{k_i-1}$ possible choices
of a head function. Here we simply assume that one such choice
is made. More realistically, one wants a head function that correctly
models the syntactic head and the structure of the Extended Projection,
as discussed in \cite{MHL}. We will return to discuss this in \S \ref{ColorSec}. 

 \subsection{Merge operations with the contraction coproduct}\label{contractMergeSec}
 
 It is natural then to ask what happens to the Merge dynamics
 $\cK=\sum_{S,S'} \fM^c_{S,S"}$, acting on workspaces in $\fF_{\cS\cO_0 \sqcup \text{\sout{ $\cS\cO_0$ }}}$, 
 if we define $\fM^c_{S,S'}$ as 
 \begin{equation}\label{McSS}
  \fM^c_{S,S'} = \sqcup \circ (\cB \otimes {\rm id}) \circ \delta^c_{S,S'} \circ \Delta^c 
 \end{equation}  
 instead of the Merge operators considered in the previous sections, with
 $$ \fM_{S,S'} = \sqcup \circ (\cB \otimes {\rm id}) \circ \delta_{S,S'} \circ \Delta^d \, .  $$

 It is important to point out that, in passing from the coproduct $\Delta^d$
 to $\Delta^c$, one is also modifying the selection $\delta^c_{S,S'}$ of which
 terms of the coproduct are used by Merge. (This is reflected in \cite{MCB}
 in the counting of the effect of Merge on the number of accessible terms,
 though it is not spelled out explicitly, so we make it more explicit here.)
 In linguistics, the trace of previous movement is no longer extracted by
 successive applications of Merge in a derivation for further movement: in
 other words, traces no longer move. To reflect this fact, the selection $\delta_{S,S'}$
 of accessible terms 
 (in the case of the coproduct $\Delta^d$), that is defined in \cite{MCB} 
 simply as $\delta_{S,S'}=\gamma_{S,S'}\otimes {\rm id}$ where
 $$ \gamma_{S,S'}(F)=\left\{ \begin{array}{ll} F & F=S\sqcup S' \\
 0 & \text{otherwise,} \end{array}\right. $$
 can be defined in the slightly modified form 
 $$ \delta^c_{S,S'}=\gamma^c_{S,S'}\otimes \tilde\gamma^c_{S,S'} \, , $$ 
 where
 \begin{equation}\label{deltaSSc}
  \gamma^c_{S,S'}(F)=\left\{ \begin{array}{ll} F & F=S\sqcup S' \,\, \text{ and } S,S'\notin  \fT_{\text{\sout{ $\cS\cO_0$ }}} \\
 0 & \text{otherwise,} \end{array}\right. 
 \end{equation}
 and, for $F=T_1\sqcup\cdots\sqcup T_r$,
 \begin{equation}\label{deltaSSc2}
  \tilde\gamma^c_{S,S'}(F)=\left\{ \begin{array}{ll}
  F & T_i \notin \fT_{\text{\sout{ $\cS\cO_0$ }}}, \, \, \forall i=1,\ldots, r \\
    0 & \text{otherwise.}
   \end{array}\right. 
 \end{equation}
 
 Note that $\gamma^c_{S,S'}$ does not prevent the extraction of any term
 $S=T_w$, where $T_w\in \fT_{\cS\cO_0\sqcup \text{\sout{ $\cS\cO_0$ }}}$
 has some of the leaves decorated by elements in $\text{\sout{ $\cS\cO_0$ }}$:
 it only prevents extraction of accessible terms that have {\em all} leaves 
 labelled by trace labels in $\text{\sout{ $\cS\cO_0$ }}$. 
 
 \smallskip
 
 The operator $\tilde\gamma^c_{S,S'}$
 correspondingly prevents the creation of components in the right-hand-side of the
 coproduct that have only trace-leaves.  The reason for including here also
 the constraint $\tilde\gamma^c_{S,S'}$ is in order to better relate the action of
 Merge with the two forms of the coproduct. For example, an Internal Merge
 transformation involving a term $T/^c T_v \in \cT_{\text{\sout{ $\cS\cO_0$ }}}$
 would have no analog in terms of the coproduct $\Delta^d$, hence in linguistic terms
 it would not be viable at Externalization. 
 
 \smallskip
 
 The presence of both constaints
 $\gamma^c_{S,S'}$  and $\tilde\gamma^c_{S,S'}$ ensures that all the Merge
 transformation that can take place in $\cK^c=\sum_{S,S'} \fM^c_{S,S'}$ also
 have an associated form in $\cK=\sum_{S,S'} \fM_{S,S'}$, expressing the
 idea that the same mechanism of structure formation always simultaneously
 maps to the Externalization interface (where the traces are not expressed) and to
 the Conceptual-Intentional (also called Syntax-Semantics) interface where the
 traces are present for interpretation. 
  
 \smallskip
 
 We address in the rest of this section the behavior of the dynamics of Merge in the form
 based on the coproduct $\Delta^c$. 
 
 \subsection{Edges of the Merge graph with the contraction coproduct}\label{traceaccSec}
 
 The arrows in the graph of the dynamical system obtained using the Merge operators as 
 defined as in \eqref{McSS} again correspond to External Merge, Internal Merge,
 and Sideward Merge operations. 
  
 \begin{enumerate}
 \item The External Merge arrows are exactly the same considered in the previous sections: 
 they remain unchanged with respect to the choice of the $\Delta^d$ or $\Delta^c$
 coproduct, since they only involve the primitive part of the coproduct that partitions
 the workspace into components, without extracting accessible terms. Since in this
 case no quotient is involved, both coproducts have the same primitive part.
 
 \item The Internal Merge arrows correspond to compositions $\fM^c_{T_v, T/^c T_v} \circ \fM_{T_v, 1}$,
 hence they differ from the case with the coproduct $\Delta^d$, due to the presence of
 the additional leaf with label $\text{\sout{$T_v$}}$ in $T/^c T_v$
 
 \item The Sideward Merge arrows also similarly differ from the case of the coproduct $\Delta^d$.
 Since again we are only considering the minimal SM arrows that extract only leaves, we
 have cases 
 \begin{itemize}
\item $\fM(\beta, \alpha) \sqcup T/^c \alpha$ (when $\beta$ and $T$
 are components of the workspace SM acts on), where the resulting component $T/^c \alpha$ 
 now still have the leaf $\alpha$ of $T$ but with label changed to $\text{\sout{$\alpha$}}$,
 \item $\fM(\beta, \alpha) \sqcup T/^c (\alpha \sqcup \beta)$, where the quotient
  $T/^c (\alpha \sqcup \beta)$ has the same leaves as $T$, but those that were respectively 
  labelled $\alpha$ and $\beta$ are now relabelled $\text{\sout{$\alpha$}}$ and
  $\text{\sout{$\beta$}}$, 
 \item $\fM(\beta, \alpha) \sqcup T/^c \alpha \sqcup T'/^c \beta$ where similarly
 the resulting components $T/^c \alpha$ and $T'/^c \beta$ maintain the same 
 number of leaves as $T$ and $T'$, with $\alpha$ and $\beta$ labels
 changed to $\text{\sout{$\alpha$}}$ and
  $\text{\sout{$\beta$}}$, respectively.
  \end{itemize}
  \end{enumerate}
 
 In all cases involving a quotient $T/^c T_v$ (or $T/^c F_{\underline{v}}$), one
 should interpret the labels $\text{\sout{$T_v$}}$ as in Definition~\ref{tracelabeldef}.
 
 \smallskip
 
 \begin{defn}\label{contractMergeGraph}
 We define the contraction Merge graph $\cG_{\fF_{\cS\cO_0 \sqcup \text{\sout{ $\cS\cO_0$ }}}}$
 as the infinite, locally finite graph with countably infinite set of vertices 
 $\fF_{\cS\cO_0\sqcup \text{\sout{ $\cS\cO_0$ }}}$ and with finitely many
 incoming and outgoing edges at each vertex given by all the possible 
 EM, IM, SM transformations listed above. 
 \end{defn}

 \subsection{Transient dynamics and projection}\label{transientKcSec}
 
 We show here that looking at the Merge dynamics in terms of the
 coproduct $\Delta^c$ does not provide any useful information, since
 the dynamics now involves a countably infinite set of states (it is
 no longer a finite Markov chain) where all states become transient, so 
 in particular there is no stationary distribution that describes the long
 term behavior of the dynamics. However, we also show that the Merge
 dynamics in the $\Delta^d$ case provides complete information to also
 determine the dynamics of the $\Delta^c$ case. The key reason for
 this lies in the fact that, in linguistics, traces of previous movement 
 are no longer extracted for further movement.
 
 \begin{prop}\label{alltransit}
 Let $\cK^c=\sum_{S,S'} \fM^c_{S,S'}$ be the action of Merge on the linear span of
 workspaces $\cV(\fF_{\cS\cO_0\sqcup \text{\sout{ $\cS\cO_0$ }}})$, with $\fM^c_{S,S'}$
 as in \eqref{McSS}.
 This action does not preserve the subspaces spanned by the finite sets $\fF_{A,n}$, and
 all $F\in \fF_{\cS\cO_0\sqcup \text{\sout{ $\cS\cO_0$ }}}$ with $F\notin \fF_{\text{\sout{ $\cS\cO_0$ }}}$ 
 are transient states for the dynamics defined by $\cK^c$.
 \end{prop}
 
 \proof As discussed in \S \ref{traceaccSec}, the EM arrows in $\cG_{\fF_{\cS\cO_0}}$
 are the same as in $\sqcup_n \cG_{n,A}$ so they do preserve the $\fF_{A,n}=V(\cG_{n,A})$. 
 However, this is not the case for IM and SM arrows. Indeed, in the case of IM arrows, if 
 $F=T\sqcup \hat F \in \fF_{A,n}$,
 its image under an IM transformation that uses the form $\Delta^c$ of the coproduct and acts
 on the $T$ component is a new workspace 
 $$ F' =\fM(T_v,T/^c T_v)\sqcup \hat F \in \fF_{B, n+1} \, ,  \ \ \ \text{ with } \ \ \ 
 B= A\cup \{ \text{\sout{$T_v$}} \} \, $$
 and with the label $\text{\sout{$T_v$}}$ as in Definition~\ref{tracelabeldef}. Similarly,
 in the case of minimal SM, as illustrated in \S \ref{traceaccSec}, we can obtain from
 $F \in \fF_{A,n}$ a new workspace 
 $$ F'=\left\{ \begin{array}{ll}  \fM(\alpha,\beta) \sqcup T/^c \alpha \sqcup \hat F  \,\,\,
 \in \fF_{A\cup\{ \text{\sout{$\alpha$}} \}, n+1} & \text{for } F=\beta \sqcup T\sqcup \hat F \\
  \fM(\alpha,\beta) \sqcup T/^c (\alpha \sqcup \beta)  \sqcup \hat F
  \,\,\, \in \fF_{A\cup\{ \text{\sout{$\alpha$}}, \text{\sout{$\beta$}} \}, n+2}   & \text{for } F=T\sqcup \hat F \\ 
\fM(\alpha,\beta) \sqcup T/^c \alpha \sqcup T'/^c \beta  \sqcup \hat F  \,\,\, \in 
 \fF_{A\cup\{ \text{\sout{$\alpha$}}, \text{\sout{$\beta$}, n+2} \}} & F=T\sqcup T'\sqcup \hat F\, . 
 \end{array}\right.  
 $$
 The resulting graph $\cG_{\fF_{\cS\cO_0 \sqcup \text{\sout{ $\cS\cO_0$ }}}}$ cannot be strongly connected, as the IM
 and SM arrows always map any subset of vertices $\fF_{A,n}$ to a subset of vertices
 $\fF_{B,n+1}$ of $\fF_{B,n+2}$ while the EM arrows map $\fF_{A,n}$ to $\fF_{A,n}$, while
 there are no directed arrows that can decrease the number of leaves. Since every
 $F\in \fF_{A,n}$ has some outgoing IM or SM arrows, this implies that every $F\in \fF_{A,n}$ 
 is a transient state for the dynamical system $\cK^c$.
  \endproof
  
  This shows that, when using the coproduct $\Delta^c$ instead of $\Delta^d$, 
  there is no analog here of the Hopf algebra Markov chain $\hat\cK$
 that we studied in the previous sections and its stationary distribution that governs the
 long term dynamics and provides information on the dynamical behavior. 
  
   \medskip
 
  On the other
 hand, this does not mean that we do not have any viable information about how
 $\cK^c$ behaves. In fact, the following result shows that the dynamics of $\cK^c$,
 acting on workspaces $F\in \fF_{\cS\cO_0}$ is completely determined by its image 
 under the projection $\Pi_{d,c}$ on $\cV(\fT_{\cS\cO_0})$ (as in \S 1.3 of \cite{MCB}).
 
 \begin{defn}\label{Pidc}
 Let $\Pi_{d,c} : \cV(\fF_{\cS\cO_0 \sqcup \text{\sout{ $\cS\cO_0$ }}}) \to \cV(\fF_{\cS\cO_0})$ be
 the projection that cuts from each component $T$ of $F$ all the trace-labelled leaves 
 and then contracts edges to obtain the maximal full binary rooted tree determined by
 what remains. 
 \end{defn}
 
 The iterates of the dynamical systems $\cK=\sum_{S,S'}\fM_{S,S'}$ and $\cK^c=\sum_{S,S'}\fM^c_{S,S'}$,
 based on the two forms of the coproduct, $\Delta^d$ and $\Delta^c$, satisfy the following relation.
 
 \begin{prop}\label{KKc}
  For $F\in \fF_{A,n}\subset  \fF_{\cS\cO_0}$, the orbit $(\cK^c)^m(F)$ under iterations of
 $\cK^c$ is contained in the fiber $\Pi_{d,c}^{-1}( \cK^m(F) ) \subset \cV(\fF_{\cS\cO_0 \sqcup \text{\sout{ $\cS\cO_0$ }}})$
 and is completely determined by the orbit of $F$ under $\cK$. 
 \end{prop}
 
 \proof  We first show that, for $F\in \fF_{\cS\cO_0}$, we have 
 \begin{equation}\label{relKKcPi}
 \Pi_{d,c} \circ \cK^c= \cK \, .
 \end{equation}
 For a forest $F\in \fF_{\cS\cO_0}$, 
 we have
  \begin{equation}\label{dcPicoprod}
({\rm id} \otimes \Pi_{d,c}) \circ \Delta^c(F) = \Delta^d(F)  \, ,
\end{equation}  
or equivalently
$$ \pi_C(T) \otimes \Pi_{d,c} (\rho^c_C(T))= \pi_C(T) \otimes \rho_C(T))\, . $$
Note that, in this case, since $F$ has no trace labels, none of the components
of the terms $F_{\underline{v}}=\pi_C(F)$ extracted from $F\in \fF_{A,n}$ can
be in $\fT_{ \text{\sout{ $\cS\cO_0$ }} }$. Thus,
we have $\delta^c_{S,S'} \Delta^c(T)= \delta_{S,S'}
\Delta^d(F)$ and we obtain in this case
$$ \fM_{S,S'}(F)=\sqcup \circ (\cB \otimes {\rm id}) \circ (\gamma_{S,S'}\otimes \Pi_{d,c}) \Delta^c(F) =
\sqcup \circ ({\rm id}\otimes \Pi_{d,c}) \circ (\cB \otimes {\rm id}) \circ \delta^c_{S,S'} \Delta^c(F)= $$
$$ \sqcup \circ (\Pi_{d,c} \otimes \Pi_{d,c}) \circ (\cB \otimes {\rm id}) \circ \delta^c_{S,S'} \circ \Delta^c(F)=
 \Pi_{d,c} \circ \fM^c_{S,S'} (F) $$  
 so that we get \eqref{relKKcPi}. 
 We now proceed inductively. Suppose that, starting with an $F\in \fF_{\cS\cO_0}$ and
 applying the $m$-th iterate of the action of $\cK^c$ we obtain that $\Pi_{d,c} \circ (\cK^c)^m(F) = \cK^m(F)$.
 Let $\tilde F_i$ be any of the components of $(\cK^c)^m(F)=\sum_i c_i \tilde F_i\in \cV(\fF_{\cS\cO_0 \sqcup \text{\sout{ $\cS\cO_0$ }}})$. When we take $\cK^c(\tilde F_i)=\sum_j c_{ij} \tilde F_{ij}$ each component $\tilde F_{ij}$
has trace-leaves that are either already trace-leaves of $\tilde F_i$ or that are added by the IM and SM 
arrows of $\cK^c$, so that we have
$$ \Pi_{d,c}(\cK^c(\tilde F_i))=\Pi_{d,c}(\cK^c( \Pi_{d,c}(\tilde F_i)) \, . $$
Using \eqref{relKKcPi} we then see that 
$$ \Pi_{d,c}(\cK^c( \Pi_{d,c}(\tilde F_i)) = \cK ( \Pi_{d,c}(\tilde F_i)) $$
hence, by linearity
$$ \Pi_{d,c}((\cK^c)^{m+1}(F)) = \Pi_{d,c}(\cK^c( \Pi_{d,c}((\cK^c)^m(F))) =  \cK ( \Pi_{d,c}( (\cK^c)^m(F))) $$
and by induction hypothesis we then get
\begin{equation}\label{projKKcPi} 
 \Pi_{d,c}((\cK^c)^{m+1}(F)) = \cK^{m+1}(F) \, .
\end{equation}
This implies that the orbit $(\cK^c)^m(F)$ of $F\in \fF_{A,n}$ is contained in
$\Pi_{d,c}^{-1}( \cK^m(F) )$. 
One can also see that the orbit $(\cK^c)^m(F)$ is completely determined by
the orbit $\cK^m(F)$, for $F\in \fF_{\cS\cO_0}$, because at each application
of $\cK^c$ the only additional information that is needed to determine the image
is the trace-leaves created by the coproduct $\Delta^c$. These are determined
by knowing the vertices $v$ where the accessible terms $T_v$ are extracted
and the value of the head function $h_T(v)$ in the component $T \supset T_v$. 
This can also be shown inductively as above. At the first step, the extracted
$T_v$ are the same for $\cK^c$ and $\cK$ since the orbit starts at a point is $F\in \fF_{\cS\cO_0}$.
Thus, comparing $F$ with the components of $\cK(F)$ identifies the terms $T_v$ involved
in each arrow. We are assuming here that all the components of the workspaces involved 
have a head function (we will return to discuss this in \S \ref{ColorSec}). So 
for $T_v\subset T$, by the properties of the head function this satisfied $h_T(v)=h_{T_v}(v)=h_{T'}(v)$ where $T'$ is the
component of $\cK(F)$ that contains the vertex $v$, so this is computable from 
$\cK(F)$. Inductively, assume that in $(\cK^c)^m(F)$ the trace-leaves and their labels
computed through the head function are all determined from data computable from
the orbit $\{ \cK^k(F) \}_{k\leq m}$. The next iterate $\cK^c((\cK^c)^m(F))$ does not
change the previous trace-leaves and labels but adds new ones through the IM and SM
arrows in $\cK^c$. Again these depend on the extracted accessible terms and the
value of the head function at their root vertices. These, as in the case of a single
iteration only depend on data in $(\cK^c)^m(F)$ that are determined by 
the orbit $\{ \cK^k(F) \}_{k\leq m}$.
\endproof 

Thus, the conclusion is that the Merge dynamics is fully described by the
Hopf algebra Markov chains $\hat\cK^{A,n}$ that we studied in the previous
sections and the change of the coproduct from $\Delta^d$ to $\Delta^c$ does
not provide any additional information. 

\smallskip

In terms of the linguistic model, the result of Proposition~\ref{KKc} can
be interpreted in the following way. The dynamics of $\cK^c$ determines
(via the projection $\Pi_{d,c}$) and is in turn determined by the dynamics of $\cK$
(via the argument given in Proposition~\ref{KKc}). This means that, while the
presence of traces is considered necessary for semantic parsing, it is not
necessary to the structure building dynamics itself, which can be fully
determined without using trace-leaves. As we have argued, this is a consequence of the fact that,
in linguistics, traces are no longer extracted for movement. 

\section{Weighted Merge dynamics: cost functions}\label{CostSec}

The conclusions we have drawn from the previous sections agree with
what is stated in the linguistics literature, in the formulation of the action
of Merge on workspaces of \cite{ChomskyUCLA} and \cite{ChomskyGK}.
Namely, some natural cost function must weight the different forms of
Merge differently, with EM and IM being optimal and SM being sub-optimal.

\smallskip

As proposed in \cite{ChomskyUCLA} and \cite{ChomskyGK}, there are
two main forms of optimality that Merge is expected to satisfy. One is what
is called Minimal Search, which optimizes the cost of searching in the
workspace for the accessible terms that Merge applies to. The other goes
under the name of Resource Restrictions, and especially Minimal Yield and
is formulated as a balance between the effect of a Merge
transformation on the number of connected components of the workspace
and on the number of accessible terms. With these two forms of optimality,
one can include an additional one, which is referred to in \cite{MCB} as
``no complexity loss", and which is usually bundled up, in the linguistics
literature, together with the ``Extension Condition": it requires that structures
are {\em grown} by Merge (with the second part of the Extension Condition
also requiring that growth happens at the root, as discussed at length in \cite{MLH}).

\smallskip

The Minimal Yield condition, in turn, has slightly different formulations in
\cite{MCB} (and in \cite{MLH}) from the original formulation of  \cite{ChomskyUCLA} 
and \cite{ChomskyGK}.
We will recall briefly these cost functions here, in \S \ref{3Costs}, following the formulation 
as in \cite{MLH}. Using the original form of the cost function does not affect the main
results stated here, so we adopt the formulation of \cite{MCB} and \cite{MLH}, since
it is more immediately transparent in terms of the tree topology.

\smallskip

We then describe in \S \ref{EstPFsec} a technique for estimating the leading order
in the weight parameter of the Perron--Frobenius eigenvalue and eigenvector of
the cost-weighted Merge dynamics, and of the stationary distribution of the
associated Hopf algebra Markov chain. 

\smallskip

The algebra of orders of magnitude is governed by the min-plus (also known
as tropical) semiring. Thus, the estimation of the leading order of the 
Perron-Frobenius data for the cost-weighted Merge is phrased in terms of
a Perron-Frobenius problem in the tropical semiring. 

\smallskip

A general approach to such tropical Perron-Frobenius problems reduces
evaluating the order of the Perron-Frobenius eigenvalue  to identifying critical 
circuits (cost-minimal circuits) in the Merge graph,
and the construction of a basis of the Perron-Frobenius eigenspace to
the construction of optimal arborescences for the Merge graph. We discuss
the details of this approach and the result in \S \ref{EstPFsec}.

\smallskip

The conclusion of this analysis is that, due to the non-uniqueness of 
critical circuits, we do not identify a unique Perron-Frobenius eigenvector,
but a basis for the eigenspace. The form of this basis shows that, 
to leading order, the weighted dynamics moves away from the more disconnected
structures that dominate in the unweighted one as an effect of SM. On
the other hand, instead of convergence to the tree structures (i.e.~the 
partition $\wp=\{ n \}$), as observed in \cite{MLH} in the weighted dynamics
for $n=3$,  one finds that any of the partitions $\wp\in \cP'(n)$ that have
at least some $k_i\geq 3$ (hence that have a non-trivial IM self-loop) 
contribute to the Perron-Frobenius eigenspace.

\smallskip

This shows that, despite what the simpler case $n=3$ described in \cite{MLH}
suggests, the total cost function that combines Minimal Search (in a bottom-up
formulation as in \cite{MLH}), Minimal Yield (as formulated in \cite{MCB}) and
Complexity-Loss, does not suffice to fully eliminate the unwanted effects of
Sideward Merge on the dynamics (while maintaining the desirable strong 
connectedness properties)  and ensure convergence to the tree structure. 

\smallskip

We describe how to correct this problem by modifying the cost
function and introducing in it a term that depends on the Shannon
entropy of a distribution $\P_\wp$ associated to the partition $\wp\in \cP'(n)$,
which simply describes the leading order of the combinatorial 
multinomial coefficient describing the partitioning of the $n$ leaves
into the different tree components of the workspace. We show that
the presence of this additional quantity ${\rm Sh}(\P_\wp)$ in the
cost function suffices to drive the dynamics toward the connected 
workspaces (consisting of a single tree). Namely, we show that
the stationary distribution of the Merge Markov chain weighted by
$t^{\fc}$ for this modified cost function $\fc$ converges, in the
limit $t\to 0$ to the uniform distribution on the trees. In other
words, the $\hat\cK^{(A,n)}$ dynamics on $\cG_{n,A}$, weighted
by this modified cost function, converges to  
the IM dynamics on $\cG^{\rm IM}_{n,A, \{ n \}}$. 
The leading order in $t$ of all the components of the stationary
distribution can be expressed explicitly in terms of costs of
certain loops of edges in the graph $\cG_{n,A}$ that lie along 
optimal arborescences in $\cG_{n,A}$.

\smallskip

In fact we show an even stronger fact: that we can use {\em only}
the cost function given by the Shannon entropy ${\rm Sh}(\P_\wp)$ 
and obtain the same result, without including the 
Minimal Search, Minimal Yield, and Complexity-Loss cost functions
at all. These cost functions, usually considered in linguistics, then
play a role only in the initial selection of only the ``minimal" 
Sideward Merge arrows that involve only atomic components,
while the influence on the dynamics is entirely coming from
the information-theoretic cost function.

\subsection{Combined cost functions}\label{3Costs}

We consider here the three main types of cost function as discussed
in \cite{MLH}, namely a bottom-up version of Minimal Search, and
two types of Resource Restriction costs: Minimal Yield and Complexity Loss.

\subsubsection{Minimal Search costs}\label{MScost}
As in \cite{MLH}, we consider a ``bottom-up" version of the Minimal
Search cost (rather than the ``top-down") version used in \cite{MCB}.
In this approach to Minimal Search we assume that, instead of the
usual description of search algorithms on rooted trees that start at
the root and descend into the tree, the search on syntactic
structures, which are built in a bottom up way, starts with the
lexical items at the leaves and searches for subsets of leaves
that form a constituency, namely that determine an accessible
term that can be extracted for use by the Merge operations.
Viewed in this way, the cost of the search depends on the
(relative) size of the set of leaves that needs to be identified.
This leads to a definition of cost function where
$$ \fc_{MS}(T_v)=\frac{\ell(T_v)}{\ell(T)}\, , $$
where $\ell(T)=\# L(T)$ is the number of leaves (the grading
in the Hopf algebra of workspaces). This gives a cost to the
extraction of an accessible term. The cost of the corresponding
quotient operation $T/T_v$ is then set to be
$$ \fc_{MS}(T/T_v)=\frac{\ell(T/T_v)}{\ell(T)} = 1- \frac{\ell(T_v)}{\ell(T)}= 1-\fc_{MS}(T_v)\, . $$
As in \cite{MCB} we also assign a resulting cost to the merging of two components or
accessible terms $A,B$ of the form
$$ \fc_{MS}(\fM(A,B))= \fb(A,B) -\fc(A)-\fc(B) \, , $$
where $\fb(A,B)=1$ if $A$ and $B$ come from the same component of the workspace
and $\fb(A, B)=2$ if they come from different components. This choice
encodes the idea that Minimal Search would first (lower cost) search for terms within a fixed
components and then (higher cost) for terms across different components. With this
cost function we see that we have for the different cases of Merge, 
EM $\fM(T,T')$, IM $\fM(T_v, T/T_v)$, SM(1) $\fM(T_v, T') \sqcup T/T_v$, 
SM(2) $\fM(T_w,T'_w)\sqcup T/T_v \sqcup T'/T'_w$, and SM(3) $\fM(T_v,T_w)\sqcup T/(T_v\sqcup T_w)$ the respective costs
\begin{center}
\begin{tabular}{|c|c|c|c|c|c|}
\hline 
                   & EM & IM & SM(1) & SM(2) & SM(3) \\
\hline
$\fc_{MS}$ &  $0$  & $0$ & $2-\fc_{MS}(T_v)$ & $2-  \fc_{MS}(T_v)- \fc_{MS}(T'_w)$ & $1-\fc_{MS}(T_v)-\fc_{MS}(T_w)$ \\
 \hline 
\end{tabular}
\end{center}

In particular, in the case of a tree $T$ with $k=\# L(T)$, if we only consider the minimal 
SM transformations, where the terms $T_v$, $T_w$, $T'_w$, $T'$ involved in the operations
of extracting and Merging are all atomic components consisting of a single leaf, we obtain the
costs
\begin{center}
\begin{tabular}{|c|c|c|c|c|c|}
\hline 
                   & EM & IM & min SM(1) & min SM(2) & min SM(3) \\
\hline
$\fc_{MS}$ &  $0$  & $0$ & $2-1/k_i$ & $2-  1/k_i -1/k_j$ & $1- 2/k_i$ \\
 \hline 
\end{tabular}
\end{center}
This means that, on a workspace $F$ with $n=\# L(F)$ partitioned across
components $F=T_1\sqcup\cdots \sqcup T_r$ with $n=k_1+\cdots +k_r$, $k_i=\# T_i$,
an SM transformation that extracts from a component $T_i$ (for SM(1) and SM(3)) or to two
components $T_i$ and $T_j$ (for SM(2)), has cost, $2-1/k_i$, $2-1/k_i-1/k_j$ and $1-2/k_i$,
for SM(1), SM(2), and SM(3), respectively.

\subsubsection{Minimal Yield costs}\label{MYcost}
For the computation of Minimal Yield costs we proceed as in \cite{MCB} and \cite{MLH} and 
we evaluate the effect of the different Merge operations (EM, IM, SM(1), SM(2), SM(3)) on
the number $b_0(F)$ of connected components  and the number $\alpha(F)$ of accessible terms of the workspace,
evaluating the change in the total number $\sigma(F)=b_0(F) + \alpha(F)$. 
We consider here only the case of the coproduct $\Delta=\Delta^d$ as we have argued
that this suffices for studying the properties of Merge as a dynamical system. 
We count for each type of Merge operation the difference between the number of
components of the resulting workspace and of the initial one. This is equal to $-1$ for EM
(which reduces the number of connected components), $0$ for IM and for SM(1), and $+1$
for SM(2) and SM(3) that increase by one the components. We also count the resulting change
on the number of accessible terms: EM increases it by $2$, IM and SM(1) leave it unchanged
and SM(2) and SM(3) reduce it, changing it by $-2$ (each quotient $T/^dT_v$ involves removal of
one vertex to obtain a full binary tree). Thus the total count of the change in the sum
of number of components and accessible terms gives a total Minimal Yield counting of
\begin{center}
\begin{tabular}{|c|c|c|c|c|c|}
\hline 
                   & EM & IM & SM(1) & SM(2) & SM(3) \\
\hline
$\fc_{MY}$ &  $1$  & $0$ & $0$ & $-1$ & $-1$  \\
 \hline 
\end{tabular}
\end{center}
This counting is independent of whether SM is minimal (extracting and merging only atomic components)
or not.

\subsubsection{Complexity Loss cost} \label{CLcost} 

Complexity loss is another possible measure of costs, discussed in \cite{MCB}. It encodes one
of the properties of the linguistic model. This is sometime grouped together with the ``Extension
Condition" in the linguistic literature, but it is better treated independently, because, as argued in 
\cite{MLH}, the Extension Condition is really an ``algebraic" constraint (hard constraint) of the
model, while this is an optimality condition (soft constraint). One can formulate it as the 
requirement that ``syntactic composition always grows or expands structures". More
precisely, we can say that there is no complexity loss if in the transformed workspace, after
the action of a Merge operation, each component of the original workspace ends up as part
of a component that is at least as complex as the original one. A simple way to measure the
complexity of a  component of the workspace is by its size, measured by its degree in the
Hopf algebra of workspaces (the number of leaves). Thus the complexity loss is measured
as the largest drop in degree over the components of the workspace. Thus, for an EM of the
form $\fM(T,T')$ the components $T$, $T'$ end up as parts of a component $\fM(T,T')$ of
largest degree (and any other components are unchanged), so there is no loss. For an IM
of the form $\fM(T_v, T/T_v)$ the original component $T$ and the resulting $\fM(T_v, T/T_v)$
have the same degree and the other components are unchanged, so also there is no complexity
loss. For an SM(1) of the form $\fM(T',T_v)$, the component $T'$ ends up in a 
larger component (so it does not incur any complexity loss) but the component $T$ from which $T_v$
is extracted results in a component $T/T_v$ of lower degree, with a complexity loss of $\ell(T_v)$.
For SM(2) there are two components $T$ and $T'$ that result, respectively, in components $T/T_v$
and $T'/T'_w$ for a loss of $\ell(T_v)+\ell(T'_w)$ and for SM(3) there is a component $T$ resulting
in $T/(T_v\sqcup T_w)$ for a loss of $\ell(T_v)+\ell(T_w)$. So the CL cost function is given by
\begin{center}
\begin{tabular}{|c|c|c|c|c|c|}
\hline 
                   & EM & IM & SM(1) & SM(2) & SM(3) \\
\hline
$\fc_{CL}$ &  $0$  & $0$ & $\ell(T_v)$ & $\ell(T_v)+\ell(T'_w)$ & $\ell(T_v)+\ell(T_w)$  \\
 \hline 
\end{tabular}
\end{center}
For the case of minimal SM operations involving only atomic (single leaf) elements, this gives
\begin{center}
\begin{tabular}{|c|c|c|c|c|c|}
\hline 
                   & EM & IM & min SM(1) & min SM(2) & min SM(3) \\
\hline
$\fc_{CL}$ &  $0$  & $0$ & $1$ & $2$ & $2$  \\
 \hline 
\end{tabular}
\end{center}

\subsubsection{Combined cost}\label{CombCostSec}
Thus, following \cite{MLH}, the costs according to these different cost functions
for the Merge operations are summarized as follows. For Merge acting 
on components $T_i$, or $T_i$ and $T_j$, of a forest with $n=k_1+\cdots+k_r$ leaves,  
one has weights $\fc_{\rm tot}$ given by:
\begin{center}
\begin{tabular}{|c|c|c|c|c|}
 \hline 
        & MS & MY & CL & total weight $\fc_{\rm tot}$ \\
        \hline 
 EM  &  $0$    &  $1$  &  $0$  & $1$ \\
 \hline 
 IM   &   $0$    &  $0$  &  $0$  & $0$ \\
  \hline 
 SM(1) & $2-1/k_i$ & $0$ & $1$  & $3-1/k_i$ \\
  \hline 
 SM(2) & $2-1/k_i -1/k_j$ & $-1$ & $2$ &  $3-1/k_i - 1/k_j$  \\
  \hline 
 SM(3)  & $1-2/k_i$ & $-1$ & $2$ & $2-2/k_i$ \\
  \hline 
\end{tabular}
\end{center}

With this counting of the total cost function of Minimal Search, Minimal Yield, and 
Complexity Loss, we consider a weighted version of the Merge dynamics. We
think of the unweighted Merge dynamics discussed so far as the value $t=1$
of a deformed dynamics, dependent on a deformation parameter $0< t \leq 1$,
which recovers the unweighted case at $t=1$ and, for $t\to 0$ has leading
order given by the least costly part of the dynamics. This is achieved by weighting
each component of the Merge graph adjacency matrix by a cost function $t^{\fc_{\rm tot}}$.
We analyze in the next section the effect of this deformation on the properties 
of the resulting dynamical system.

\section{Tropical Perron-Frobenius and eigenvector estimates}\label{EstPFsec}

We have seen in Theorem~\ref{MergeDynStatDistr} that the stationary distribution for
the Merge Hopf algebra Markov chain is given by
$$ \pi(F) = \frac{\eta_{\cP'(n)}(\wp) \, \, \psi_{\cP'(n)}(\wp)}{\Lambda_{\wp,n}} \ \ \ \text{ for } \wp=p(F) $$
where $\eta_{\cP'(n)}$ and $\psi_{\cP'(n)}$ are, respectively, the left and right Perron--Frobenius eigenvectors
of the adjacency matrix of the graph $\cG_{\cP'(n)}$. In the case where we consider the weighting by
the cost functions, the same result holds for the weighted Merge Hopf algebra Markov chain dynamics.
Namely, we have the following.

\begin{prop}\label{MergeDynStatDistrW}
Consider the directed Merge graph $\cG_{n,A}$, where the arrows
are weighted by a weight $t^{\fc_{\rm tot}}$, according to the total cost $\fc_{\rm tot}$ 
(combining MS, MY, CL) as computed in \S \ref{CombCostSec}, and with $t\in \R^*_+$ a
parameter. We write $\cK^{A,n,w}(t)$ for the resulting non-negative matrix with these weights,
\begin{equation}\label{KtFF}
 \cK^{A,n,w}(t)_{F,F'}=t^{\fc_{{\rm tot},F,F'}}\,\, \cK^{A,n}_{F,F'} \, 
\end{equation} 
seen as a function of the parameter $t$. 
Consider again the projection map $p:\cG_{n,A} \to \cG_{\cP'(n)}$, where in the 
graph $\cG_{\cP'(n)}$ we consider the EM and SM arrows weighted by $t^{\fc_{\rm tot}}$, and the
IM loops weighted by $t^{\fc_{\rm tot}}=1$. 
We write $\cK^{\cP'(n),R,w}(t)$ for this resulting weighted adjacency matrix of $\cG_{\cP'(n)}$.
Then the stationary distribution of the weighted Merge Hopf algebra Markov chain $\hat\cK^{A,n,w}(t)$ 
determined by this weighted $\cK^{A,n,w}(t)$ has stationary distribution 
\begin{equation}\label{piKt}
 \pi^w(F,t) = \frac{\eta_{\cP'(n),w}(\wp,t) \, \, \psi_{\cP'(n),w}(\wp,t)}{\Lambda_{\wp,n}} \ \ \ \text{ for } \wp=p(F) \, , 
\end{equation} 
where $\eta_{\cP'(n),w}$ and $\psi_{\cP'(n),w}$ are the left and right Perron--Frobenius eigenfunctions of
$\cK^{\cP'(n),R,w}(t)$.
\end{prop}

\proof The argument is as in Proposition~\ref{propPFpi}, but taking into account the weights at the
edges. All the IM edges in the $\cG_{n,A}$ have $\fc_{\rm tot}=0$, so all the edges in the fibers
of the projection $p$ do not carry any weight. We have the same numerical factors $\Lambda_{\wp,n}=\# V_{n,A,\wp}$
in the fiber multiplicities at the vertices.
\endproof

Note that the form \eqref{KtFF} is analogous to the weighted adjacency matrices that we
discussed in \S \ref{FreeEnSec}, after 
setting $t=e^{-\beta}$ for the thermodynamic parameter of the Boltzmann distribution.

The form \eqref{piKt} of the stationary distribution reduces, as in the unweighted case, the problem to the computation
of the left and right Perron-Frobenius eigenfunctions of the matrix $\hat\cK^{\cP'(n),R,w}(t)$ for the weighted
graph $\cG_{\cP'(n)}$.

While the graph $\cG_{\cP'(n)}$ is much smaller than the Merge graph $\cG_{n,A}$, it is still 
computationally difficult to obtain a closed expression (as a function of the parameter $t$) for
the Perron-Frobenius eigenvalue and eigenvector, as done in the simple case of $n=3$ in \cite{MLH}.  
This means that, in the general case of arbitrary $n\geq 4$, one needs to obtain an estimate of the
leading order behavior in $t$ for the Perron-Frobenius eigenvalue and left/right eigenvector. 

There are general methods (see for instance \cite{Akian})
for evaluating the leading order of the Perron-Frobenius
eigenvalue and eigenvector of a non-negative matrix that depends on a real parameter,
which is based on solving an associated Perron-Frobenius problem in a semiring  that
provides the appropriate  ``algebra of orders of magnitudes". 

We follow the setting of \cite{Akian}, adapted to our case. In particular, the large real
parameter $p>0$ of \cite{Akian} will be here the parameter $-\log t$, for $0< t <1$. 
Our parameter $t$ is designed to give back the unweighted case discussed in the
previous sections when $t=1$ and give the limiting form of the dynamics, which we
will be discussing in this section, when $t\to 0$, as in the $n=3$ example
analyzed in \cite{MLH}. Thus, our $t\to 0$ limit will indeed correspond to the
$p=-\log t \to \infty$ limit of \cite{Akian}. The matrix denoted by $\cA_p$ in \cite{Akian}
will be our $t$-dependent matrix $\cK^{\cP'(n),R,w}(t)$, which has, by construction,
the following form.

\begin{lem}\label{matrixKt}
The entries of the matrix $\cK^{\cP'(n),R,w}(t)$ are of the form
\begin{equation}\label{Kteq}
 \cK^{\cP'(n),R,w}_{\wp,\wp'}(t) = \cK^{\cP'(n),R}_{\wp,\wp'} \,\, 
t^{\fc_{\wp,\wp'} } \, , 
\end{equation}
where $\cK^{\cP'(n),R}_{\wp,\wp'}$ are the entries of the unweighted 
matrix $\cK^{\cP'(n),R}$, and where 
$$ \fc_{\wp,\wp'}:= \fc_{\rm tot}( \upsilon_{\wp\wp'} )\, , $$
with the label $\upsilon_{\wp,\wp'}$ indicating 
whether the corresponding entry $\cK^{\cP'(n),R}_{\wp,\wp'}$ is an
EM, SM(1), SM(2), or SM(3) arrow and $\fc_{\rm tot}( \upsilon_{\wp\wp'} )$
the associated cost according to the table of \S \ref{CombCostSec}. 
\end{lem}

\smallskip
\subsection{Tropical semiring and Perron-Frobenius problem} \label{TropPFsec}
We want to extract the dominant term in $t$ of the Perron--Frobenius eigenvalue
and eigenvector of $\cK^{\cP'(n),R,w}(t)$, in the limit where $t\to 0$, as in \cite{MLH}.
To this purpose, we turn the Perron--Frobenius problem for $\cK^{\cP'(n),R,w}(t)$
into an associated Perron--Frobenius problem in the min-plus semiring
$S=(\R\sqcup \{ \infty \}, \min, +, \infty, 0)$ that evaluates orders of magnitudes
as $t\to 0$.

\begin{lem}\label{PFminplus}
For $\cK^{\cP'(n),R,w}(t)$ as in \eqref{Kteq}, the Perron-Frobenius problems
\begin{equation}\label{PFKtL}
\sum_{\wp'} \cK^{\cP'(n),R,w}_{\wp,\wp'}(t) \, \eta_{\wp'}(t) = \lambda_t \, \eta_\wp(t) 
\end{equation}
\begin{equation}\label{PFKtR}
\sum_{\wp} \psi_\wp(t) \, \cK^{\cP'(n),R,w}_{\wp,\wp'}(t) = \lambda_t \, \psi_{\wp'}(t) 
\end{equation}
with leading orders $\lambda_t=\lambda \, t^\ell + o(t^\ell)$, $\eta_\wp(t)= \eta_\wp\, t^{u_\wp} + o( t^{u_\wp} )$
and $\psi_\wp(t) =\psi_\wp\, t^{q_\wp} + o(t^{q_\wp})$ determines Perron-Frobenius problems
for the orders of magnitudes $u_\wp$, $q_\wp$, $\ell$ in the min-plus semiring of the form
\begin{equation}\label{PFKtLmp}
\min_{\wp'} \, \{ \fc_{\wp,\wp'} + u_{\wp'} \} = \ell + u_\wp
\end{equation}
\begin{equation}\label{PFKtRmp}
\min_\wp\,  \{ \fc_{\wp,\wp'} + q_\wp \} = \ell + q_{\wp'} \, . 
\end{equation}
\end{lem}

\proof
We write \eqref{PFKtL} as
$$ \sum_{\wp'}\,\, ( \cK^{\cP'(n),R}_{\wp,\wp'} t^{\fc_{\wp,\wp'} } \eta_{\wp'}\, t^{u_{\wp'}} + o( t^{u_{\wp'}+\fc_{\wp,\wp'} } ) ) =
\lambda\, \eta_\wp\, t^{\ell + u_\wp} + O(t^{\ell + u_\wp}) \, . $$
Due to the positivity of the matrix and the PF data, there are no cancellations. 
The leading term when $t\to 0$ in the sum on the left-hand-side is then equal to 
$\min \{ \fc_{\wp,\wp'} +u_{\wp'} \}$ and this has to match
the leading term $\ell + u_\wp$ on the right-hand-side, giving \eqref{PFKtLmp}.
The case of \eqref{PFKtR} and \eqref{PFKtRmp} is analogous.
\endproof

\smallskip
\subsection{Optimal arborescences and Perron--Frobenius eigenvalue}\label{ArborSec} 
To give an explicit description of the Perron--Frobenius eigenvectors of the min-plus
Perron-Frobenius problem \eqref{PFKtLmp} and \eqref{PFKtRmp}, we first recall the
notion of {\em optimal arborescence} in a weighted directed graph.

\begin{defn}\label{Arbor}
Let $G$ be a directed graph with edges weighted by a cost function $\fc: E(G) \to \R$. 
The weight of a subgraph $G'\subset G$ is given by the total weight of its edges,
$$ \fc(G')=\sum_{e\in E(G')} \fc(e) \, . $$
For a choice of a vertex $v_r\in V(G)$, an {\em arborescence} $T$ with sink at $v_r$ is a subgraph
$T\subset G$ that satisfies 
\begin{enumerate}
\item $T$ is a spanning tree of the underlying undirected graph of $G$,
\item for every $v\in V(G)$ there is a directed path of edges in $E(T)$ from $v$ to the root vertex
$v_r$.
\end{enumerate}
An arborescence $\tilde T$ with source at $v_r$ satisfies condition (1) above and
\begin{itemize}
\item[$(2^\prime)$] for every $v\in V(G)$ there is a directed path of edges in $E(T)$ from the root vertex $v_r$ to $v$.
\end{itemize}
An arborescence $T$ or $\tilde T$ is an {\em optimal arborescence} if it also satisfies the third property:
\begin{enumerate}
\setcounter{enumi}{2}
\item $T$ (or $\tilde T$) has the minimal weight among all arborescences with sink (or with source) at $v_r$. 
\end{enumerate}
\end{defn}

If the graph $G$ is strongly connected, an arborescence with sink at $v_r$ exists for any
choice of the vertex $v_r$. When this is the case, an optimal arborescence can be
constructed through the Chu-Liu/Edmonds/Bock Algorithm, which first constructs a
subgraph $G'$ of $G$ obtained by retaining at each vertex of $G$ the edges of minimal
cost among all the incoming edges, and then proceeding with a recipe for 
eliminating loops in $G'$ to obtain $T$. The case with source at $v_r$ is analogous
but selecting outgoing edges with minimal cost. (For a discussion of algorithms to find
minimal arborescences, see \cite{Gabow}.)

\smallskip

We have the following description of the asymptotic orders of magnitude as $t\to 0$ for the 
Perron--Frobenius eigenvalue and eigenspace, in the case of the Merge dynamics weighted
with the total cost function $\fc_{\rm tot}$.

\begin{thm}\label{totcPFminplusell}
Consider the $\fc_{\rm tot}$-weighted Merge dynamics with matrix $\cK^{\cP'(n),R,w}(t)$ as in \eqref{Kteq},
and with asymptotics given by the solutions to the min-plus Perron-Frobenius problem \eqref{PFKtLmp} and
\eqref{PFKtRmp}. The min-plus Perron--Frobenius eigenvalue is $\ell=0$. 
\end{thm}

\proof
We compute what the order of magnitude $\ell$ for the Perron--Frobenius
eigenvalue should be, so that the PF systems \eqref{PFKtLmp} and \eqref{PFKtRmp}
become equations for the orders of magnitudes $u_\wp$ and $q_\wp$, with $\fc_{\wp,\wp'}$ and
$\ell$ as known quantities. 

There is a general method for computing the Perron--Frobenius eigenvalue over semirings,
which is used for instance in \cite{Akian}, \cite{Akian2}, \cite{Cuni}, \cite{Fried}. For a semiring
$(S,\oplus,\odot,e,\epsilon)$ (in particular the min-plus semiring we are interested in here), 
one defines the {\em permanent} of an $S$-valued $n\times n$ matrix as
$$ {\rm per}(K) = \bigoplus_{\sigma \in S_n} \bigodot_{i=1}^n K_{i\sigma(i)} \, . $$
In the case of the min-plus semiring this gives
$$ {\rm per}(K) = \min_{\sigma\in S_n} (K_{1,\sigma(1)}+ \cdots + K_{n,\sigma(n)}) \, . $$
One can use this to define an analog of the characteristic polynomial
$$ {\rm per}(Y\, {\rm Id} \oplus K) = \bigoplus_{\sigma \in S_n} \bigodot_{i=1}^n (Y \delta_{i,\sigma(i)} \oplus K_{i,\sigma(i)}) 
\in S[Y] \, , $$
whose coefficients are the semiring k-th traces
$$ \Tr^S_k(K) := \bigoplus_{J \subset \{ 1, \ldots, n \}, \, \#J=k} \left( \bigoplus_{\sigma\in S_J} \bigodot_{j\in J} A_{j,\sigma(j)}
\right) $$
with $\sigma \in S_J$ the permutations of the subset $J$. The roots of this characteristic polynomial
are called ``algebraic eigenvalues" of $K$ (see \cite{Akian2}). Note that in general not all of these algebraic 
eigenvalues will have eigenvectors in this semiring setting: those that do are referred to as ``geometric eigenvalues".
For the min-plus semiring, in particular, the minimal algebraic eigenvalue is given by
\begin{equation}\label{rhominK}
 \rho_{\min} (K)= \bigoplus_{k=1}^n \bigoplus_{i_1,\ldots, i_k} (K_{i_1 i_2}\odot \cdots \odot K_{i_k i_1})^{1/k} 
= \min_{k=1, \ldots, n} \min_{c=(i_1,\ldots,i_k)} \frac{1}{k} (K_{i_1 i_2} + \cdots + K_{i_k i_1} )\, . 
\end{equation}
where $c=(i_1,\ldots,i_k)$ runs over all circuits of length $k$ in the graph $G(K)$ with $n$ vertices
and an edge $e_{ij}$ for each nontrivial entry $K_{ij}$.  A circuit $c=(i_1,\ldots,i_k)$ that realizes the
minimum is a ``critical circuit".  The critical graph $G_{\rm crit}(K)\subset G(K)$ consists
of vertices and edges that belong to critical circuits. If $G(K)$ is strongly connected, then $\rho_{\min} (K)$ is the unique
geometric eigenvalue and its multiplicity is equal to the number of strongly connected components
of $G_{\rm crit}(K)$. 

\smallskip.

The graph $\cG_{\cP'(n)}$ with the edges $e_{\wp,\wp'}$ 
weighted by the entries $\cK^{\cP'(n),R}_{\wp,\wp'} t^{\fc_{\wp,\wp'} }$ of the matrix 
$\cK^{\cP'(n),R}(t)$ is strongly connected and 
has IM self-loops at all the partitions $\wp$ that contain at least
one row of size $k_i\geq 3$. Since the cost $\fc_{\rm tot}({\rm IM})=0$, these
IM diagonal entries of $\cK^{\cP'(n),R}(t)$ are independent of $t$ and equal to 
the corresponding entry $\cK^{\cP'(n),R}_{\wp,\wp}$. 

\smallskip

Let $K$ be the min-plus valued matrix with entries the costs
$K_{\wp,\wp'}=\fc_{\wp,\wp'}$ of the Merge transformations
associated to the edges, as described in \S \ref{CombCostSec}. 
Let $G(K)$ be the graph associated to the min-plus matrix $K$
as above. The graph $G(K)$ has the same vertices and edges
of $\cG_{\cP'(n)}$, with the edges weighted by $\fc_{\wp,\wp'}$,
hence it is also strongly connected, since none of the 
$\fc_{\wp,\wp'}$ equals $\infty$ (the additive unit of the min-plus semiring).
The IM self-loops of $\cG_{\cP'(n)}$ have weight $\fc_{\wp,\wp'}=0$ (the
multiplicative unit of the max-plus semiring). Since all other
costs coming from non-IM (EM or SM) arrows are larger, $\fc_{\wp,\wp'}>0$,
the circuits $c=(\wp_1,\ldots, \wp_k)$ that realize the minima of
$(\fc_{\wp_1,\wp_2}+\cdots+\fc_{\wp_k,\wp_1})/k$ are the single-edge IM loops
$\fc_{\wp,\wp}=0$ (with $k=1$). Thus, there are many critical circuits, each consisting of
a single edge, and a single vertex (any partition $\wp$ with at least one
row of length at least three is such a vertex). Thus, the critical graph $G_{\rm crit}(K)$ consists
of a disjoint union of single-edge loops, one for each such partition $\wp$. 

\smallskip

Thus, we have $\ell=0$ in \eqref{PFKtLmp} and \eqref{PFKtRmp}, as this $\ell$
corresponds to the unique geometric eigenvalue \eqref{rhominK}. 
\endproof

\smallskip
\subsection{Tropical Perron-Frobenius eigenspace}\label{TropEigSpaceSec}

We then analyze the Perron-Frobenius eigenspace of the tropical
Perron-Frobenius problem  \eqref{PFKtLmp} and \eqref{PFKtRmp} with $\ell=0$.

\smallskip

Since there
are multiple critical circuits, we cannot use the eigenvalue asymptotics of \cite{Akian}, \cite{Akian2}
that relies on the uniqueness of the critical circuits to relate the min-plus (or max-$\times$) Perron-Frobenius
eigenvector to the Perron-Frobenius eigenvector of the original problem \eqref{PFKtL}, \eqref{PFKtR}. 
However, we can describe the Perron-Frobenius eigenspace in terms of its generators, as in \cite{Akian3}. 

\smallskip

\begin{thm}\label{totcPFminplus}
Consider the cost function $\fc_{\rm tot}$ of the dynamics $\cK^{\cP'(n),R,w}(t)$ as in \eqref{Kteq},
with asymptotics determined by the min-plus Perron-Frobenius problem \eqref{PFKtLmp} and
\eqref{PFKtRmp}. The right Perron-Frobenius
eigenspace is generated by a family $\{ u^{(\wp)} \}$ of min-plus eigenvectors as in \eqref{PFKtLmp},
where $\wp$ ranges over those partitions $\wp\in \cP'(n)$ whose Young diagram has at least one row 
of length at least $3$. The entry $u^{(\wp)}_{\wp'}$ of the basis eigenvector 
$u^{(\wp)}=(u^{(\wp)}_{\wp'})_{\wp'\in \cP'(n)}$ is given by the cost $\fc(\gamma_{\wp'\wp})$ of the
directed path in an optimal arborescence $T_\wp$ in $\cG_{\cP'(n)}$ with sink at $\wp$. Similarly, 
the left Perron-Frobenius eigenspace is generated by a family $\{ q^{(\wp)} \}$ of solutions of \eqref{PFKtRmp},
for $\wp$ in the same set of partitions $\wp$ with at least one row of length at least $3$, where for a
for $\wp'\in \cP'(n)$, the entry $q^{(\wp)}_{\wp'}$ is the cost $\fc(\tilde\gamma_{\wp\wp'})$ of the directed
path $\tilde\gamma_{\wp\wp'}$ in an optimal arborescence $\tilde T_\wp$ in $\cG_{\cP'(n)}$ with source at $\wp$. 
\end{thm}

\proof
One defines the Kleene star $K^\star$ of a min-plus valued matrix $K$ in the form
\begin{equation}\label{Kleene}
K^\star := \cI \oplus K \oplus K^2 \oplus \cdots \oplus K^{n-1}\, , 
\end{equation}
where $\cI$ is the identity matrix (in the min-plus sense). The entries are given by
$$ (K^\star)_{ij} = K_{i=i_0,i_1}\odot \cdots \odot K_{i_{k-1}, i_k=j} $$
along the directed path of minimal weight from $i$ to $j$ in $G(K)$, which is 
of some length $1\leq k \leq n-1$.  As shown in Proposition~2.6 of \cite{Akian3}
(specialized to our case with $\rho_{\min} (K)=0$ and
with each critical cycle $C_j$ consisting of a single vertex $j$ and a looping edge), 
there is an eigenvector $u_j$ of $K$ for each critical cycle $C_j$, proportional to 
the $j$-th column of the matrix $K^\star$. These eigenvectors are the generators of
the eigenspace of $\rho_{\min} (K)$. 

\smallskip

In our case, this means that, for each partition $\wp$ with at least one row of length $k_i \geq 3$
there is an eigenvector $u^{(\wp)}=(u^{(\wp)}_{\wp'})_{\wp'\in \cP'(n)}$ 
that satisfies \eqref{PFKtLmp} with $\ell=0$. Correspondingly,
the same argument applied to the transpose matrix shows that for each such $\wp$ there is an
eigenvector $q^{(\wp)}=(q^{(\wp)}_{\wp'})_{\wp'\in \cP'(n)}$ that satisfies \eqref{PFKtRmp} with $\ell=0$. 
The set $\{ u^{(\wp)} \}$ (respectively, $\{ q^{(\wp)} \}$) of these eigenvectors spans the eigenspace.

\smallskip

Moreover, the description of these generators in terms of columns of the Kleene star matrix
can be made geometrically more explicit in our case, in terms of the graph $\cG_{\cP'(n)}$. 
To see this, we use the notion of optimal arborescence in a weighted directed graph, recalled
in Definition~\ref{Arbor} above.

\smallskip

Let $\wp$ be one of the
partitions in $\cP'(n)$ that have at least one row with $k_i\geq 3$ so that they are a non-trivial
IM self-loop edge in $\cG_{\cP'(n)}$. And let $T_\wp$ denote an optimal arborescence with sink 
at $\wp$, with respect to the cost function $\fc_{\rm tot}$. This is a directed spanning tree of
$\cG_{\cP'(n)}$, where each vertex $\wp' \neq \wp$ has a single outgoing edge 
in $T_\wp$, with minimal cost among all the outgoing edges at $\wp'$ in $\cG_{\cP'(n)}$,  
while the vertex $\wp$ has no outgoing edges in $T_\wp$. 

\smallskip

In the construction of the Kleene star matrix, we have
$$ (K^\star)_{\wp', \wp}= \fc(e_{\wp',\wp_1})+ \cdots + \fc(e_{\wp_{k-1},\wp}) $$
where $e_{\wp',\wp_1} \ldots e_{\wp_{k-1},\wp}$ is a minimal cost path from
$\wp'$ to $\wp$. In particular, this means (Bellman's optimality principle) that
for any $0\leq i < j \leq k$ the subpath $e_{\wp_i, \wp_{i+1}} \ldots e_{\wp_{j-1}, \wp_j}$
is also a cost minimizing path between $\wp_i$ and $\wp_j$, so in particular at
each vertex $\wp_i$ the next outgoing edge in the path can be taken to be an edge in $T_\wp$
and we can identify $(K^\star)_{\wp', \wp}$ with the $\fc$-cost of the unique path 
$\gamma_{\wp',\wp}$ in $T_\wp$ from $\wp'$ to $\wp$. 
Thus, the eigenvectors $u^{(\wp)}$ satisfying \eqref{PFKtLmp} 
that span the eigenspace are associated to the optimal arborescences $T_\wp$
and have components $u^{(\wp)}=(u^{(\wp)}_{\wp'})_{\wp'\in \cP'(n)}$ given
by the costs $u^{(\wp)}_{\wp'}=\fc(\gamma_{\wp',\wp})$ of the paths in the
arborescence. 

\smallskip

In the case of the eigenvalues $q^{(\wp)}$ of the \eqref{PFKtRmp} min-plus Perron-Frobenius problem,
we are considering the same setting but with the transpose $K^\tau$ of the $K$ matrix. Thus, the vector
$q^{(\wp)}$ is the $\wp$-column of the Kleene star matrix $(K^\tau)^\star$ of the transpose $K^\tau$. 
The entries of this Kleene star matrix are given by
$$  (K^\tau)^\star_{\wp', \wp}= \fc(e_{\wp,\wp_1})+ \cdots + \fc(e_{\wp_{k-1},\wp'}) = \fc(\tilde\gamma_{\wp,\wp'})\, , $$
where $\tilde\gamma_{\wp,\wp'}=e_{\wp,\wp_1}\ldots e_{\wp_{k-1},\wp'}$ is a minimal cost path from $\wp$ to $\wp'$. 
Thus, in this case the eigenvectors $q^{(\wp)}$ correspond to optimal arborescences $\tilde T_\wp$ with source at $\wp$,
with the components $q^{(\wp)}_{\wp'}=\fc(\tilde\gamma_{\wp',\wp})$ given by the costs of the paths
in the arborescence. 
\endproof

\smallskip

Moreover, the result of \cite{Akian3} shows that the left/right Perron-Frobenius eigenvector of
$\cK^{\cP'(n),R,w}(t)$ has asymptotics $u=(u_\wp)_{\wp\in \cP'(n)}$ (respectively, $q=(q_\wp)_{\wp\in \cP'(n)}$)
in the eigenspace, although not necessarily equal to one of the above generators of the eigenspace.
In fact, explicit examples computed in \cite{Kustar} show that it is in general not the case that the
asymptotics of the Perron-Frobenius eigenvector is one of the generators of the eigenspace of the
min-plus Perron-Frobenius problem.

\smallskip
\subsection{Optimization with respect to cost functions}\label{FreeEnCostSec}

The introduction of optimality with respect to Minimal Search, Minimal Yield, and Complexity Loss
in the linguistics model of Minimalism, \cite{ChomskyUCLA}, \cite{ChomskyGK}, \cite{MCB}, formulated
in the setting we are considering here, of the Merge Hopf algebra Markov chain dynamical system, 
becomes the minimization property of the Markov chain 
$$ \hat\cK^{(A,n,(\fc_{\rm tot},\beta))}  = \hat\cK^{(A,n,\fc_{\rm tot})}(t) \ \ \ \text{ for } t=e^{-\beta} \, , $$ 
with $\fc_{\rm tot}$ the cost function of \S \ref{CombCostSec}, with respect to the free energy function  
$$ \bF(\hat\cK^{(A,n,(\fc_{\rm tot},\beta))})=\bar\bE_{\bP_{(\fc_{\rm tot},\beta)}} -\beta^{-1} {\rm Sh}(\bP_{(\fc_{\rm tot},\beta)})\, $$
where $\bP_{(\fc_{\rm tot},\beta)}$ is the Boltzmann distribution with $\bE= \fc_{\rm tot}$, as 
explained in \S \ref{FreeEnSec}. 

\smallskip

We discuss in \S \ref{SMproblemSec2} to what extent the use of the cost function $\fc_{\rm tot}$
achieves the goal of eliminating the unwanted properties of the Sideward Merge arrows, allowing
for convergence via External Merge to fully formed sentences (tree structures). We will then present
in \S \ref{EntCostSec} an alternative choice of a simple cost function which corrects for the remaining
problems and formulates optimization entirely in information theoretic terms.

\smallskip
\subsection{The Sideward Merge problem in the weighted dynamics}\label{SMproblemSec2}

The result of our discussion of the weighted dynamics with cost function $\fc_{\rm tot}$ and
weights $t^{\fc_{\rm tot}}$ on the adjacency matrix of the Merge graph $\cG_{n,A}$ shows
that the cost function corrects for some of the unwanted properties of the SM arrows in the unweighted dynamics.
However, the surprise with respect to the simple case $n=3$ discussed in \cite{MLH} is that
the cost counting $\fc_{\rm tot}$ alone is not sufficient to completely correct for the SM pull 
on the dynamics in the direction of less connected structures.

\smallskip

More precisely, we found that, to leading order when $t\to 0$, the stationary distribution for the
weighted dynamics lies in the Perron--Frobenius eigenspace described above. Since the
basis of this eigenspace only involves partitions $\wp$ with at least one of the $k_i\geq 3$,
we no longer see the phenomenon typical of the unweighted dynamics, where the stationary
distribution assigns highest probability to the most disconnected structure $\wp_{1^n}$
due to the dominance of SM. However, the leading order of the stationary distribution can 
now have contributions from all the $\wp$ with at least one of the $k_i\geq 3$ (all those
with IM dynamics), and is not supported only on $\wp=\{ n \}$ (convergence to the 
connected structures, workspaces consisting of a single tree, with the IM dynamics), as
one expects when EM dominates the SM--EM tension in the dynamics, and as happens in the
case $n=3$ of \cite{MLH} with the same weights and cost function. 

\smallskip

This reveals that Minimal Search, Minimal Yield, and Complexity Loss 
as usually considered in the linguistics literature (the slightly different MY formulation
of \cite{ChomskyUCLA}, \cite{ChomskyGK} makes no difference) are not sufficient
cost functions to eliminate the unwanted effects of Sideward Merge.

\smallskip

As we discuss in the next section, what is missing is the role of {\em information
optimization}, which as we have already observed in a different form in \S  \ref{RWvsHMCsec}, and especially \S \ref{MaxEntSec}, 
that the Merge Hopf algebra Markov chain satisfies. We show here that adding a 
Shannon information count to the cost function achieves the result of 
concentrating the dynamics, in the leading order part of the stationary
distribution, on just the workspaces consisting of a single tree with the IM dynamics.

\section{Shannon entropy optimization and the Merge dynamics} \label{EntCostSec}

\begin{defn}\label{PartProb}
To a partition $\wp\in \cP'(n)$, of the form $n=a_1 k_1+\cdots + a_r k_r$, we associate the 
probability distribution 
\begin{equation}\label{Ppart}
\P_\wp := (\underbrace{\frac{k_1}{n}, \ldots, \frac{k_1}{n}}_{a_1\text{-times}}, \ldots, \underbrace{\frac{k_r}{n}, \ldots, \frac{k_r}{n}}_{a_r\text{-times}} ) \, .       
\end{equation}
\end{defn}

We extend the cost function $\fc_{\rm tot}$ considered in the previous section, by including a term
that depends on the Shannon entropy ${\rm Sh}(\P_\wp)$  of the distribution \eqref{Ppart}. 

\begin{defn}\label{costSH}
On the graph $\cG_{\cP'(n)}$ we consider a cost function $\fc: E(\cG_{\cP'(n)})\to \R$ of the form
\begin{equation}\label{totcostSh}
\tilde\fc_{\rm tot}(e) = \fc_{\rm tot}(e) + {\rm Sh}(\P_{s(e)}) \, , 
\end{equation}
where $\fc_{\rm tot}$ is the total cost function of MS, MY, CL, as described in \S \ref{CombCostSec} and
$\P_{s(e)}$ is the distribution \eqref{Ppart} at the partition $\wp=s(e) \in \cP'(n)$, the source vertex of the edge $e$
in the directed graph $\cG_{\cP'(n)}$. 
\end{defn} 

The same argument used in Theorem~\ref{totcPFminplusell} and
Theorem~\ref{totcPFminplus}, applied to the modified cost function \eqref{totcostSh}, 
leads to the following result. 

\begin{prop}\label{costShPF}
With the modified cost function $\tilde\fc_{\rm tot}$ of \eqref{totcostSh}, the dynamics 
$\cK^{\cP'(n),R,w}(t)$ as in \eqref{Kteq} has asymptotics given by the solution of the min-plus 
Perron-Frobenius problem \eqref{PFKtLmp} and \eqref{PFKtRmp} 
with $\ell=0$ and with PF L/R-eigenvectors $\upsilon=\upsilon^{(\wp=\{ n \})}$,
wuth $\upsilon_{\wp'}=\tilde\fc(\gamma_{\wp' \{ n \}})$, the cost of the directed path in an optimal arborescence
$T_{\{ n \} }$ with sink at $\wp=\{ n \}$, and $q=q^{(\wp=\{ n \})}$, with $q_{\wp'}=\tilde\fc(\tilde\gamma_{\{ n \} \wp'})$,
the cost of the directed path in an optimal arborescence  $\tilde T_{\{ n \} }$ with source at $\wp=\{ n \}$.
\end{prop} 

\proof The main difference with respect to Theorem~\ref{totcPFminplusell} and
Theorem~\ref{totcPFminplus} is that, with the modified cost function \eqref{totcostSh}, 
there is only one critical circuit, consisting of the IM self-loop at the partition $\wp=\{ n \}$.
Indeed, the added term ${\rm Sh}(\P_{s(e)})$ in the cost function is minimal, 
${\rm Sh}(\P_\wp) = 0$ iff $\wp=\{ n \}$ and is strictly positive at all other $\wp\in \cP'(n)$. 
Thus, the only minimal cost loop, with $\tilde\fc_{\rm tot}(\gamma) =0$ is the single IM edge
with source and target $\wp=\{ n \}$. This again gives $\ell=0$ for the min-plus PF eigenvalue
as in Theorem~\ref{totcPFminplusell}. Since in this case we do have a unique  critical circuit, 
we can obtain not only a basis of the PF eigenspace of the min-plus Perron-Frobenius problem,
but also directly the PF eigenvector which gives the asymptotic behavior of the original
Perron-Frobenius problem. In particular, this gives the form $\upsilon_{\wp'}=\tilde\fc(\gamma_{\wp' \{ n \}})$
and $q_{\wp'}=\tilde\fc(\tilde\gamma_{\{ n \} \wp'})$ or the R/L-eigenvector of the min-plus problem.
\endproof

Consequently we obtain the asymptotic behavior of the stationary distribution of the
weighted Merge dynamics on $\cG_{n,A}$.

\begin{thm}\label{pidistrSh}
The asymptotic behavior for $t\to 0$ of the stationary distribution of the weighted form
$\hat \cK^{A,n,w}(t)$ of the Hopf algebra Markov chain of Merge acting on workspaces,
with weights $t^{\tilde\fc}$ and cost function $\tilde\fc$ of \eqref{totcostSh} is given by
the total cost of the loop along the optimal arborescences $T_{\{ n \} }$ and  $\tilde T_{\{ n \} }$,
\begin{equation}\label{asympiFt}
{\rm order}_{t\to 0}(\pi^w(F,t))= \tilde\fc( \gamma_{p(F)\, \{ n \}} \cup \tilde\gamma_{\{ n \}\, p(F)}) \, .
\end{equation}
This implies that in the limit $t\to 0$ the stationary distribution of $\hat \cK^{A,n,w}(t)$ is
the uniform distribution supported on the workspaces that consist of a single tree, $F=T \in p^{-1}(\{ n \})$. 
Thus, the Merge dynamics weighted by the cost function $\tilde\fc$ converges to fully formed
trees with the remaining action of Internal Merge.
\end{thm}

\proof We know by Proposition~\ref{MergeDynStatDistrW} that the Hopf algebra Markov chain
$\hat \cK^{A,n,w}(t)$ for the dynamics of Merge on $\cG_{n,A}$, with the weighted adjacency
matrix 
$$ \cK^{A,n,w}_{F,F'}(t)= t^{\tilde\fc_{{\rm tot},F,F'}}\,\, \cK^{A,n}_{F,F'} \, , $$
has stationary distribution as in \eqref{piKt}, given by 
$$  \pi^w(F,t) = \frac{\eta_{\cP'(n),w}(\wp,t) \, \, \psi_{\cP'(n),w}(\wp,t)}{\Lambda_{\wp,n}} \, , $$
with $\wp=p(F)$ under the projection map $p:\cG_{n,A} \to \cG_{\cP'(n)}$, where
$\eta_{\cP'(n),w}$ and $\psi_{\cP'(n),w}$ are the left and right Perron--Frobenius eigenfunctions of
$\cK^{\cP'(n),R,w}(t)$ on $\cG_{\cP'(n)}$. Thus, when looking at the asymptotic behavior for $t\to 0$
of $\pi^w(F,t)$, we see from Proposition~\ref{costShPF} that this is given by
$$ {\rm order}_{t\to 0}(\pi^w(F,t))= \upsilon_{p(F)} + q_{p(F)} $$
with $\upsilon=\upsilon^{(\wp=\{ n \})}$ and $q=q^{(\wp=\{ n \})}$ as in Proposition~\ref{costShPF}.
This gives \eqref{asympiFt} as $\upsilon_{p(F)} + q_{p(F)}$ is the $\tilde\fc$-cost of the loop
$\gamma_{p(F)\, \{ n \}} \cup \tilde\gamma_{\{ n \}\, p(F)}$.
All the costs of these loops are strictly positive, except for the case where $p(F)=\{ n \}$ where
the cost is zero. Thus, in the limit $t\to 0$ the $\pi^w(F,t)$ will converge to the uniform
distribution on $p^{-1}(\{ n \})$, the workspaces with a single connected component. 
\endproof

\smallskip

In this form, the weighted dynamics has exactly the expected effect, namely it compensates
for the fact that the presence of Sideward Merge tends to push the dynamics to less
connected structures preventing convergence of the structure formation process to a single tree.
At the same time, one maintains the good properties of Sideward Merge, the strong connectedness
of the graph (irreducibility of the dynamics). 

\smallskip
\subsection{Shannon entropy cost function}

In fact, while we maintained in Proposition~\ref{costShPF} and Theorem~\ref{pidistrSh}
a combination, as in Definition~\ref{costSH} of the cost function $\fc_{\rm tot}$ coming
from the usual notions of Minimal Search, Minimal Yield, and Resource Restriction 
developed in the linguistics literature (see \cite{ChomskyUCLA}, \cite{ChomskyGK}), 
in our weighted dynamics 
it is really the Shannon entropy part of the cost function that has the desired role
of forcing the dynamics to converge to the workspaces consisting of a single tree as
one wants in the structure formation process of syntax. The following result is proved
exactly as Theorem~\ref{pidistrSh}. 

\begin{thm}\label{aloneSh}
Consider the Merge dynamics $\hat \cK^{A,n,{\rm Sh}}(t)$
weighted only by the Shannon entropy cost, namely using the 
weighted adjacency matrix on $\cG_{n,A}$ of the form
\begin{equation}\label{costonlySh}
 \cK^{A,n,{\rm Sh}}_{F,F'}(t): = t^{{\rm Sh}(\P_{p(F)})}\,\, \cK^{A,n}_{F,F'} \, .
\end{equation}
Then the stationary distribution has asymptotic behavior
\begin{equation}\label{piShonly}
{\rm order}_{t\to 0}(\pi^{\rm Sh}(F,t))= \sum_{e \in \gamma_{p(F)\, \{ n \}}} {\rm Sh}(\P_{s(e)}) \,\, + \sum_{e'\in \tilde\gamma_{\{ n \}\, p(F)}} {\rm Sh}(\P_{s(e')}) \, . 
\end{equation}
In particular, ${\rm order}_{t\to 0}(\pi^{\rm Sh}(F,t))>0$ for all $F \in p^{-1}(\wp)$ with $\wp\neq \{ n \}$, while
${\rm order}_{t\to 0}(\pi^{\rm Sh}(F,t))=0$ for $F=T\in p^{-1}(\{ n \})$. Thus, in the limit $t\to 0$ the stationary
distribution converges to the uniform distribution on the trees $T\in p^{-1}(\{ n \})$, where only Internal Merge is acting.
\end{thm}

\begin{rem}\label{remShalone} {\rm The result of Theorem~\ref{aloneSh}
shows that the Minimal Search and Resource Restrictions cost functions
play a role in simplifying the graph $\cG_{n,A}$ by identifying the
``minimal" Sideward Merge arrows that suffice for strong connectedness
while significantly reducing the size of the graph, making it sparse.  Optimization
of information-theoretic data then drives the dynamics towards the connected workspaces, 
compensating for the presence of Sideward Merge. 
}\end{rem}

\smallskip

This observation leads to the natural question on the meaning of the added 
cost function $\fc_{\rm Sh}(e)={\rm Sh}(\P_{s(e)})$ in 
the context of the geometry of the projection map $p: \cG_{n,A} \to \cG_{\cP'(n)}$. 

\smallskip

We have seen in \S \ref{MaxEntSec} how the Merge Hopf algebra Markov chain $\hat\cK^{A,n}$
has an entropy maximization property for the Shannon entropy of the probability on paths determined
by the probability on edges defined by the transition probabilities of a Markov chain. We have also
seen in \S \ref{FreeEnSec} how, introducing a weight by a cost function on edges, one can reformulate this optimization
as a minimization of the free energy functional $\bar\bF = \bar \bE - \beta^{-1} {\rm Sh}$, with
the energy $\bE$ determined by the cost function and the probability distribution on edges determined
by the transition probabilities of the (weighted) Markov chain. 

\smallskip

We are considering here an additional datum of a probability distribution $\P_\wp$ associated to 
each vertex $\wp\in \cP'(n)$. The meaning of the quantity ${\rm Sh}(\P_\wp)$ can be described as follows.

\begin{rem}\label{ShPmean}{\rm
By definition, the Shannon entropy ${\rm Sh}(\P_\wp)$ measures the expected amount of information
(or the expected surprisal) in observing the outcome of an event distributed according to $\P_\wp$.
If we fix a workspace $F=T_1 \sqcup \cdots \sqcup T_r$ in $p^{-1}(\wp) \subset V(\cG_{n,A})$, we
can observe how the lexical material (the elements of the set $A$) at the leaves is distributed among
the components $T_i$ of the workspace $F$. Since we are here dealing only with {\em free} structure
formation (the free symmetric Merge), there are no constraints on how the elements of $A$ can be 
distributed among the leaves of the $T_i$, and this implies that the only law governing the possible
assignments is the probability $\P_\wp$. (We will see in \S \ref{ColorSec} and \S \ref{ExtSec} that
filtering for theta roles assignments, phase structure, and syntactic parameters will impose constraints
that change this probability distribution.)}
\end{rem}

\begin{rem}\label{multiSh}{\rm 
The probability distribution $\P_\wp$ and the Shannon entropy ${\rm Sh}(\P_\wp)$ are
directly related to the combinatorial counting of the assignments of the elements of $A$ to
the components of the workspace, via the well known relation of ${\rm Sh}(\P_\wp)$  to the
multinomial coefficient $\mu_{\wp,n}$ of \eqref{multinomial}, \eqref{multinomial2}, namely
the asymptotic behavior
\begin{equation}\label{muShP}
 \mu_{n,\wp} \sim e^{n \, {\rm Sh}(\P_\wp) + o(n)}\, .  
\end{equation} 
}\end{rem}

This can be derived directly from the Stirling approximation of factorials, which
also shows that the lower order term behaves like $o(n)=O(\log(n))$.
The behavior in \eqref{muShP} can be interpreted in terms of 
the fact that $\mu_{n,\wp}$ counts words of length $n$ in an alphabet of $r$ symbols, where
the $i$-th symbol occurs $k_i$ times. 

\smallskip

\begin{rem}\label{multiCount}{\rm
In the counting of the multiplicities of the projection map $p: \cG_{n,A} \to \cG_{\cP'(n)}$,
namely the sizes of the fibers $\# p^{-1}(\wp)=\Lambda_{\wp,n}$ we have $\Lambda_{\wp,n}=
\mu_{n,\wp}\cdot \Gamma_{\wp,n}$ with $\Gamma_{\wp,n}$ as in \eqref{multinomialGamma}, where
the factor $\Gamma_{\wp,n}$ accounts for permuting rows of equal length in the partition $\wp$
and for counting all the different possible tree topologies of a tree $T_i$ with $k_i$ leaves. 
However, unlike in the counting of Proposition~\ref{dynIMunif1}, when we have fixed a 
workspace $F=T_1 \sqcup \cdots \sqcup T_r$ in $p^{-1}(\wp)$,
and we look at the distribution of the lexical material $A$ among the $T_i$, only the $\mu_{n,\wp}$
part of the multiplicity $\Lambda_{\wp,n}$ remains relevant, as the other choices encoded in the
$\Gamma_{\wp,n}$ factor have been fixed. }
\end{rem}

\smallskip

The effect on the partition of various Merge transformations has a
corresponding effect on the probability distribution $\P_\wp$ and
its Shannon entropy that can be phrased in terms of natural properties
of the entropy functional. For example, we have the following observation.

\begin{prop}\label{EMchainruleSh}
External Merge, acting on partitions by $$\wp=\{ k_1, \ldots, k_r \} \mapsto \wp'=\{ k_1, \ldots, k_i+k_j, \ldots, k_r \}$$
acts on ${\rm Sh}(\bP_\wp)$ as the chain rule of the Shannon entropy as
\begin{equation}\label{EMgroupSh}
{\rm EM}: {\rm Sh}(\P_\wp) \mapsto {\rm Sh}(\P_{\wp'}) =  {\rm Sh}(\P_\wp) 
- \frac{k_i+k_j}{n} {\rm Sh}( \frac{k_i}{k_i+k_j}, \frac{k_j}{k_i+k_j})
\end{equation}
\end{prop}

\proof By the symmetry property of the Shannon entropy, we can assume without loss of 
generality that the pair $k_i, k_j$ is $k_1,k_2$. 
The chain rule for the Shannon entropy (also known as extensivity condition) gives
\begin{equation}\label{groupingSh}
{\rm Sh}(\frac{k_1}{n}, \ldots, \frac{k_r}{n}) = {\rm Sh}(\frac{k_1+k_2}{n}, \frac{k_3}{n}, \ldots, \frac{k_r}{n}) +
\frac{k_1+k_2}{n} {\rm Sh}( \frac{k_1}{k_1+k_2}, \frac{k_2}{k_1+k_2}) \, ,
\end{equation}
which translates into the action of an EM transformation that merges the components $T_i$ and $T_j$ 
of the workspace to $\fM(T_i,T_j)$ on the Shannon entropy by \eqref{EMgroupSh}.
\endproof

\begin{rem}\label{muShchain}{\rm
The identity of multinomial coefficients 
$$ \binom{n}{ m , k_3, \ldots, k_r } \cdot \binom{m}{k_1, k_2} = \binom{n}{k_1, k_2, k_3, \ldots, k_r } $$
implies, in the asymptotic behavior the chain rule identity \eqref{groupingSh} for the Shannon entropy.}
\end{rem}

\begin{rem}\label{EMspecialSh}{\rm
Note how this singles out External Merge as being associated to the
chain rule of the Shannon entropy, which is the key property in the 
axiomatic characterization (Khinchin axioms) of information.  
It seems interesting to investigate whether Internal Merge and Sideward Merge also
have a meaning in terms of information theoretic properties. 
}
\end{rem}

\smallskip

When we use ${\rm Sh}(\P_\wp)$ as an energy/cost function, by setting
$\fc_{\rm Sh}(e)={\rm Sh}(\P_{p(s(e))})$, we weight the adjacency matrix
$\cK^{A,n}$ of $\cG_{n,A}$ as 
$$ \cK^{A,n, {\rm Sh}}_{F,F'} = t^{{\rm Sh}_{p(F)}}\, \, \cK^{A,n}_{F,F'} = e^{-\beta{\rm Sh}_{p(F)}}\, \, \cK^{A,n}_{F,F'} \, , $$
with $t=e^{-\beta}$, with $\beta$ the inverse temperature thermodynamic parameter of the Boltzmann distribution, so
that the range $0< t <1$ corresponds to the usual range $0 < \beta < \infty$ for the inverse temperature, with $t\to 0$
corresponding to the zero-temperature limit $\beta \to \infty$. 
Using this notation, we write here $\cK^{A,n, {\rm Sh},\beta}_{F,F'}$ rather than $\cK^{A,n, {\rm Sh}}_{F,F'}(t)$ as
before, for consistency with \S \ref{FreeEnSec}. 

\smallskip

\begin{prop}\label{maxShunifIM}
In the thermodynamic limit of zero temperature, $\beta \to \infty$ the free energy optimization for 
$\hat\cK^{A,n, {\rm Sh}}_{F,F'}$ is achieved by the uniform distribution on edges of the Internal Merge
graph $\cG^{\rm IM}_{n,A,\{ n \}}$ that maximizes the Shannon entropy.
\end{prop}

\proof
As shown in \S \ref{FreeEnSec}, the corresponding Markov chain $\hat\cK^{A,n, {\rm Sh},\beta}$ minimizes
the free energy function
$$ \bF(\hat\cK^{A,n, {\rm Sh},\beta}) =  \bar\bE(\hat\cK^{A,n, {\rm Sh},\beta}) 
-\beta^{-1} {\rm Sh}(\hat\cK^{A,n, {\rm Sh},\beta}) \, . $$
In order to look at the $\beta \to \infty$ behavior, it is more convenient to consider the
maximization of the form of the energy function as in \eqref{freeendefMax},
$$  \tilde\bF(\cS)={\rm Sh}(\cS)  - \beta \overline{\bE}(\cS) \, . $$ 

As shown in Lemma~\ref{RWfreeEn}, this free energy is given as in \eqref{BolzEbar} by 
$$ \tilde\bF(\hat\cK^{A,n, {\rm Sh},\beta}) = -\beta \overline{\bE}_{\bP_{({\rm Sh},\beta)}} + {\rm Sh}(\bP_{({\rm Sh},\beta)})\, , $$
where 
$$ \bP_{({\rm Sh},\beta)}(e)= \frac{1}{\lambda} \psi(s(e)) e^{-\beta {\rm Sh}(\P_{p(s(e))})} \eta(t(e)) $$
with $\lambda$ and $\psi$, $\eta$ the PF eigenvalue and L/R PF eigenvector of $\cK^{A,n, {\rm Sh},\beta}$. 
Thus, we have
$$ \bar\bE(\hat\cK^{A,n, {\rm Sh},\beta}) =  \sum_F \pi(F) {\rm Sh}(\P_{p(F)}) $$
with $\pi$ the stationary distribution of $\hat\cK^{A,n, {\rm Sh},\beta}$. 
As discussed in  \S \ref{FreeEnSec}, we have $\bF(\hat\cK^{A,n, {\rm Sh},\beta})=-\beta^{-1}(\log(\lambda) + {\rm Sh}(\pi))$
as in \eqref{BoltzOpt3}, 
for $\lambda=\lambda_{({\rm Sh},\beta)}$ the PF eigenvalue of $\cK^{A,n, {\rm Sh},\beta}$. 

Theorem~\ref{aloneSh} further shows that, in the zero-temperature limit $\beta \to \infty$ ($t\to 0$),
the stationary distribution $\pi(F)=\pi^{({\rm Sh},\beta)}(F)$ converges to the uniform distribution on
the fiber $p^{-1}(\{ n \})$. Thus, the terms in the sum
$$ \sum_F \pi^{({\rm Sh},\beta)}(F)\,\, {\rm Sh}(\P_{p(F)}) $$
have $\pi^{({\rm Sh},\beta)}(F) \to 0$ for $\beta \to \infty$, for all $F\notin p^{-1}(\{ n \})$, where
${\rm Sh}(\P_{p(F)})\neq 0$ and $\pi^{({\rm Sh},\beta)}(F) \to 1$ for $F\in p^{-1}(\{ n \})$, where
${\rm Sh}(\P_{p(F)})=0$ so that the sum vanishes in the limit. 

Moreover, the distribution $\P_{({\rm Sh},\beta)}(e)$ on $\cG_{n,A}$ reduces to the distribution
$\bP_{\rm IM}(e)=\hat \cK^{\rm IM}_{T,T'}$ on the set of edges
of the Internal-Merge graph $\cG^{\rm IM}_{n,A,\{ n \}}$, with fixed $T=s(e)$ and with $T'=t(e)$.
By Proposition~\ref{dynIMunif} this distribution  is the uniform distribution
$$ \bP_{\rm IM}(e) =\frac{1}{\lambda_{\rm IM}}  \frac{\eta^{\rm IM}(t(e))} {\eta^{\rm IM}(s(e))}  = \frac{1}{2n-4} $$
where $\eta^{\rm IM}(T)=1$ for all $T\in V_{n,A,\{ n \}}$ and $\lambda_{\rm IM}=d_{\{ n \}}=2n-4$. 
Thus, in the $\beta \to \infty$ limit, where the maximization reduces to the maximization of
the Shannon entropy of the probability distribution on edges of $\cG^{\rm IM}_{n,A,\{ n \}}$, achieved at
the uniform distribution. 
\endproof

\smallskip

\section{Parameterizing cube and dynamics} \label{CubeSec} 

In this and the following sections we mention some further directions in which it may be interesting to
explore the dynamical properties of the action of Merge on workspaces. One of the questions posed
in \cite{MCB} concerns the extent to which the structure formation process of Merge can be encoded
in the image of syntactic objects inside some semantic spaces, in the syntax-semantics interface.

\smallskip

The model of syntax-semantics interface developed in \cite{MCB} is based on the
idea that the computational structure of language resides entirely on the side of syntax
(syntax-first model), while the semantics side can be reduced to a minimal kind of
topological structure, capturing notions of proximity in some ambient ``semantic space"
(a topological or metric space, such as a manifold or stratified space). Since this
model of semantics is not, in itself, endowed with a computational structure-formation
process like Merge, it is natural to ask whether Merge would still be somehow
encoded in, and at least partially reconstructible from, the image of syntactic
objects in a semantic space. This question also has relevance to the current 
simulations of human language in machine architectures like large language models,
where one can ask whether there are circuits (in the sense of mechanistic
interpretability) involving the weights of attention modules of transformer architectures,
that would correspond to the computation of syntax via the Merge action on
workspaces. 

\smallskip

Another related question is about the realization of the process of structure formation
via Merge in neurocomputational architectures: the broad question of neuroscience
correlates of structure formation in syntax. A recent proposal \cite{MarBer} on possible 
ways of encoding the Merge operation in function spaces (such as wavelets) shows that
a good encoding of the Merge operation is possible in a deformation of the tropical 
semiring determined by the R\'enyi entropy. 

\smallskip

In both this setting and in the model of the syntax-semantics interface of \cite{MCB}, 
one considers a parameterizing space of non-planar full binary rooted trees with
labelled leaves, where continuous parameters are introduced, reflecting a metric
embedding for such trees inside some ambient semantic space. In the model of
the syntax-semantics interface of \cite{MCB} one uses the BHV moduli spaces \cite{BHV}
of non-planar trees to obtain this parameterization, following \cite{DeMo}. While
the BHV moduli spaces appear adequate to describe the embeddings of syntactic
objects in semantic spaces and the interaction with the planarization of trees in
externalization (via the relation to the associahedra and the moduli space of
real curves of genus zero as in \cite{DeMo}), this geometric setting does not seem
ideal to capture a possible image of the Merge operation. The
BHV moduli spaces have an interesting associated operad structure, the mosaic
operad discussed in \cite{Deva}. However, there is no immediate direct relation
between the boundary strata of the BHV moduli
spaces, which correspond to operations of shrinking edges to zero length, and the
coproduct operation that gives rise to the Merge action, which corresponds to
extracting accessible terms via admissible cuts. 

\smallskip

On the other hand, in \cite{MarBer}, continuous parameters are also
associated to syntactic objects, in the form of parameters $\lambda_v \in [0,1]$
associated to the non-leaf vertices of the tree, that determine a 
probability distribution at the leaves of the tree, as a polynomial function
in these variables, which encodes the tree structure and is optimized
over in the Merge operation represented as a sum operation in the
deformed semiring. Since in that setting the Merge operation becomes
expressed in terms of an optimization over continuous parameters,
this also raises the question of how the discrete Hopf algebra Markov chain
describing the action of Merge on workspaces can be mapped to a
setting where continuous parameterizations (metric trees) are involved. 

\smallskip

We do not address this question in the present paper, but we make
some general considerations in this section on continuous parameterizations 
of metric trees, as a way of illustrating the problem and some possible approaches. 

\smallskip

It was shown in \cite{MarBer} that, given a syntactic object $T\in \fT_{\cS\cO_0}$,
endowed with a head function $h_T$ as in Definition~\ref{headdef} (see 
also \cite{MCB}, \S 1.13.3), one can assign a metric structure to the edges
of $T$ by a set of parameters $\lambda=(\lambda_v)_{v\in V^o(T)}$, where each
variable $\lambda_v$ is valued in $[0,1]$, so that $\lambda \in [0,1]^{n-1}$,
with $n=\# L(T)$ and $V^o(T)$ the set of non-leaf vertices of $T$. The
edge lengths are assigned so that  the two edges below a given $v\in V^o(T)$
get lengths $\lambda_v$ and $1-\lambda_v$, as for example:
\begin{center}\textcolor{white}{[]}
\xymatrix@R=1.2cm@C=0.8cm{
  & *{} \ar@{-}[dl]_{\lambda_1} \ar@{-}[dr]^{1-\lambda_1} & \\
  \alpha & & *{} \ar@{-}[dl]_{\lambda_2} \ar@{-}[dr]^{1-\lambda_2} \\
  & \beta & & *{} \ar@{-}[dl]_{\lambda_3} \ar@{-}[dr]^{1-\lambda_3} \\
  & & \gamma & & \delta
}
\end{center}

\begin{lem}\label{headepsilon}
The datum of a head function is
equivalent to a choice of which of the two edges gets weight $\lambda_v$
and which $1-\lambda_v$. 
\end{lem}

\proof As recalled in Definition~\ref{headdef}, 
a head function $h_T: V^o(T)\to L(T)$ assigns a leaf to each
internal (non-leaf) vertex, with the rule that if $T_w \subset T_v$ and 
$h_T(v)\in L(T_w)$ then $h_T(v)=h_T(w)$.
Thus, if $h_T(v)$ is the leaf assigned to
$v\in V^o(T)$ by the head function $h_T$ and $\gamma_{h_T(v)}$ is
the path in $T$ from $v$ to the leaf $h_T(v)$, we can assign the
weights $\lambda_v$ and $1-\lambda_v$ according to the following rule:
to $e_v$ one of the two edges below a vertex $v\in V^o(T)$ we assign weight
\begin{equation}\label{epsilonvA}
\epsilon_{e_v} =\left\{ \begin{array}{ll} \lambda_v & e_v \in \gamma_{h_T(v)} \\ 1-\lambda_v & e_v \notin \gamma_{h_T(v)} 
\end{array}\right.
\end{equation}
Conversely, given an assignment of $\lambda_v$ and $1-\lambda_v$ to the pair of edges below $v$, for each 
$v\in V^o(T)$, the function $h_T: V^o(T)\to L(T)$ that assigns to $v$ the leaf obtained by following the path
from $v$ consisting of edges marked by the variable $\lambda_v$ (as opposed to $1-\lambda_v$) is a head function.
\endproof

\smallskip

The labeling of the edges by these continuous variables
induces a probability distribution $A^T=(a^T_\ell)_{\ell\in L(T)}$ on the leaves of $T$,
as shown in \cite{MCB}, where each $a^T_\ell$ is a polynomial function 
$a^T_\ell =P^T_\ell(\lambda)$ of the variables $\lambda_v$, of the form
\begin{equation}\label{Aprob}
a^T_\ell = \prod_{v\in \gamma_{v_0,\ell}} \epsilon_{e_v} \, ,
\end{equation}
where $\gamma_{v_0,\ell}$ is the path in $T$ from the root vertex $v_0$ to the leaf $\ell$. 

\smallskip

The probability $A$
is used in \cite{MarBer} to combine representations of lexical data $\alpha_\ell\in \cS\cO_0$
in an appropriate function space, to obtain a representation of the syntactic object $T$. 

\smallskip

We now consider the operation of restricting the point $\lambda=(\lambda_v)$ to
one of the faces in the the boundary of the cube $[0,1]^{n-1}$, namely restricting
one of the variables $\lambda_v$ to $0$ or $1$, as shown in the figure below.
\begin{center}
\includegraphics[scale=0.35]{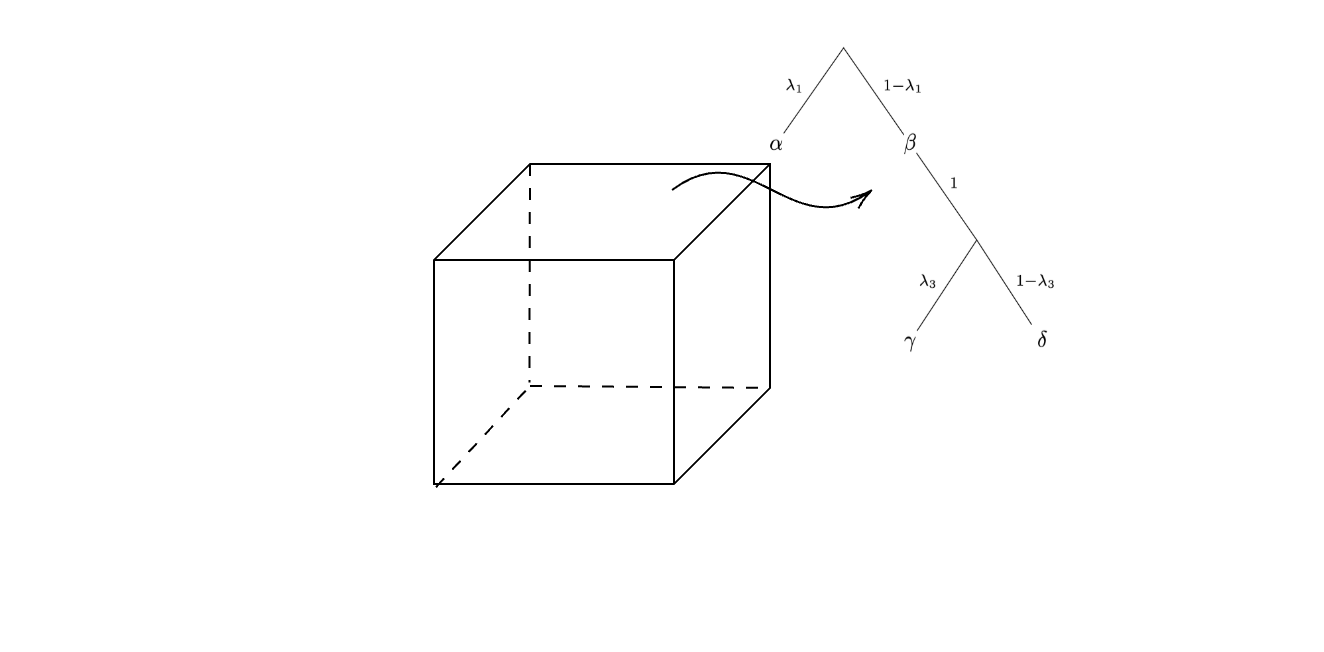}
\end{center}

\smallskip

\begin{lem}\label{delcontrA}
Taking one of the variables $\lambda_v$ to $0$ or $1$ corresponds to the
contraction of $T/^d T_w$, for $w$ the vertex below $v$ connected by
the edge $e_v$ with $\epsilon_{e_v}\to 0$.
\end{lem}

\proof  Consider the probability distribution $A^T=(a^T_\ell)_{\ell\in L(T)}$ and the
effect of taking one of the $\lambda_v$ to be $0$ or $1$. If $\lambda_v=0$, then
the edge $e_{v,1}$ with $e_{v,1} \in \gamma_{h_T(v)}$ has length $\lambda_v=0$
and the other edge $e_{v,2}$ below $v$ has length $1-\lambda_v=1$. Thus the probability
distribution $A^T$ changes in the following way
\begin{equation}\label{Acontr}
 a^T_\ell |_{\lambda_v=0} =\prod_{w\neq v \in \gamma_{v_0,\ell}} \epsilon_{e_w} \, , 
\end{equation} 
if the path $\gamma_{v_0,\ell}$ contains the edge $e_{v,2}$, and
\begin{equation}\label{Adel}
 a^T_\ell |_{\lambda_v=0} =0 \, , 
\end{equation} 
if the path $\gamma_{v_0,\ell}$ contains the edge $e_{v,1}$, and is unchanged,
\begin{equation}\label{Anone}
a^T_\ell |_{\lambda_v=0} = a_\ell \, ,
\end{equation} 
if $\gamma_{v_0,\ell}$ contains neither $e_{v,1}$ or $e_{v,2}$. 
Observe then that \eqref{Acontr}, together with \eqref{Anone} (and \eqref{Adel} since those
leaves are no longer present) gives exactly the resulting probability distribution
$A^{T/^d T_w}$, on the quotient tree $T/^d T_w$, for $w=t(e_{v,1})$ the vertex below $v$ along $e_{v,1}$, as
shown in the example
\endproof

This result extends from a single edge to an admissible cut as follows. 

\smallskip

\begin{cor}\label{delcontrAcut}
Let $C$ be an admissible cut on $T$ and let $F_{\underline{v}}=\pi_C(T)$ be the forest
extracted by the admissible cut, with $\rho^{(d)}_C(T)=T/^d F_{\underline{v}}$ the
corresponding quotient, in the coproduct term $\pi_C(T) \otimes \rho^{(d)}_C(T)$. Then
setting $\epsilon_{e}=0$ for each edge $e\in C$ of the cut gives
\begin{equation}\label{AcutC}
A^T |_{\epsilon_{e}=0\, , e\in C} = A^{T/^d  F_{\underline{v}}} \, . 
\end{equation}
\end{cor}

\proof As in Lemma~\ref{delcontrA} we have that $A^T |_{\epsilon_{e}=0\, , e\in C}$
is zero on all leaves in $F_{\underline{v}}$; is unchanged on all leaves where the
path to the root does not contain any vertex above an edge of the cut $C$; and,
for leaves where the path contains an edge sister to an edge of $C$, 
is the product of the $\epsilon_e$ of all the other edges of the path that are not sister to a cut edge. 
A path from a leaf to the root cannot contain more than 
one edge of $C$ because $C$ is an admissible cut. Thus the resulting 
$A^T |_{\epsilon_{e}=0\, , e\in C}$ is exactly $A^{T/^d  F_{\underline{v}}}$. 
\endproof

\smallskip

We can also obtain the probability distribution of the extracted accessible terms $T_v$
from $A^T$ by specialization to some of the deeper strata of the boundary of the
parameterizing cube $[0,1]^{n-1}$. 

\begin{lem}\label{lemATv}
Let $T_v \subset T$ be an accessible term of a syntactic object $T$, and let $\gamma_{v_0,v}$
be the path in $T$ from the root vertex $v_0$ to $v$. Then specialization of the variables $\epsilon_e$
for all edges $e\in \gamma_{v_0,v}$ to $\epsilon_e=1$ gives
\begin{equation}\label{ATv}
A^T|_{\epsilon_e=1\, , e\in \gamma_{v_0,v}} = A^{T_v} \, . 
\end{equation}
\end{lem}

\proof We have, for $\ell \in L(T_v)$
\begin{equation}\label{aellTv}
 a^T_\ell |_{\epsilon_e=1\, , e\in \gamma_{v_0,v}} = \prod_{w\in \gamma_{v,\ell}} \epsilon_w \, ,
\end{equation} 
while for $\ell \in L(T)\smallsetminus L(T_v)$, the path $\gamma_{v_0,\ell}$ will
contain a sister edge of one of the edges in $\gamma_{v_0,v}$. Since the weights $\epsilon_e$
of the edges in $\gamma_{v_0,v}$ are specialized to $\epsilon_e=1$, then their sister edges
are specialized to $0$. This implies that for all $\ell \in L(T)\smallsetminus L(T_v)$
\begin{equation}\label{aellTv2}
 a^T_\ell |_{\epsilon_e=1\, , e\in \gamma_{v_0,v}} = 0 \, .
\end{equation} 
This then gives exactly the probability distribution $A^{T_v}$. 
\endproof

\smallskip

We can then associate the following geometric locus in the parameterizing cube,
to the extraction of an accessible term $T_v\subset T$ and the corresponding quotient $T/^d T_v$.

\begin{prop}\label{TvquotCube}
Let $T_v \subset T$ be an accessible term. To the term $T_v \otimes T/^d T_v$ in the coproduct $\Delta^d(T)$
we can associate the pair $(A^{T_v}, A^{T/^dT_v})$ of probability distributions on $L(T_v)$ and
$L(T/^d T_v)=L(T)\smallsetminus L(T_v)$, respectively. When expressed in terms of the variables 
$\lambda=(\lambda_v)_{v\in V^o(T)}$, this pair of probabilities corresponds to a point 
$(\lambda^{T_v},\lambda^{T/^dT_v})$ in the product $\cQ_{T_v} \times \cQ_{T/^dT_v}$ of two 
sub-cubes of the parameterizing cube $[0,1]^{n-1}$, with dimensions $\dim \cQ_{T_v} + \dim \cQ_{T/^dT_v} = n-2$,
with $\cQ_{T_v} \simeq [0,1]^{V^o(T_v)}\subset [0,1]^{V^o(T)}$ and $ \cQ_{T/^dT_v} \simeq 
[0,1]^{V^o(T/^d T_v)}\subset [0,1]^{V^o(T)}$.
\end{prop}

\proof
In Lemma~\ref{delcontrA}, the probability $A^{T/^d  T_v}$ is a point in a stratum $\Delta_{L(T/^d  T_v)}$
of dimension $\# L(T/^d  T_v)-1$ inside the $(n-1)$-simplex $\Delta_{L(T)}$. On the other hand, the locus
inside the cube $[0,1]^{n-1}$ cut out by the condition $\lambda_v=0$ (or $\lambda_v=1$)
is a face $[0,1]^{n-2}$ of codimension $1$. The difference lies in the fact that the 
weights $\lambda_w$ of vertices $w\in V^0(T_v)$ do not contribute to $A^{T/^d  T_v}$,
because of \eqref{Adel}, while they are coordinates in $[0,1]^{n-2}$. Thus one can
obtain a sub-locus of the face $[0,1]^{n-2}$ that corresponds to points $A^{T/^d  T_v}$ by
fixing the value of the variables $\lambda_w$ with $w\in V^0(T_v)$ (because of \eqref{Adel},
it does not matter to which value: we will use $\lambda_w=0$). We denote by 
$\cQ_{T/^d  T_v}  \subset [0,1]^{n-1}$ for the resulting $(\# L(T/^d  T_v)-1)$-dimensional
sub-cube, which we can identify with $\cQ_{T/^d T_v} \simeq [0,1]^{V^o(T/^dT_v)}\subset [0,1]^{V^o(T)}$.

\smallskip

Similarly, in the case of Lemma~\ref{lemATv}, the probability $A^{T_v}$ corresponds 
to a point in a stratum $\Delta_{L(T_v)}$
of dimension $\# L(T_v)-1$ inside the $n-1$ dimensional simplex $\Delta_{L(T)}$.
On the other hand, the locus $\{ \epsilon_e=1\, , e\in \gamma_{v_0,v} \}$ 
determines a sub-cube $[ 0,1 ]^{n-r-1} \subset \partial [0,1]^{n-1}$ in the
boundary of the parameterizing cube, where $r={\rm length}( \gamma_{v_0,v})$,
which is of larger codimension, since in $A^{T_v}\in \Delta_{L(T_v)}$ the variables
$\epsilon_w$ of all vertices $w$ lying on paths from $L(T)\smallsetminus L(T_v)$ to the path $\gamma_{v_0,v}$
no longer contribute. Thus, we can restrict to a cube $[ 0,1 ]^{\#L(T_v)-1} \subset [ 0,1 ]^{n-r-1}$
obtained by setting all the $\lambda_w$ of these vertices to a fixed value. Again, it does not matter which one,
since these weights do not contribute to $A^{T_v}$ by \eqref{aellTv2}. For later convenience
we will set these values to $0$.  We denote by $\cQ_{T_v}\subset [0,1]^{n-1}$ 
the resulting $(\# L(T_v)-1)$-dimensional sub-cube, which we identify with 
$\cQ_{T_v} \simeq [0,1]^{V^o(T_v)}\subset [0,1]^{V^o(T)}$.

\smallskip  

Thus, for any accessible term $T_v \subset T$, we have obtained, by specialization of the variables 
in the parameterizing cube $[0,1]^{n-1}$, a pair of lower dimensional cubes $\cQ_{T_v}$ and $\cQ_{T/^d T_v}$ 
inside the boundary $\partial [0,1]^{n-1}$, with $\dim \cQ_{T_v} + \dim \cQ_{T/^dT_v} = n-2$.  
\endproof

\smallskip

\begin{cor}\label{Ajon}
Under the Internal Merge operation  $\fM(T_v, T/^d T_v)$, the probability distribution $A^{\fM(T_v, T/^d T_v)}$
at the leaves is a point in the join 
\begin{equation}\label{joinDelta}
 \Delta_{L(T_v)}\star \Delta_{L(T/^d T_v)} 
\end{equation}
of the simplexes $\Delta_{L(T_v)}$ and $\Delta_{L(T/^d T_v)}$, with
$A^{\fM(T_v, T/^d T_v)}$ a polynomial function as in \eqref{Aprob} of 
$\lambda\in [0,1]\times \cQ_{T_v} \times \cQ_{T/^dT_v}$.
\end{cor}

\proof
When we form the Internal Merge $\fM(T_v, T/^d T_v)$, a new root vertex $u$ is added, and correspondingly
a new variable $\lambda_u\in [0,1]$. The probability distribution at the leaves then satisfies
\begin{equation}\label{Amerge}
A^{\fM(T_v, T/^d T_v)} = \left\{ \begin{array}{ll} 
\lambda_u A^{T_v} + (1-\lambda_u) A^{T/^d T_v} & h_T(v_0)\in L(T_v)  \\ (1-\lambda_u) A^{T_v} + \lambda_u A^{T/^d T_v}  & h_T(v_0)\notin L(T_v)  \, ,
\end{array}\right.
\end{equation}
with $v_0$ the root of $T$. Namely $A^{\fM(T_v, T/^d T_v)}$ is a point in the join \eqref{joinDelta}, 
which can again be identified with the simplex $\Delta_{L(T)}$. In terms of the parameterizing cube, 
by Proposition~\ref{TvquotCube}, the probabilities $A^{\fM(T_v, T/^d T_v)} \in  \Delta_{L(T_v)}\star \Delta_{L(T/^d T_v)}$ 
can be written as polynomial functions $a^{\fM(T_v, T/^d T_v)}_\ell=\P_\ell(\lambda)$ 
of the form \eqref{Aprob}, for $\lambda\in [0,1]\times \cQ_{T_v} \times \cQ_{T/^dT_v}$.
\endproof

\smallskip

We have obtained above a description of the Internal Merge operation in terms
of a parameterizing cube $[0,1]^n=[0,1]\times [0,1]^{n-1}$. Starting with a
point $\lambda^T \in \{ 0 \} \times [0,1]^{n-1}$, with $A^T=P(\lambda^T)\in \Delta_{L(T)}$, 
one first assigns to it the pair of points $(\lambda^{T_v}, \lambda^{T/^d T_v})\in \cQ_{T_v} \times \cQ_{T/^dT_v}$,
with $(A^{T_v}=P(\lambda^{T_v}), A^{T/^d T_v}=P(\lambda^{T/^d T_v})))\in  \Delta_{L(T_v)}\times \Delta_{L(T/^d T_v)}$, 
by the specialization of variables as above. One then combines these to obtain a point as in \eqref{Amerge}
with $A^{\fM(T_v, T/^d T_v)}=P(\lambda_u, \lambda^{T_v}, \lambda^{T/^d T_v}) \in 
\Delta_{L(T_v)}\star \Delta_{L(T/^d T_v)}$, for $(\lambda_u, \lambda^{T_v}, \lambda^{T/^d T_v})\in [0,1]\times [0,1]^{V^o(T_v)} \times [0,1]^{V^o(T/^dT_v)} \subset [0,1]^n$. Thus, we can think of Internal Merge as the
result of a continuous flow inside $[0,1]^n$ that goes from a point in the sub-cube $\{ 0 \} \times [0,1]^{V^o(T)}$
to a point in the sub-cube $[0,1]\times [0,1]^{V^o(T_v)} \times [0,1]^{V^o(T/^dT_v)}$. 
One can then ask whether these Internal Merge transformations can be realized as part of 
a continuous flow on the parameterizing cube $[0,1]^n$. Similarly, one can look at the EM and
minimal SM transformations in terms of their effects on the probabilities $A^T=P(\lambda^T)$ and
the parameters $\lambda^T\in [0,1]^{V^o(T)}$, as a possible way to embed the Merge dynamics
into a realization of syntactic objects as metric trees inside some semantic space (as in Chapter~3 of \cite{MCB})
or in some function space representation as in \cite{MarBer}. Any such continuous flow would have to
satisfy some compatibility condition with the discrete Markov chain that we described in the previous sections
for the action of Merge on workspaces, in order to represent a realization of the Merge dynamics. 

\smallskip

A related question arises when one considers the syntactic objects, endowed with
an assignment of continuous parameters along the edges, as points in the BHV moduli
spaces of \cite{BHV}. The walls between the open cells of the BHV moduli space correspond
to the trees where one (or more in the higher codimensional strata) of the edge variables are
equal to $0$. A wall crossing between two adjacent open cells separated by a codimension one
boundary stratum corresponds to a transformation (at one of the
internal vertices of $T$) of the form
$$ \Tree[ [ $T_{v,1}$ $T_{v,2}$ ] $T_{v'}$ ] \mapsto \Tree[ $T_{v,1}$ [ $T_{v,2}$ $T_{v'}$ ]] \, . $$
Since wall crossings in BHV happen in this way, by moving an edge across one of the non-leaf vertices,
operations like IM that extract an accessible term, possibly from deep inside the tree, and merges
it back at the root, do not seem directly compatible with the form of the BHV boundary strata, in the
sense that the terms $T_v \otimes T/^d T_v$ do not seem to naturally correspond to codimension-one
strata (or $F_{\underline{v}} \otimes T/^d F_{\underline{v}}$ to higher codimensional strata). 

\smallskip

\begin{ques}\label{Q1}
Can the discrete Merge Hopf algebra Markov chain $\hat\cK^{(A,n)}$ (and its weighted version)
be embedded as a continuous Markov chain inside either the parameterizing cube $[0,1]^n$
or the moduli space BHV$_n$, or in topological models of semantic spaces? 
\end{ques}

We leave these considerations about continuous parameters, parameterizing spaces,
and flows compatible with the discrete Markov chain as questions for future work.

\section{I-language constraints and Colored Merge dynamics}\label{ColorSec}

The Hopf algebra Markov chain of the action of Merge on workspaces is the
dynamical model of {\em free} structure formation in syntax through the free 
symmetric Merge.

In the Minimalist Model of syntax, freely formed structures produced by the
free symmetric Merge are then filtered for compliance with theta role assignments
(theta criterion) and for well formed phases (head and complement structure,
extended projection). Both of these systems of filtering can be modeled
mathematically in terms of bud generating systems of colored operads, as
explained in \cite{ML} and \cite{MHL}. 

It is also shown in \cite{ML} and \cite{MHL} that filtering freely formed
structures by these coloring rules is equivalent to imposing constraints
on the Merge action and form only those structures that are compatible
with the constraints. The reason for this equivalence is that all the bud
generators of the relevant colored operads (see \S \ref{ColOpSec} for the
definitions) are corollas with a single
node, two input leaves, and one output root. One refers to this constrained
form of Merge as ``colored Merge".

When interpreted in this way, one sees that the effect of selecting only
those syntactic objects with well formed phases and theta roles corresponds
to pruning the free Merge dynamics given by the Hopf algebra Markov
chain $\hat\cK^{A,n}$ on the Merge graph $\cG_{n,A}$, so that the only
arrows remaining are those that have coloring corresponding to one
of the bud generators and the only vertices remaining are those workspaces
that are obtainable through these colored Merge operations, 
with the resulting coloring, so that all the component trees of the workspace
are elements of the given colored operad. 

In this colored Merge formulation, instead of the free Merge operations
described, as we have seen above, in the form
$$ \fM_{S,S'}=\sqcup \circ (\cB\otimes {\rm id}) \circ \delta_{S,S'} \circ \Delta \, , $$
one considers a version of the form
\begin{equation}\label{coloredM}
\fM^c_{S,S'} = \sqcup \circ (\cB^c \otimes {\rm id}) \circ\delta^c_{c_S,c_{S'}}  \delta_{S,S'} \circ \Delta \, ,
\end{equation}
where $c,c_S,c_{S'}\in \Omega$ are colors and
\begin{equation}\label{colorgen}
 \Tree[.$c$ [ $c_S$ $c_{S'}'$ ] ] 
\end{equation}  
belong to the set $\bB$ of bud generators of the colored operad $\fO_{\Omega,\bB}$ (see \S \ref{ColOpSec}
for more details), and $S,S'$ are 
syntactic objects that belong to the colored operad, with $c_S$ and $c_{S'}$
the colors assigned to the root vertex of $S$ and $S'$ respectively. The resulting
$\fM^c(S,S')$ has root vertex colored by the color $c$ and the remaining vertices
colored according to the coloring of $S$ and $S'$.

Let $\Omega$ be the set of colors ($\Omega_\Theta$,
$\Omega_\Phi$ for theta roles and phases, or $\Omega_{\Theta,\Phi}$ 
for the combined case as discussed in \cite{MHL}). 

The colored Merge graph $\cG_{n,A}^{\Omega}$ 
(in particular $\cG^{\Theta}_{n,A}$ or $\cG^{\Phi}_{n,A}$,
or $\cG^{\Theta, \Phi}_{n,A}$), can be obtained from $\cG_{n,A}$ in
two possible ways, reflecting the equivalent ways of thinking of
coloring as filtering of already formed structures or as building via
colored Merge. We describe these two constructions in \S \ref{ColMSec1}
and \S \ref{ColMSec2} after recalling briefly the notions of colored
operad and bud generating systems in \S \ref{ColOpSec}.

\subsection{Colored operads and bud generating systems}\label{ColOpSec}

A colored operad (in the category of sets) is a collection of sets 
 $$ \fO=\{ \fO(c,\underline{c})\,|\, \underline{c}=(c_1,\ldots, c_n), \,\, n\geq 1, \,\, c,c_i \in \Omega \} $$
 with $\Omega$ a (finite) set of colors, with composition operations
 \begin{equation}\label{colopergamma}
\begin{array}{rl}
\gamma: & \fO(c, (c_1,\ldots, c_n)) \times \fO(c_1, (c_{1,1}, \ldots, c_{1,k_1})) \times \cdots \times
\fO(c_n, (c_{n,1}, \ldots, c_{n,k_n}))  \\ & \to \fO(c, (c_{1,1}, \ldots, c_{1,k_1}, \ldots, c_{n,1}, \ldots, c_{n,k_n}))  \end{array}
\end{equation}
subject to an associativity property (which we do not write out explicitly here).
An algebra over a colored operad is a collection of sets $A=\{ A_c \}_{c\in \Omega}$
 with an action
 $$ \gamma_A: \fO(c, (c_1,\ldots,c_n)) \times A_{c_1}\times \cdots \times A_{c_n} \to A_c $$
 satisfying a compatibility condition with the composition $\gamma$ (which we also do not
 write explicitly here). 
 
 A convenient way of constructing colored operads, which is the basis for the
 formalization of theta roles and phase structure in \cite{ML} and \cite{MHL},
 is through {\em bud generating systems}, as in  \cite{Giraudo}.
 
 One starts with an ordinary operad $\fO$  (non-colored) and a finite set of colors $\Omega$.
 First one forms a colored operad by all possible assignments of colors,
 \begin{equation}\label{BudO}
\bB_\Omega(\fO)(n) := \Omega \times \fO(n) \times \Omega^n \, .
\end{equation}
Then one selects in $\bB_\Omega(\fO)(n)$ a finite set $\cR  \subset \bB_\Omega(\fO)$ 
of generators for the  actual (more restrictive) colored operad one wishes to obtain. 

One defines a {\em bud generating system} $\bB=(\fO,\Omega,\cR, \cI,\cT)$ as the data 
of the (non-colored) operad $\fO$,  the finite set of 
colors $\Omega$, separated into a set $\cI$ of initial colors and a set $\cT$ of terminal colors, 
and the choice of the set of generators $\cR  \subset \bB_\Omega(\fO)$ (the local coloring rules). 

One then obtains a colored operad $\fO_{\Omega, \bB}$ generated by the bud system 
$\bB=(\fO,\Omega,\cR, \cI,\cT)$ as 
\begin{equation}\label{opBO}
\fO_{\Omega,\bB}(n):= \{ x=(c,T,\underline{c})\in 
\bB_\Omega(\fO)(n) \,|\, {\bf 1}_c \to_{\bB} x, \, \, c\in \Omega \smallsetminus \cT, \, \, 
\underline{c}\in (\Omega \smallsetminus \cT)^n \} \, ,
\end{equation}
where ${\bf 1}_c \to_{\bB} x$ means that $x$ is obtained from the colored unit ${\bf 1}_c$ through
a sequence of operad compositions of generators $r_1,\ldots, r_N\in \cR$. 

The language $\cL(\bB)$ generated by the bud system $\bB=(\fO,\Omega,\cR, \cI,\cT)$ is an
algebra over the colored operad $\fO_{\Omega, \bB}$ given by 
\begin{equation}\label{LBlang}
 \cL(\bB)=\{ x=(c,T,\underline{c})\in \bB_\Omega(\fO)\,|\, {\bf 1}_c \to_{\bB} x, \, \, c\in \cI, \, \, \underline{c}\in \cT^n \} \, . 
\end{equation}

\subsection{Colored Merge graph via filtering}\label{ColMSec1}

For each 
$\alpha \in \cS\cO_0$ there will be a subset $\Omega_\alpha \subset\Omega$ of
colors that are compatible with $\alpha$. Thus, for $A=\{\alpha_1, \ldots, \alpha_n\}$
the lexical items at the leaves, we have corresponding choices of coloring
\begin{equation}\label{OmegaA}
\Omega_A:= \Omega_{\alpha_1}\times \cdots \times \Omega_{\alpha_n}.
\end{equation}

Given a choice of coloring $(c_1,\ldots, c_n)\in \Omega_A$
of the leaves $\alpha_1, \ldots, \alpha_n$ of $F$, for each vertex 
$F=T_1 \sqcup \cdots \sqcup T_r \in \fF_{A,n}$ of the graph $\cG_{n,A}$, and for each component
$T_i$ of $F$, consider all the possible ways in which $T_i$ can be colored as an element
of $\cL(\bB)(k_i)\subset \fO_{\Omega,\bB}(k_i)$, with $k_i=\# L(T_i)$. 

If at least one of the components 
does not have any possible coloring in $\cL(\bB)$ (for any available choice of
coloring of the leaves), then the vertex $F$ does not contribute any vertices 
to  $\cG_{n,A}^{\Omega}$. If there are colorings, then $F$ contributes a number of vertices of 
$\cG_{n,A}^{\Omega}$ ``above $F$", equal to the number of different possible colorings. 

For each of the coloring of an $F\in \fF_{A,n}$, and each Merge arrow (namely each $\cK^{A,n}_{F,F'}=1$
entry of the adjacency matrix of $\cG_{n,A}$), there is an arrow in $\cG_{n,A}^{\Omega}$ to one of the
colorings of $F'$, since we know that the colorings can be described via colored Merge
operations, and each coloring of $F$ is the target of an arrow from one of the coloring of some
$F''$ with $\cK^{A,n}_{F'',F}=1$. 

The graph consisting of these vertices and arrows is the Colored Merge graph $\cG_{n,A}^{\Omega}$.
Correspondingly this gives a Colored Merge Markov chain $\hat\cK^{A,n,\Omega}$, where 
$\cK^{A,n,\Omega}$ is the adjacency matrix of $\cG_{n,A}^{\Omega}$, and $\hat\cK^{A,n,\Omega}$ is
obtained as in \eqref{hatKmat}. 

Note that the relation between $\cG_{n,A}^{\Omega}$ and $\cG_{n,A}$ is not
straightforward, because there is no projection map between them, to which we could
apply our previously developed technique. 
There is rather a projection $p^\Omega: \cG_{n,A}^{\Omega} \to \cG^{\rm pruned}_{n,A}$ to a pruned version
of $\cG_{n,A}$, where one removes all the vertices $F$ that admit no coloring in $\cL(\bB)$ compatible with 
the available colorings of the leaves and of the vertices adjacent to them. 

The pruning that gives $\cG^{\rm pruned}_{n,A}$ will in general break the symmetry of $\cG_{n,A}$
with respect to the projection $p: \cG_{n,A} \to \cG_{\cP'(n)}$, making the dynamics of the
Colored Merge Markov chain $\hat\cK^{A,n,\Omega}$ more difficult to analyze than the
Hopf algebra Markov chain $\hat\cK^{A,n}$ of the free Merge action. Moreover, while the
behavior of $\hat\cK^{A,n}$ is independent of the specific choice of the lexical data
$A=\{ \alpha_1, \ldots, \alpha_n \}$ at the leaves, the dynamics of $\hat\cK^{A,n,\Omega}$
is sensitive to these data, and will be different for different choices of $A$, since the
set $\Omega_A \subset \Omega^n$ of available colorings at the leaves changes with $A$ 
and this set $\Omega_A$ determines which vertices are pruned in $\cG^{\rm pruned}_{n,A}$
and what the multiplicities (the sizes of the fibers) are in $p^\Omega: \cG_{n,A}^{\Omega} \to 
\cG^{\rm pruned}_{n,A}$.

This has a significant implication for the linguistic model. We know it is in principle equivalent 
to think of structure formation in syntax as being performed by the free symmetric Merge (the
dynamical system $\hat\cK^{A,n}$, followed by filtering by coloring for phases and theta roles
of the resulting fully formed structures $T \in p^{-1}(\{ n \})$, or else only form structures via
the Colored Merge operations, so that filtering by colorability happens at each step and the
dynamics is governed by $\hat\cK^{A,n,\Omega}$. These processes are equivalent in the
sense that they give rise to the same set of well-colored objects in $\cL(\bB)$ with
underlying trees $T \in p^{-1}(\{ n \})$. 

However, $\hat\cK^{A,n}$ and $\hat\cK^{A,n,\Omega}$ are not equivalent as dynamical systems. 
This difference can be viewed as an analog of the difference between having the same
weak generative capacity or the same strong generative capacity in the setting of formal languages.
Both $\hat\cK^{A,n,\Omega}$ or $\hat\cK^{A,n}$ followed by colorability-filtering of the $T \in p^{-1}(\{ n \})$
produce the same elements in $\cL(\bB)$, and this can be seen as analogous to having the same weak 
generative capacity for formal languages. However, the dynamical processes of structure
building involved are not in themselves equivalent.

From the dynamical systems point of view, it is more natural to think of structure formation as achieved by the
free Merge via $\hat\cK^{(A,n)}$ (and the weighted version) followed by filtering, so that filtering can be
understood as conditioning on the stationary distribution $\pi$ of the weighted free Merge dynamics.

\subsection{Colored Merge graph via dynamical building} \label{ColMSec2}

We also describe here how one can equivalently think of the construction of $\cG_{n,A}^{\Omega}$  as a
recursive process of 
replacing Merge arrows $\fM_{S,S'}$ with colored Merge arrows
$\fM^c_{S,S'}$.  One first needs to assume
that $S$ and $S'$ have been colored, and in particular one needs
to know the colors $c_S$ and $c_{S'}$ assigned to their root
vertices, and then this determines that the original arrow $\fM_{S,S'}$
is replaced by a collection arrows $\fM^c_{S,S'}$, one for each 
 $c$ with the property that \eqref{colorgen} is in the set of
 colored operad generators. When new arrows out of the
 vertices $\fM^c_{S,S'}$ are similarly analyzed, the SM and IM
 arrows will add new coloring to the original $S$ and $S'$,
 and so on until all allowed coloring and all allowed colored
 Merge operations are obtained.

One can start with the workspaces $F\in p^{-1}(\wp=\{ 2, \underline{1}^{n-2} \})$,
namely the workspaces with a single cherry component and all the other components
that are single leaves. Without loss of generality, say that $F=\fM(\alpha_1,\alpha_2) \sqcup
\alpha_3 \sqcup \cdots \sqcup \alpha_n$. 

For the single leaf components $\alpha_3 \sqcup \cdots \sqcup \alpha_n$
we have $\# (\Omega_{\alpha_3}\times \cdots \times \Omega_{\alpha_n})$
possible coloring, while for $\fM(\alpha_1,\alpha_2)$ the set of possible
coloring $$\Omega_{\alpha_1,\alpha_2}\subset \Omega\times \Omega_{\alpha_1}\times \Omega_{\alpha_2}$$
is the subset
$$ \Omega_{\alpha_1,\alpha_2}=\{ (c, c_1, c_2) \,|\,  \Tree[ .$c$ [ $c_1$ $c_2$ ] ] \in \bB \} \, . $$
Thus, over each vertex $F\in p^{-1}(\wp=\{ 2, \underline{1}^{n-2} \})$ in $\cG_{n,A}$ there will
be $N_{\alpha_1,\alpha_2} \cdot \prod_{i=3}^n N_{\alpha_n}$ vertices in $\cG_{n,A}^{\Omega}$,
where $N_{\alpha_i}=\# \Omega_{\alpha_i}$ and $N_{\alpha_i,\alpha_j} =\# \Omega_{\alpha_i,\alpha_j}$.

One can then consider the EM arrows from the vertices of $\cG_{n,A}^{\Omega}$ above the 
vertices in $p^{-1}(\wp=\{ 2, \underline{1}^{n-2} \})$ in $\cG_{n,A}$ 
and color the resulting component $\fM(\fM(\alpha_i.\alpha_j), \alpha_k)$ by assigning color $c'$
to the root vertex of $\fM(\fM(\alpha_i.\alpha_j), \alpha_k)$ if the colors $c$ and $c_k$ of
the root of $\fM(\alpha_i,\alpha_j)$ and of $\alpha_k$ satisfy
$$ \Tree[ .$c'$ [ $c$ $c_k$ ] ] \in \bB \, , $$
while no EM arrow exists if there is no $c'\in \Omega$ with this property. 

Iterating this process, for a given $\wp\in \cP'(n)$, one considers all the
workspaces $F\in p^{-1}(\wp)$ as vertices in $\cG_{n,A}$. Since each such $F$
is the target of some EM map, we obtain a collection of vertices in 
$\cG_{n,A}^{\Omega}$ above $F$, obtained by the coloring of targets of EM
maps as above, assuming that the source of these EM arrows have already
been colored by the same method. 

After obtaining all the possible coloring by colored EM, above each $F\in p^{-1}(\wp)$
for all $\wp\in \cP'(n)$, we consider the IM arrows at each $F\in p^{-1}(\wp)$
where $\wp$ has at least one $k_i\geq 3$ (so that there are IM arrows). For each
of the colored vertices above a given $F\in p^{-1}(\wp)$, we consider all the
possible colorings of IM arrows.

As discussed in \cite{ML} and \cite{MHL}, in
this case the coloring of the tree $\fM(T_v, T/^d T_v)$ is obtained as two
colored Merge operations: $\fM^c_{T_v,1}$ with coloring
\begin{equation}\label{MS1color}
 \Tree[ .$c_1$ [ $c_v$ $c_0$ ] ] \in \bB 
 \end{equation}
with $c_v$ the color of the vertex $v$ in $T$, $c_0$ the coloring assigned to 
the unit element $1$, and $c_1$ is another specified color for this kind
of generator (like the edge-of-phase color ${\mathfrak s}^\downarrow$ in the case
of phases or the non-theta position color for theta roles). This is then
followed by a colored Merge $\fM^{c'}(\fM^c_{T_v,1}, T/^d T_v)$, with
\begin{equation}\label{MS1col2}
 \Tree[ .$c'$ [ $c_1$ $c$ ] ] \in \bB 
\end{equation} 
for $c$ the color of the root of $T$. Since the leaf marked with $1$ is
not counted as a leaf in $\fM(T_v, T/^d T_v)$, this coloring is
equivalent to just the condition \eqref{MS1col2}
with $c$ the color of the root of $T$ and $c_1$ the color newly assigned to $v$
(as replacement of its original color $c_v$). 

Whenever a triple of colors satisfying the required condition 
does not exist, then the corresponding IM arrow of $\cG_{n,A}$ is pruned and 
does not exist as an edge of $\cG_{n,A}^{\Omega}$. These IM arrows
will add new colored vertices above $F\in p^{-1}(\wp)$, as they may produce
colorings not yet achieved by EM arrows. These vertices are added to 
$\cG_{n,A}^{\Omega}$.

Similarly, one can consider from all the obtained vertices of $\cG_{n,A}^{\Omega}$ 
all the possible 
minimal SM arrows and consider the corresponding coloring conditions.
As shown in \cite{MLH}, these will use additional generators in $\bB$ not
used by EM and IM, hence they will create new possible colorings, hence adding
new vertices to $\cG_{n,A}^{\Omega}$. 

One can then take again all the EM, IM, SM arrows from all the colored
structures obtained, and see if they give rise to new colorings. The process
stabilizes because the underlying tree structures are unchanged hence there
is a maximal number of possible colorings.  This produces a colored Merge
graph $\cG_{n,A}^{\Omega}$ from the Merge graph $\cG_{n,A}$.

\section{Parameter setting and Externalization} \label{ExtSec} 

There is another very important form of filtering on the syntactic objects formed
by the free symmetric Merge, in addition to the ``universal" (I-language)
filtering for phases and theta roles, namely the filtering by language-specific
syntactic parameters. The free symmetric Merge is seen as the core computational
structure of the Universal Grammar, while the parametric variation of individual
languages is modeled (since the original Principles and Parameters formulation)
through a procedure of parameter-setting, which happens during language acquisition,
and specializes the universal computational mechanism to the specifics of one
particular language. 

\smallskip

This process can be seen as further filtering the structures in $\cL(\bB)$ obtained
through the $\hat\cK^{A,n,\Omega}$ dynamics (or the $\hat\cK^{A,n}$ dynamics
followed by colorability-filtering) according to language-specific filters, as well
as a simultaneous assignment of a planar structure to the trees (also dependent
on language-specific word order parameters). This process is referred to
as Externalization.

\smallskip

Some discussion of Externalization in the context of the mathematical model
of Minimalism is included in \S 1.12 of \cite{MCB}, and a more extensive discussion
will appear elsewhere (\cite{MHB}, \cite{SM2}), so we will not discuss here a general model of
Externalization and syntactic parameters. We will, however, make some
comments on how the process of parameter setting can be looked at from
the perspective of the Markov chains $\hat\cK^{A,n}$
and $\hat\cK^{A,n,\Omega}$. 

\smallskip

We give only a very impressionistic view here of how the
process of parameter setting can be seen from this
dynamical system viewpoint. A full investigation will have to
incorporate a more detailed model of syntactic parameters
and the process of parameter setting, some of which is
being developed in \cite{MHB}, \cite{SM2}, so it will be
left to future work.

\smallskip

We work here within the model of Externalization proposed in \S 1.12
of \cite{MCB}, where we regard Externalization as a two-step process:
a section $\sigma_\cL$ that assigns planar structure to trees (in a language-dependent way)
and a projection $\Pi_\cL$ that filters out structures (also in a language dependent way),
thus resulting in a diagram of the form
$$ \xymatrix{ & \fT_{\cS\cO_0}^{pl} \ar[dr]^{\Pi_\cL} \ar[dl] & \\ \fT_{\cS\cO_0}\ar@/^2.0pc/[ur]^{\sigma_\cL} & & \fT^{pl}_{\cS\cO_0, \cL} } $$
with $\fT_{\cS\cO_0}^{pl}$ the set of planar binary rooted trees, and $\fT^{pl}_{\cS\cO_0, \cL}$ the set
of those planar binary rooted trees that are not filtered out by the syntactic parameters of the language $\cL$.

\smallskip

Behind this picture we assume that there is a space of syntactic parameters (which we take here to be 
binary valued for simplicity, though it is often useful to allow a third ``undefined" value). We
identify the space of all possible parameter values with the binary 
cube $\F_2^N =\{ 0,1 \}^N$, where $N$ is the total number of parameters. 
We assume here, also for simplicity of discussion, that these $N$ parameters are independent
binary variables, although in fact the problem of relations between syntactic parameters is one of
the crucial questions in the modeling of Externalization, but for our purposes here we can 
temporarily ignore it, and assume that an actual sufficient set of independent parameters can be identified.  
Among these $N$ parameters, a certain number $M<N$ deal with word order properties and 
have the effect of determining the section $\sigma_\cL$. The remaining parameters determine the projection 
$\Pi_\cL$. (For our purposes here it does not matter how exactly $\sigma_\cL$ and $\Pi_\cL$ are
obtained from the parameter values, except for a few considerations that we make below.) 

\smallskip

As part of the model, we assume that the lexical items $\alpha_\ell\in \cS\cO_0$ at the leaves
of the syntactic objects carry associated data $\varphi_\ell$ from a finite set $\bm{\Phi}$ of syntactic features. 
The syntactic parameters $\fp=\fp(\varphi \in \bm{\Phi})$, are functions of the features $\varphi \in \bm{\Phi}$, where
typically each syntactic parameter depends only on a certain subset of features. Thus, a syntactic
parameter $\fp$ affects (as planarization or as filtering) a certain structure $F^\Omega$ (a vertex 
of $\cG_{n,A}^{\Omega}$) iff those features $\varphi \in \bm{\Phi}$ that $\fp$ depends on
are present at the leaves of the components of $F^\Omega$.  We write $\bm{\Phi}_\fp \subset \bm{\Phi}$
for the subset of features that the parameter $\fp$ depends on.

\smallskip

Parameter values are set in the process of language acquisition. As noted in \cite{LongoTre}, these
two kinds of parameters behave differently in that respect. For syntactic parameters of the type 
that we associate to determining the projection $\Pi_\cL$, there is empirical evidence for a default
value of $-$ (or $0\in \F_2$), and positive evidence in language acquisition is required to flip the
value to $+$ ($1\in \F_2$). Thus, we can assume that the subspace $\F_2^{N-M}$ is initialized at
the $\underline{0}$ vector. The case of the word-order parameters, on the other hand, does not appear 
to have a default state. So it is suggested in \cite{LongoTre} that these parameters should have an
``unspecified" initial state, out of which the two possible $-/+$ states emerge as a symmetry
breaking phenomenon. In the model we present here, constructing the
section $\sigma_\cL$ can be framed as a setting where word-order parameters are initialized in 
an undefined state and positive evidence sets them (equivalently, recursively constructs $\sigma_\cL$).
On the other hand, the analysis in \cite{MHB} of the Final-Over-Final-Condition 
suggests that one may assume a default harmonic head-final planar embedding determined by the 
head function, with respect to which other embeddings are obtained via a group action, with 
the FOFC representing a minimality property with respect to generators of that group action. Note that
all the syntactic objects that are colorable for phases have a head function, so the default assumption
on the embedding is possible after filtering by well-structured phases. 
For simplicity, we just assume here that there is no default value for word order parameters (hence no 
default choice of planar embedding). 

\smallskip

Assume that we start with the graph $\cG_{n,A}^{\Omega}$, so that the filtering by theta roles and
phases have already been taken care of. 

\smallskip

As a first step, we need to perform a language-specific choice of planar embeddings (the section $\sigma_\cL$)
on the component trees $T_i^\Omega$ in the workspaces $F^\Omega$ (the vertices of $\cG_{n,A}^{\Omega}$, constructed as
described in the previous section). The $\Omega$ superscript is used here to remember that these trees
are colored, $T_i^\Omega\in \cL(\bB)$. 

\smallskip

To that purpose, we can take, above each vertex $F^\Omega=T_1^\Omega\sqcup \cdots \sqcup T_r^\Omega$ 
of $\cG_{n,A}^{\Omega}$, a collection of 
$$ 2^{n-r}=\prod_{i=1}^r 2^{k_i-1} $$
vertices, given by all the possible choice of planar embeddings for each of the components $T_i^\Omega$.
(A full binary rooted tree $T$ with labelled leaves has $2^{\# L(T)-1}$ planar embeddings.)

\smallskip

A priori, none of these vertices is selected. However, exposure to positive examples in language acquisition
will select some $\sigma_\cL(T^\Omega)$, namely selects a particular planar embedding for some
syntactic objects $T^\Omega\in \cG_{n,A}^{\Omega}$. This selection automatically causes further selection of the same 
planar embedding wherever $\sigma_\cL(T^\Omega)$ occurs: as a component of some $F^\Omega$
in some other $\cG_{m,A'}^{\Omega}$ with $A'$ another set of lexical data, sharing the same features
that the word-order parameters depend on as $A$. It also selects $\sigma_\cL(T^\Omega)$ when
$T^\Omega$ occurs as an accessible term inside some larger $\tilde T^\Omega$, unless some
other word-order parameter causes a change of planar embedding of $T^\Omega$ caused by
other data present in the rest of $\tilde T^\Omega$.  A locality and finiteness assumption on
syntactic parameters (see \cite{MHB}, \cite{SM2} for further discussion) implies that all 
$\sigma_\cL(T^\Omega)$ of arbitrary size, are determined by knowing the embeddings $\sigma_\cL$
on the graphs $\cG_{n,A}^\Omega$ with $n\leq n_0$ small, by identifying the planarization of 
accessible terms and quotient terms of small size, as discussed below. (See \cite{MHB}, \cite{SM2} for specific bounds 
on the size $n_0$ needed.)

\smallskip

This process selects a set of vertices $\sigma_{\cL}(F^\Omega)$, which has the same
cardinality as the set of vertices $F^\Omega$ of $\cG_{n,A}^{\Omega}$, since $\sigma_\cL$
assigns a unique planar structure to each $F^\Omega$ (we are not discussing here scrambling phenomena and
how to account for languages with a relatively free word order). We can equivalently write
these vertices as pairs $(F^\Omega, \pi_\cL)$ where $\pi_\cL$ is a choice of one among the 
$2^{n-r}$ planar embeddings.

\smallskip

The remaining $N-M$ parameters that determine the projection $\Pi_\cL$ are initialized as 
$\underline{0}\in \F_2^{N-M}$. Unless positive examples during language acquisition flip
the values of some of these parameters to $1$, they will remain in the $0$ default value. 
A syntactic parameter $\fp$ having
default value $0$ means that all structures $F^\Omega\in \cG_{n,A}^{\Omega}$ where the 
syntactic features $\varphi_\ell$ associated to lexical data $\alpha_\ell$ at the leaves $\ell \in L(T)$ 
are involved in determining the parameter $\fp$ are filtered out in the projection $\Pi_\cL$.

\smallskip

We assume two basic properties of the projection $\Pi_{\cL}$.
\begin{itemize}
\item The projection $\Pi_{\cL}$ as decomposes into a product
$\Pi_{\cL}=\prod_{\fp} \Pi_{\cL, \fp}$. Thus, a structure $F^\Omega$
passes the filter iff it passes it with respect to each parameter, 
$\Pi_{\cL, \fp}(F^\Omega)=F^\Omega$ for all $\fp$. 
\item For $F^\Omega=T^\Omega_1 \sqcup \cdots \sqcup T^\Omega_r$, we have
$\Pi_{\cL}(F^\Omega)=0$ iff $\Pi_{\cL}(T^\Omega_i)=0$ for some $1\leq i \leq r$. 
\end{itemize}

\smallskip

Initializing all the $N-M$ parameters to $0$ means that at the start of the process {\em all}
syntactic objects are filtered out, so nothing passes the Externalization filter. Suppose then
that positive examples are encountered in the language acquisition process that cause
one of the parameters $\fp$ to switch to value $1$. This change will affect only those 
structures $F^\Omega\in \cG_{n,A}^{\Omega}$ where the features at the leaves satisfy 
$\{ \varphi_\ell \}_{\ell\in L(F^\Omega)}\supset \bm{\Phi}_\fp$. 

\smallskip

By the considerations above on the projection $\Pi_\cL$, we can assume that
$F^\Omega=T^\Omega$ has a single component, and that the set of
features $\{ \varphi_\ell \}_{\ell\in L(T^\Omega)}$ contains the set $\bm{\Phi}_\fp$.

\smallskip

To model the way in which changing the value of $\fp$ to $1$
changes $\Pi_{\cL,\fp}$, we can first assume that positive examples suffice to determine
which among the tree $T^\Omega$ with sufficiently small size $\# L(T^\Omega)\leq n_0$
and with $\{ \varphi_\ell \}_{\ell\in L(T^\Omega)}\supset \bm{\Phi}_\fp$ 
are now in the range of $\Pi_{\cL,\fp}$ and no longer in the kernel. 
In other words, the set 
\begin{equation}\label{Tfp}
\cT_\fp :=\{ T^\Omega \,|\, \# L(T^\Omega)\leq n_0 \text{ and } \{ \varphi_\ell \}_{\ell\in L(F^\Omega)}\supset \bm{\Phi}_\fp \} 
\end{equation}
is separated out into two subsets
\begin{equation}\label{pmTfp}
\cT_\fp = \cT_\fp^- \sqcup \cT_\fp^+
\end{equation}
where $\Pi_{\cL,\fp}|_{\cT_\fp^+}={\rm id}$ and $\Pi_{\cL,\fp}|_{\cT_\fp^-}= 0$. 

\smallskip

For trees $T^\Omega$ of larger size, a locality condition for the parameter $\fp$ (see
\cite{MHB}) would ensure that, if all the accessible terms $T^\Omega_v \subset T^\Omega$
and respective quotients $T^\Omega /^d T^\Omega_v$ of size at most $n_0$, 
whose leaves contain the relevant features, are in the range of the projection $\Pi_{\cL,\fp}$, 
then so is $T^\Omega$. (Note that $T^\Omega$ could still be in the kernel of other projections 
$\Pi_{\cL,\fp'}$ for other parameters $\fp'\neq \fp$ hence in the kernel of $\Pi_\cL$.) 

\smallskip

With this procedure of parameter setting, small positive examples suffice to ``activate"
structures $F^\Omega$ of arbitrary size through the locality condition. If all the accessible
terms and quotients of components of $F^\Omega$ that are of size at most $n_0$ are
in $\cT_\fp^+$ then $F^\Omega \in \cT_\fp^+$, so that $\Pi_{\cL,\fp}(F^\Omega)=F^\Omega$
and $\Pi_{\cL,\fp}(F^\Omega)=0$ otherwise. 

\smallskip

When all the parameters have been set in the language learning process,
to a new vector $\underline{\fp}\in \F_2^N$, the process describes above
results in activating a certain subset of the vertices $(F^\Omega, \pi_\cL)$
of the graphs  $\cG_{n,A}^{\Omega}$ for any choice of $n$ and $A$. 
We denote this set of vertices as $V(\cG_{n,A}^{\Omega,\cL})$. The
graph $\cG_{n,A}^{\Omega,\cL}$ is then the graph induced by 
$\cG_{n,A}^{\Omega}$ on the set of vertices $V(\cG_{n,A}^{\Omega,\cL})$.
The selection of $\cT_\fp^+$  and $\cT_\fp^-$ and the resulting construction
of $\cG_{n,A}^{\Omega,\cL}$ completes the process of parameter setting. 

\smallskip

For a given $\cL$, the collection of the graphs $\cG_{n,A}^{\Omega,\cL}$, each 
with their associated MERW Markov chain, can be viewed as a model of the 
generative grammar of the language $\cL$.

\smallskip

\subsection*{Acknowledgments} The first author was supported for this project by NSF grants DMS-2104330 and 
DMS-2506176 and by Caltech’s T\&C Chen Center for Systems Neuroscience. The second
author was supported for this project by an Epistea fellowship and hosted by Caltech's Visiting 
Undergraduate Research Program. We thank Paolo Aluffi, Yassine El Maazouz,  Sita Gakkhar, Henry Gustafson, Riny Huijbregts, 
Amy Pang, Juan Pablo Vigneaux, and Elizabeth Xiao for helpful conversations.

\end{document}